\documentclass [12pt]{article}      
\usepackage{amsmath} 
\usepackage{amsthm}
\usepackage{latexsym}
\usepackage{pb-diagram}
\usepackage{amssymb}
\usepackage{bbm}
\usepackage{color}
\usepackage{geometry}
\usepackage{graphicx}
\usepackage{hyperref}
\usepackage{cancel}
\usepackage{pdfpages}
\hypersetup{
    colorlinks,
    citecolor=blue,
    filecolor=blue,
    linkcolor=blue,
    urlcolor=blue
}
\hypersetup{linktocpage}
\usepackage{enumerate}
\usepackage{todonotes}
\makeatletter
\let\r@eqnnum\@eqnnum
\input{leqno.clo}
\let\l@eqnnum\@eqnnum
\newcommand{\leqnos}{\let\@eqnnum\l@eqnnum}
\newcommand{\reqnos}{\let\@eqnnum\r@eqnnum}
\reqnos
\makeatother

\newtheorem{theorem}{Theorem}
\newtheorem*{theorem*}{Theorem}
\newtheorem{corollary}[theorem]{Corollary}
\newtheorem*{corollary*}{Corollary}

\newtheorem{definition}[theorem]{Definition}
\newtheorem*{definition*}{Definition}
\newtheorem{prop}[theorem]{Proposition}
\newtheorem*{prop*}{Proposition}
\newtheorem{proposition}[theorem]{Proposition}

\newtheorem{example}[theorem]{Example}
\newtheorem*{example*}{Example}
\newtheorem{fact}[theorem]{Fact}
\newtheorem*{fact*}{Fact}
\newtheorem{remark}[theorem]{Remark}
\newtheorem*{remark*}{Remark}

\newtheorem*{claim*}{Claim}

\newtheorem{open}{Open Problem}

\newcommand{\rest}{\ensuremath{\upharpoonright}}

\newcommand{\la}{\langle}
\newcommand{\ra}{\rangle}

\newcommand{\Chi}{\ensuremath{\mbox{\Large{$\chi$}}}}
\newcommand{\bsigma}{\mbox{\boldmath$\Sigma$}}
\newcommand{\bpi}{\mbox{\boldmath$\Pi$}}
\newcommand{\bdelta}{\mbox{\boldmath$\Delta$}}

\newcommand{\trees}{\ensuremath{\mathcal{T}\!\! rees}}

\newcommand{\xbmt}{\ensuremath{(X,\mcb,\mu,T)}}
\newcommand{\ycns}{\ensuremath{(Y,\mcc,\nu,S)}}

 \newcommand{\poN}{\mathbb N}
\newcommand{\poZ}{\mathbb Z}
\newcommand{\nn}{{\mathbb N}}

\newcommand{\bk}{{\mathbb K}}

\newcommand{\bt}{\mathbb T}

\newcommand{\poR}{\ensuremath{\mathbb R}}

\newcommand{\mct}{{\ensuremath{\mathcal T}}}

\newcommand{\mca}{\ensuremath{\mathcal A}}
\newcommand{\mcg}{\ensuremath{\mathcal G}}

\newcommand{\mcx}{{\mathcal X}}
\newcommand{\mcc}{{\mathcal C}}
\newcommand{\mck}{{\mathcal K}}
\newcommand{\mcb}{\mathcal B}

\newcommand{\bfni}[1]{\noindent {{\bf{#1}}}}

\renewcommand{\qed}{{\nopagebreak \hfill $\dashv$ 
 \par\bigskip}}

\newcommand{\pf}{{\par\noindent{$\vdash$\ \ \ }}}

\newcommand{\oto}{one-to-one\ }

 \title{The complexity of the Structure and Classification of Dynamical Systems}
 \author{Matthew Foreman}
 \date{}
 
\begin{document}
 \maketitle

\bfni{\underline{Glossary}}

\begin{description}
 \item[Analytic set] A subset $A$ of a Polish space $X$ is  \hyperlink{aandcoa}{\emph{analytic}} if there is a Polish space $Y$ and a Borel set $B\subseteq X\times Y$ such that $A=\{x:$ for some $y\in Y, (x,y)\in B\}$.  Equivalently $A$ is the projection of $B$ to the $X$-axis.
 \item[Anosov diffeomorphism] Definition \ref{anosovdiff}.
 \item[Benchmarking] An informal term describing the location of a set or an equivalence relation in terms of reducibility to other sets or equivalence relations. 
 \item[Bernoulli Shift] Definition \ref{Bshift}.
 \item[Borel hierarchy] Definition \ref{borel hi}
 \item[Co-analytic set] A subset $C$ of a Polish space $X$ is \hyperlink{aandcoa}{\emph{co-analytic}} if there is a Polish space $Y$ and a Borel set $B\subseteq X\times Y$ such that $C=\{x$ for all $y\in Y, (x, y)\notin B$. Equivalently, the complement of $C$ is analytic.
 \item[Descriptive complexity] This term is used for the placement of a structure or a classification problem among the benchmarks of reducibility.
 \item[Distal] \emph{Topologically distal} is definition \ref{topodistal}, \emph{measure distal} is definition \ref{measdist}.
 \item[$=^+$ equivalence relation] Introduced and discussed in section \ref{fsbs}.
 \item[$E_0$ equivalence relation] Definition \ref{E0}.
 \item[Factor map: measurable] Let $X, Y$ be standard measure spaces and $T:X\to X$, $S:Y\to Y$ be measure preserving transformations.  Then $S$ is a \emph{factor} of $T$ if there is a (not necessarily invertible) measure preserving $\pi:X\to Y$ such that $\pi\circ T=S\circ \pi$ almost everywhere. The map $\pi$ is the \emph{factor map}.
 \item[Factor map: topological] Let $X, Y$ be topological spaces and $T:X\to X$,\\ $S:Y\to Y$ be homeomorphisms. Then $S$ is a \emph{factor} of $T$ if there is a continuous surjective map $\pi:X\to Y$ such that $\pi\circ T=S\circ \pi$.  The map $\pi$ is the \emph{factor map}.
  \item[Ill-founded, Well-founded] If $T\subseteq X^{<\nn}$ is a tree, then $T$ is \emph{well-founded} if and only if there is no function $f:\nn\to X$ such that for all $n, f\rest\{0, 1, 2, \dots n\}\in T$.  Equivalently, $T$ has no infinite paths. If $T$ is not well-founded, then $T$ is \emph{ill-founded}.
 \item[$\mck$-automorphisms] Definition \ref{Kaut}.
 \item[Kakutani Equivalence] Definition \ref{Kakeqdef}.
 \item[Minimal homeomorphism] A homeomorphism $h:X\to X$ is minimal if for every $x\in X$, $\{h^n(x):n\in\nn\}$ is dense.
 \item[Measure Conjugacy] Let $S, T$ be measure preserving transformations defined on standard measure spaces $(X, \mcb, \mu)$ and 
 $(Y,\mcc,\nu)$. Then $S, T$ are \emph{measure conjugate} if there is an invertible measure preserving transformations $\phi:X\to Y$ such that
 $\phi\circ T= S\circ \phi$ almost everywhere.
 \item[Morse-Smale diffeomorphism]  Definition \ref{morsesmaledef}.
 \item[$\bpi^1_1$-norms] Definition \ref{pi11normsdef}.
 \item[$\boldsymbol{\phi^*(\mu)}$] Let $X, Y$ be Polish spaces, $\mu$  a measure on $X$ and $\phi:X\to Y$ be a measurable map. Then $\phi$ and $\mu$ induce a measure $\phi^*(\mu)$ on the Borel subsets of  $Y$ by setting $\phi^*(\mu)(B)=\mu(\phi^{-1}(B)$.
 \item[Polish Space] A Polish space is a topological space $(X,\tau)$ such that there is some complete separable metric $d$ on $X$ inducing the topology $\tau$. The topology $\tau$ is called a \emph{Polish Topology}.
 \item[Polish Group Action] A \emph{Polish group} is a topological  group with a Polish topology.  A Polish group action is a Polish group acting $G$ acting jointly continuously on a Polish space $X$. See section \ref{pgas}.
 \item[Reduction] If $A$ and $B$ are subsets of Polish spaces $X$ and $Y$, then a one dimensional  \emph{reduction} of $A$ to $B$ is a function $f:X\to Y$ such that for all $x$, $x\in A$ if and only if $f(x)\in B$. 
 
 If  $A\subseteq X\times X$ and $B\subseteq Y\times Y$, then a two dimensional reduction is a function 
 $f:X\to Y$ such that for all $(x_1, x_2)\in X\times X, (x_1, x_2)\in A$ if and only if $(f(x_1), f(x_2))\in B$.
 The reduction is \emph{continuous} or {Borel} if the function is continuous of Borel.  One and two dimensional Borel reductions are notated as $\preceq^1_\mcb$ and $\preceq^2_\mcb$. (Section \ref{reductions}.)
 \item[Rotation Number] Definition \ref{rotationnumber}
 \item[Separable and complete measure spaces] A measure space $(X,\mcb,\mu)$ is \emph{separable} if there is a countable subset $\mca\subseteq\mcb$ such that for all $S\in \mcb$ and $\epsilon>0$, there is an $A\in \mca$ such that $\mu(A\Delta S)<\epsilon$. A measure space is \emph{complete} if whenever $A\in \mcb$ has $\mu(A)=0$ and $B\subseteq A$, then $B\in \mcb$.
 \item[Smooth equivalence relation] Definition \ref{smoothdef}.
 \item[Smooth manifold] A manifold with a $C^k$ structure for some $k\ge 1$.
 \item[Smooth transformation] A $C^k$ map $f:M\to M$ where $M$ has a $C^k$-structure  for some $k\ge 1$.
 \item[Standard measure space] A measure space $(X,\mcb, \mu)$ is \emph{standard} if it is separable and complete.
 \item[The group $S_\infty$] The group $S_\infty$ is the group of permutations of the natural numbers. (Section \ref{S-inf}.)
 \item[Topological Conjugacy]  Let $X, Y$ be topological spaces and $f:X\to X, g:Y\to Y$. Then $f$ and $g$ are topologically conjugate if there is a homeomorphism $\phi:X\to Y$ such that $\phi\circ f=g\circ \phi$.
 \item[Trees] A (downward branching) tree is a partial ordering $(T, \le_T)$ such that for all $t\in T$, there is a well-ordering $\le^*$ of  $\{s\in T:s\ge t\}$  such for $s_1, s_2\ge t$, 
 \[s_2\le_T s_1 \mbox{ if and only if } s_1\le^* s_2.\]
 A \emph{path} through a tree is a one-to-one function $p$ from an ordinal $\alpha$ whose range is $\le_t$-upwards closed and if 
 $\beta<\gamma<\alpha$, $p(\gamma)\le_Tp(\beta)$. A \emph{branch} through $T$ is a maximal path through $T$. A tree is \emph{well-founded} if and only if it has no infinite paths. 
 \item[Trees of finite sequences] Fix a set $X$. Let $X^{<\nn}$ be the collection of finite sequences of elements 
 of $X$.  Order $X^{<\nn}$ by setting $\sigma\le \tau$ if and only if $\sigma$ is an initial segment of $\tau$.  A 
 set $T\subseteq X^{<\nn}$ is a \emph{tree of finite sequences of elements of $X$} if and only if whenever 
 $\sigma\in T$ and $\tau$ is an initial segment of $\sigma$ then $\tau\in T$. Trees of finite sequences can be viewed as elements of $\{0, 1\}^{X^{<\nn}}$.  Putting the discrete topology on $\{0,1\}$ and the product topology on $\{0, 1\}^{X^{<\nn}}$, the space of trees is a compact topological space.  If $X$ is countable then the space of trees on $X$ is homeomorphic to the Cantor set. 
  
 Trees are discussed in section \ref{reducingillfounded}.   %
  \item[Topologically transitive] If $X$ is a metric space, a homeomorphism  $h:X\to X$ is topologically transitive if for some $x\in X$, $\{h^n(x):n\in\nn\}$ is dense. 
   \item[Turbulent] Turbulence is a property of a continuous action of a Polish group on a Polish space. It gives a method for showing that an equivalence relation is not the result of an $S_\infty$-action. It is discussed in section \ref{turbulence} where the notion is formally defined.
  \item[Well-founded, Ill-founded] If $T\subseteq X^{<\nn}$ is a tree, then $T$ is \emph{well-founded} if and only if there is no function $f:\nn\to X$ such that for all $n, f\rest\{0, 1, 2, \dots n\}\in T$.  Equivalently, $T$ has no infinite paths. If $T$ is not well-founded, then $T$ is \emph{ill-founded}.

 \end{description} 

\tableofcontents
 \section{Introduction: What is a dynamical system? What is structure? What is a classification?}\label{whatisdyn}
 This article is a survey of the \emph{complexity} of structure, anti-structure,  classification and anti-classification results in dynamical systems. The paper focusses primarily on ergodic theory, with excursions into topological dynamical systems. Part of the intention is to  suggest methods and problems in related areas. It is not meant to be a comprehensive study in any sense and makes no attempt to place results in a historical context.  Rather it is intended to suggest that classification and anti-classification results that are likely to be common in many areas of dynamical systems.  Hence the focus is on giving particular examples where descriptive set theory techniques have been successful, and mention related work in passing to suggest that the examples given here are not isolated.\footnote{Apologies for leaving out important results, this reflects my ignorance and space limitations.}
 
 The target audience for this paper is two fold.  In no particular order, one audience consists of researchers whose primary interest is in some area of dynamics who may not be aware of, or who want to learn more about the methods for studying the complexity of structure and classification offered by descriptive set theory. For these people, the attempt is to give familiar examples from dynamical systems where the complexity is understood relative to known benchmarks.  Logical technicalities are avoided as much as possible, and the background necessary is given in the appendix.  This is \emph{naive} descriptive set theory--with the rough meaning that it omits the technicalities of quantifiers. A self contained source on basic descriptive set theory, requiring very minimal background is 
 \cite{naivedst}. There are many excellent, more advanced sources for this material (\cite{kechris}, \cite{moscho}, \cite{marker}). The author makes no attempt at historical attributions of the theorems in descriptive set theory proper, these are well-covered in \cite{kechris}, \cite{moscho}.

 The second audience are those descriptive set theorists who are likely to be very familiar with the language of \emph{reductions} and \emph{quantifier complexity}.  For these people, the attempt is to give a very high level overview of some of the known results specifically in dynamical systems.    The choices of examples reflect the authors taste and background, but include extremely famous results and the problems left open by them. There is a list of open problems at the end of the paper and more will appear in  a future arXiv submission with several co-authors, particular F. G. Ramos.
 
 The upshot is that in any given paragraph in the paper one of the audiences will be questioning the point of what is being presented.  I ask the reader's patience in understanding the dual mission.
  \medskip
 
 \paragraph{What is the difference between a structure theory and a classification?}
 The setting for all of these results is a Polish space $X$, and we study classes $\mcc\subseteq X$.  In some cases $\mcc=X$, in which case the structure question is less cogent. 
 
 The distinction between structure and classification is a bit vague. Very roughly a structure theorem for 
 $\mcc$ gives criteria for membership of elements that belong to $X$ to be in $\mcc$. The structure theorem gives more information than the definition of $\mcc$.  
 
 Often this additional information gives simply gives a very explicit test for belonging to $\mcc$ that may use ideas slightly different than the definition. Another form of a structure theorem is 
  a  method for building an element of $\mcc$.  The form of this type of classification is:
	\begin{center}
	$x\in \mcc$ \\
	if and only if\\
	$x$ can be built using a construction with properties 1)-n)	
	\end{center} 
 where properties 1-n) are concrete and explicit.

While a structure theory involves one object at a time and answers the question whether that object is in 
  $\mcc\subseteq X$, a classification theory involves two objects known to be in $\mcc$ and asks whether they are equivalent with respect to an equivalence relation in question. It is usually concerned with assigning invariants to the equivalence classes.
 
 Examples of common equivalence relations can include being in the same orbit by a group action, or being isomorphic, or have some other property in common.
  \medskip

 \medskip
 
 This paper is concerned with two cogent questions.  The first is 
 
	 \begin{quotation}
 	\noindent Can the structure or classification theory be 
 	done with inherently countable information (is it Borel)? 
 	\end{quotation}

 The second question is:
 	\begin{quotation}
 	\noindent Where does its complexity sit relative to existing benchmarks? For example: Are there numerical invariants?
	Can one assign countable structures to elements of the class so that isomorphism is the invariant?
	\end{quotation}

Many classical structure theorems were proved in the 1960's and 1970's.  They fit well into the rubric suggested by this paper.

\paragraph{What is anti-classification?} Reductions are tools that allow the establishment of  lower bounds on 
the complexity of equivalence relations. These lower bounds are often the established benchmarks discussed in this paper.  The complexity of some of some benchmarks is extremely high--often they are  not even Borel.  

Why is this important? An underlying thesis of this paper is that a general solution to a question about a Polish space that is not Borel cannot be solved with inherently countable information.\footnote{This is not original to this paper.} If an equivalence relation is not Borel a general  question whether $x_1\sim x_2$ requires \emph{some} uncountable resource--usually an application of the uncountable Axiom of Choice. Thus a general theme of this paper consists of determining when a question can be answered using inherently countable information: it focusses on the Borel/non-Borel distinction.

This language can be confusing:  for example if an equivalence relation $E$ can be ``classified by countable structures," it sounds like it is ``classifiable."  However even for countable structures, the isomorphism relation may not be Borel.  An example of this is the equivalence relation of \emph{isomorphism} for \emph{countable groups}.  (See section \ref{S-inf}.)
 
\paragraph{Dynamical Systems} Broadly speaking a \emph{dynamical system} is \emph{any} group action $\phi: G\times X\to X$. However, to be interesting, the group actions are taken to preserve some structure on $X$. Commonly studied types of structure include  topological, measure theoretic, smooth and complex. 
 
 There is a deep and well-developed theory for general groups, both amenable and non-amenable, discrete, topological, and carrying a differentiable or complex structure.  Accordingly the actions can be discrete, continuous or smooth actions. 
 
 Typical groups include $\poZ^n$, $\poR^n$,  more general Lie groups, free groups $\mathbb F_n$, and groups arising as automorphism groups of natural structures.  However much of the theory is modeled on initial successes with the groups $\poZ$ and $\poR$, so that will be the focus in this paper. In general the theory gets progressively more difficult as the groups move from $\poZ$ to $\poZ^d$ to general amenable groups, to free groups, and then to general non-amenable groups. Since part of the story being told here is that classifications can be ``impossible," showing this in the simplest, most concrete situation illustrates the point most dramatically.
 
The roots of the theory of much of dynamical systems can be traced to the study of vector fields on smooth manifolds, naturally linked to smooth $\mathbb R$-actions. (See \cite{Smale1963}.) As argued by Smale in \cite{Smale1967} these actions have smooth cross sections that give significant information about the corresponding solutions to ordinary differential equations.  The smooth cross sections are $\poZ$-actions by diffeomorphisms of the manifold.  In short, $\mathbb R$-actions give rise to interesting $\poZ$-actions. In ergodic theory, a class of $\poZ$-actions are induced by $\poR$-actions using the method of \emph{first returns}. This is studied via Kakutani equivalence and discussed in section \ref{kaku}.

Turning this around, $\poZ$-actions can induce $\poR$-actions using the method of \emph{suspensions}.  This technique allows lifting the complexity results from $\poZ$-actions to $\poR$-actions. The upshot is that $\poZ$ and $\poR$-actions are closely related.  For this reason Smale (\cite{Smale1967}) and others argue for studying 
$\poZ$-actions. Since the $\poZ$-actions are determined by their generator, this amounts to studying single transformations.

In summary: this paper will focus on $\poZ$-actions, which are determined by the generating transformation--in effect we are studying single transformations.

\paragraph{What structure and equivalence relations will be considered?}
While we will give examples of difficulties with structure theory in other contexts, the 
 structure  and  classification theory of transformations described in this paper breaks very roughly into the \emph{quantitative} theory and the \emph{qualitative} theory.  The ergodic theorem gives a framework for studying a function by repeated  sampling along an orbit of a transformation $T$.  As the number of samples grows the averages converge to the average of the function over the whole space provided the given transformation is ergodic. Hence it can be viewed as the quantitative theory. The appropriate equivalence relation is conjugacy by measure-preserving transformations.

In \cite{Smale1967}, Smale describes  the study of transformations (in particular diffeomorphisms of manifolds) up to topological conjugacy as the \emph{qualitative theory}.  Symbolic shifts are also studied up to homeomorphisms (and up to their automorphism groups).

The understanding of the complexity of classifications in the quantitative theory is better developed than the understanding of the complexity of classifications in the qualitative theory.  Ergodic theory classifications are discussed in sections \ref{measureiso}.  The nascent connections with the qualitative theory is discussed in sections \ref{homeoclass}, along with the many open problems. 

To repeat, the two equivalence relations this survey will focus on are:
	\begin{itemize}
		\item measure conjugacy,
		\item topological conjugacy.
	\end{itemize}
The first equivalence relation involves transformations that act on probability measure spaces $X, Y$, perhaps with other structural restrictions. (For example volume preserving diffeomorphisms of a compact manifold.) Two such transformations $S$, $T$ are  \emph{measure conjugate} (or measure isomorphic) if there is a measure isomorphism $\phi:X\to Y$ such that $S\circ\phi=\phi\circ T$ almost everywhere.

\hypertarget{topconintro}{The second equivalence relation} involves homeomorphisms acting on  topological spaces  
$X, Y$, perhaps with other structural restrictions ($S$ and $T$ might be required to preserve smooth structures on $X$ and $Y$.)  A pair $S, T$ are 
\emph{topologically conjugate} (or topologically equivalent) if there is a homeomorphism $h:X\to Y$ such that $S\circ h=h\circ T$.

\paragraph{Choices of examples} A theme of the survey is that many natural questions have intractable complexity, such as being infeasible using inherently countable information. For this to be most convincing the classes should be taken to be as clear and concrete as possible in each context.  In the context the qualitative behavior of diffeomorphisms, for example, this is the motivation for focussing on  single transformations, taking the manifold to be as simple as possible (say the 2-torus) and the diffeomorphisms to be $C^\infty$. In many situations, the results are more general, but that is not the goal.

\paragraph{Homeomorphisms} Often the class being given a structure theory consists of homeomorphisms with additional properties.  In addition to smooth structure, one can consider minimal homeomorphisms, or topologically transitive homeomorphisms.  An important example in section \ref{intodistal} is the collection of topologically distal transformations.  

\hypertarget{minimaltrans} {We define the first two classes as they arise in many contexts. Let $X$ be a topological space and $h:X\to X$ be a homeomorphism. Then $h$ is \emph{minimal} if every forward $h$-orbit of an $x\in X$ is dense. A weaker property is \emph{topological transitivity} which means that \emph{for some $x\in X$} the forward $h$-orbit of $x$ is dense.}

\paragraph{What a Classification  isn't.} 
\hypertarget{notclassy}{To make sense of what a classification is, it  might best be motivated by giving an example what a classification \emph{isn't}. }   For concreteness let us consider the space of measure-preserving transformations of the unit interval with the equivalence relation of measure isomorphism (conjugacy).  

Because there are only continuum many measure-preserving transformations, there are at most continuum many equivalence classes $[T]$.  By the Axiom of Choice, one can build a one-to-one function 

\[F:\{[T]: T\mbox{ is a measure-preserving transformation}\}\to \mathbb R.\]
By setting $C(T)=F([T])$,  we get a map giving complete numerical invariants to the equivalence relation of measure isomorphism.
\smallskip

This is of course NOT what is meant by a classification since simply evoking the Axiom of Choice gives no useful information whatsoever.

\paragraph{Requirements to be a structure theory or a classification}  A  structure theory or a classification   must be  effective or computable in some sense. 
There are at least \hypertarget{bigthree}{three common versions of effectiveness}:

	\begin{itemize}
	\item Computable \emph{feasible} algorithms.
	\item Computable using inherently finite techniques.
	\item ``Computable" using inherently countable techniques.
	\end{itemize}
We will say nothing about the first notion. Structure and Classification results that use inherently finite techniques are called \emph{computable} or \emph{recursive} classifications. This topic is well-covered in other literature (\cite{brayamp}, \cite{Weihrauch}) and we are content with giving a couple of examples.  Structure and Classification that uses inherently countable techniques are  \emph{Borel}. Techniques to show a classification is Borel, or not Borel are the main topic of this paper.

We note that showing a class does \emph{not} have a Borel structure or does \emph{not} have a Borel classification is a stronger result than showing that there is no such computable theory. For this reason this paper focusses on Borel classifications. There is much earlier work showing that various kinds of classifications are impossible using finitary techniques. These appear in places such as \cite{Mitin} and \cite{Ryzhikov}.

\smallskip
\paragraph{Benchmarking} The complexity of the existing classifications can be measured by  comparison with  benchmark equivalence 
relations coming from outside dynamical systems.  This benchmarking process can also give a rigorous method of saying that one classification is more complex than a second one. 
This benchmarking is called the theory of analytic equivalence relations and is discussed at length in section \ref{benchmarks}.  A very brief introduction  can also be found in  \cite{whatis} (several of the diagrams in this paper were motivated by very similar diagrams in \cite{whatis}).  There is an extensive literature about this topic, many of the results can be found in notes and books such as \cite{GaoLecture}, \cite{Gaobook}, \cite{Hj}, and \cite{kechris}.

\section{Examples of structure and classification results in  dynamical systems} 
To try to clarify the notions we begin by giving examples from dynamical systems.\footnote{Much of this section is modeled on the introductory  chapter of K. Petersen's book (\cite{peterson}).}

\subsection{Hamiltonian Dynamics} These flows were developed by Hamilton to describe the evolution of physical systems. (See \cite{wikipedia:HamSys}.)  A Hamiltonian system is described by a  twice differentiable 
($C^2$) function $H(\mathbf{q}, \mathbf{p})$ from  $\mathbb R^{6N}$ to $\poR$, giving the energy of the system.  The system is described by \emph{Hamilton's Equations}:

\begin{eqnarray*}
\frac{d\mathbf{p}}{dt}=-\frac{\partial H}{\partial \mathbf{q}},\\
\frac{d\mathbf{q}}{dt}=+\frac{\partial H}{\partial \mathbf{p}}.
\end{eqnarray*}
The variables $\mathbf{p}, \mathbf{q}\in \poR^{3N}$ are interpreted as the generalized position and momentum variables and   the solution $r(t)$ is viewed as the trajectory of a point
 in an initial position $r(0)\in \poR^{6N}$.  Putting this the context of group actions, the solution is a one-parameter flow
 \[T:\mathbb R\to \poR^{6N}.\]

What makes ergodic theory relevant to this example is \emph{Liouville's Theorem} that the Hamiltonian system $\{T_t\}_t$ preserves Lebesgue measure. For many fixed energies $E$ the surface $H=E$ is a compact manifold and thus, on this manifold, the flow $\{T_t\}_t$ preserves an invariant probability measure (Khinchin (\cite{Khintchine1949}).  

The Birkhoff ergodic theorem suggests the appropriate equivalence relation for studying systems with invariant measures:
\begin{theorem}
Suppose that $\xbmt$ is a measure-preserving system and $f\in L^1(\mu)$. Then the averages $\left(\frac{1}{N}\right)\sum_0^{N-1}f(T^n(x))$ converge almost everywhere to a $T$-invariant function $\bar{f}$ with 
	\[\int\bar{f}d\mu=\int fd\mu.\]  In particular if $T$ is ergodic then the averages converge to the constant value $\int\! fd\mu$.
\end{theorem}
An informal interpretation of this theorem is that if one repeatedly samples  an ergodic stationary system and then averages the results, the limit is the mean of the system. This suggests that the natural equivalence relation is 
\emph{conjugacy by a measure-preserving transformation}, because this equivalence relation takes repeated sampling along trajectories to repeated sampling along trajectories.  In different terminology, because \emph{measure isomorphism} take averages to averages, it is the natural equivalence relation to study in this context.

\subsection{Rotations of the unit circle}
In this example the class $\mcc$ is the collection of orientation preserving homeomorphisms of the unit circle.  The structure is evident (as $\mcc$ itself forms the relevant Polish space $X$).

 Let $S^1$ be the unit circle. It is convenient to view $S^1$ as the unit interval $[0,1]$ with $0$ and $1$ identified. (An equivalent approach is to view $S^1$ as $[0, 2\pi]$ with the endpoints identified.) 
 Let $\pi:\mathbb R\to S^1$ be the map $x\mapsto [x]_1$ where $[x]_1$ is the positive fractional part of $x$.

 If $f:S^1\to S^1$ is  an  orientation preserving homeomorphism, then a \emph{lift} of $f$ is an increasing function $F:\mathbb R\to \mathbb R$ where $[F(x)]=h([x])$. Using the notion of lift we can define the rotation number of $f$.\footnote{This discussion is largely based on  \cite{Barreira}, where the reader can see  details.}
	\begin{definition}\label{rotationnumber}
	Let 
	\begin{equation}\rho(f)=\lim_{n\to\infty}\frac{F^n(x)-x}{n}.\label{rotnum}
	\end{equation}
	Then $[\rho(f)]_1$  is called the \emph{rotation number} of $f$
	\end{definition}
 The following facts are necessary for this to make sense:
 	\begin{itemize}
		\item $\rho(f)$ exists and is independent of $x$,
		\item If $\rho^*$ is defined by equation \ref{rotnum} using a different lift $F^*$ of $f$, then 
		for all $x$, $[\rho^*(f)]_1=[\rho(f)]_1$.
	\end{itemize}
 Granting these facts, the rotation number is a well defined member of $[0,1)$ depending only on $f$.  We will denote this number by $\rho(f)$ rather than $[\rho(f)]_1$.

 The equivalence relation on homeomorphisms of the circle is conjugacy by a homeomorphism.  So $f\sim g$ if there is a homeomorphism $h$ such that 
   \[h\circ f=g\circ h.\]
 If $f$ is an arbitrary homeomorphism of the circle, then it is conjugate to a rotation preserving homeomorphism, so we focus on those.
   
 One can check that:
 	\begin{proposition}\label{numinv} If $f, g$ are orientation preserving homeomorphisms of $S^1$ that are conjugate by an orientation preserving homeomorphism, then $\rho(f)=\rho(g)$.
 	\end{proposition}
We will refer to statements like Proposition \ref{numinv} as saying that \emph{$\rho(f)$ is a {numerical invariant} for the equivalence relation of orientation preserving topological conjugacy}. Indeed having  numerical invariants is one of the key notions of classifiability.

The classification question for $\rho$ can then be phrased as:

	\begin{quotation}
	\noindent If $\rho(f)$ a \emph{complete} invariant?  	
	\end{quotation}
\noindent Explicitly: can inequivalent $f$ and $g$ have $\rho(f)=\rho(g)$?	
	\smallskip
	
Clearly the rotation of the circle by $2\pi\theta$ radians, $R_\theta$, has rotation number 
$\theta$. So the question of whether the rotation number is a complete invariant can be restated as asking whether $\rho(f)=\theta$ implies that $f$ is conjugate to $R_\theta$.

The following results show that $\rho(f)$ \emph{IS} a complete numerical invariant for functions with dense orbits.  The first is due to Poincar\'e.

	\begin{theorem}(Poincar\'e) \label{poincthm}
	Suppose that $f$ is a homeomorphism of $S^1$ with a dense 
	 orbit.  Then $f$ is topologically conjugate to $R_\theta$ for an irrational $\theta$.
	\end{theorem}

Denjoy showed that, if $f\in C^2$, then having an irrational rotation number implies having a dense orbit, however this is not true for general homeomorphisms. (See \cite{wikipedia:rotation}.)

Thus homeomorphisms of the circle with a dense orbit is a story of a \emph{successful} classification of the strongest kind:  it is a classification by a numerical invariant.  Moreover equation \ref{rotnum} shows that with a sufficiently computable $f$, the numerical invariant can itself be computed and that the rotation number is a continuous function from $X$ to $[0,1)$. (See \cite{HasKat}.)

The next theorem gives information about the situation when $f$ has a rational rotation number: 
 	\begin{theorem} 
	Let $f$ be an orientation preserving homeomorphism from $S^1$ to $S^1$. Then:
		\begin{enumerate}
		\item $f$ has a rational rotation number if and only if $f$ has a periodic point.
		\item If $\rho(f)=p/q$ and $(p,q)=1$ then $f$ has periodic points and 
		each periodic point has period $q$. 
		\end{enumerate}
	\end{theorem}
For homeomorphisms with period points the issue relates to the classification of order preserving homeomorphisms of the unit interval, which is studied in \cite{Hj}.
\paragraph{After this paper was written}   Joint work of the author and Gorodetski shows that even for diffeomorphisms, the equivalence relation of conjugacy-by-homeomorphisms turns out to be a maximal among equivalence classes reducible to $S_\infty$-actions and thus quite complex.  The complexity lies in diffeomorphisms with fixed points. 
 
\medskip

 
 \subsection{Symbolic Shifts} Let $\Sigma$ be a finite or countable alphabet.\footnote{Several authors refer to an alphabet as a ``language," a practice dating to Tarski.} We let 
 $\Sigma^\poZ$ stand for the collection of bi-infinite countable sequences $\la f(n):n\in \poZ\ra$ 
 written in the alphabet $\Sigma$.  
 
 The space $\Sigma^\poZ$ carries the natural product 
 topology.  With respect to this topology the shift map $sh:\Sigma^\poZ\to\Sigma^\poZ$ defined by:
 	\[sh(f)(n)=f(n+1)\]
	is a homeomorphism.
 A closed, shift-invariant $\bk\subseteq \Sigma^\poZ$ is called a \emph{subshift}.

 The structure and classification theories for symbolic shifts splits naturally into topological classifications and, when the symbolic shift has a specific measure associated with it, the classification up to measure isomorphism.  The latter is discussed in the next section.
 \medskip
 
While many classes of symbolic shifts are studied, perhaps the most common are those shifts in a finite alphabet. Among the classes studied is the collection of   \emph{Shifts of Finite Type}. This is the class of subshifts that are determined by a fixed finite collection of words $w_1, w_2, \dots w_k\in \Sigma^{<\nn}$. Given this collection the subshift consists of all $f\in \Sigma^\poZ$ such that for no  interval is $f\rest[n, n+m)$  equal to some $w_i$.

It is straightforward to verify that the collection of subshifts of finite type is a countable
subset of the compact subsets of $\Sigma^\poZ$. Hence the structure aspect is easy to understand.

Each shift of finite type is determined by finite information (the forbidden words) and hence it makes sense to ask whether it is possible to determine whether two subshifts are topologically conjugate using inherently finite information.  Various invariants exist which are derivable from \emph{adjacency matrices}  and two subshifts of finite type are conjugate by a homeomorphism if and only if their adjacency matrices are \emph{strong shift equivalent}. This is an example of a \emph{recursive reduction}. It reduces  the equivalence relation of conjugacy by homeomorphism to the equivalence relation on the collection of adjacency matrices of being strong shift invariant (See section \ref{reductions}).

However the question of whether there is an algorithm to determine conjugacy by homeomorphisms for shifts of finite type remains an open problem (see Open Problem \ref{OP2}).

While other invariants exist for aperiodic subshifts such as \emph{topological entropy},  these are not complete invariants for conjugacy by homeomorphism. (In \cite{CG} it is shown that the classification of homeomorphisms of the Cantor set up to conjugacy has no Borel invariants--in fact it is maximal for $S_\infty$-actions. See section \ref{camerlogao}.) 

In fact, Clemens \cite{clemens} was able to show that there are no Borel computable complete numerical invariants for conjugacy of subshifts.  Clemens determined the exact Borel complexity of conjugacy of subshifts.

 \subsection{Measure structure}  
 There is an important and very well developed study of transformations that 
 preserve Borel measures that are not probability measures \cite{Aaronson}, but to illustrate the connections with descriptive set theory,  we consider probability measures. 
  
 We will discuss two types of questions about the structure of measure-preserving transformations:  understanding a single transformation, and trying to understand the structure of factors and isomorphisms. In section \ref{measureiso} we discuss the classification problem for measure-preserving systems and smooth measure-preserving systems. 
 \footnote{In \cite{kechrisglobal}, Kechris  discusses a more general setting: Let $G$ be a countable non-amenable group. Fix a standard non-atomic probability space $(X, \mu)$ and let $A(G,X,\mu)$ be the space of measurable and $\mu$-preserving actions of $G$ on $X$. Note that $A(G,X,\mu)$ is a Polish space when equipped with the weak topology. In this survey we will be primarily concerned with $\poZ$ actions.}

    A basic object in ergodic theory is  $(X,\mcb,\mu, T)$ where 
 $(X,\mcb,\mu)$ is a standard probability measure space and $T:X\to X$ is an invertible measure-preserving transformation. Following Furstenberg \cite{FUbook}, we will call these measure-preserving systems (and $T$ a measure-preserving transformation
 \footnote{We will only be considering invertible measure preserving systems so we drop the adjective \emph{invertible}.} 
 and reserve that term for probability measure-preserving systems. 
 
 \hypertarget{mpt}{While there are many models for the space of measure-preserving transformations, one is the collection of Lebesgue measure-preserving transformations of the unit interval equipped with Lebesgue measure.  These are identified if they agree almost everywhere, and have the weak topology (the \emph{Halmos topology}), which carries a complete separable metric (\cite{halmos1944}, \cite{halmosbook}). For a standard probability space we will call the group of measure preserving transformations $MPT(X)$. With the weak topology it is a Polish group, with the operation of composition.} 
 
 Two measure-preserving systems \xbmt\ and \ycns\ are (depending on the author) \emph{isomorphic} or \emph{conjugate} if there is an invertible measure-preserving transformation 
 $\phi:X\to Y$ whose range has $\nu$-measure one such that for $\mu$-almost all $x\in X$
 \begin{equation}\label{conjugal}
 S\circ \phi(x)=\phi\circ T(x).
 \end{equation}
  If $\phi$ is not invertible but still satisfies equation \ref{conjugal} almost-everywhere, then $\phi$ is a \emph{factor map} 
  and $S$ is a factor of $T$.
  
If $M$ is a compact metrizable space and $T:M\to M$ is a homeomorphism, then the Krylov-Bogoliubov Theorem shows that there is a T-invariant probability measure $\mu$ on the Borel subsets of $M$.  The space $P_T$ of $T$-invariant measures on $M$ is convex and  compact with the weak topology.  Hence the Krein-Milman theorem says that it is spanned by its extreme points.\footnote{The interested reader can find explanatory discussions of these topics in Walters \cite{Walters} and Mirzakhani-Feng \cite{Mirza}.}

A measure-preserving transformation is \emph{ergodic} if the only invariant measurable sets  either have
measure zero or measure one.  The extreme points of the simplex of $T$-invariant measures on $M$ are exactly the measure for which $T$ is ergodic. (We will call these the \emph{ergodic measures} when $T$ is clear from context.)
The Ergodic Decomposition Theorem says that if $\mu$ is a $T$-invariant measure then $\mu$ 
can be written as a direct integral of the ergodic measures. For the statement of the theorem we let 
$\mathcal M(X)$ to be the space of invariant probability measures on $X$ with the weak topology.

\begin{theorem}
Let $T$ be a measure-preserving transformation of a standard probability space $(X,\mcb,\mu)$.  Then there is a measure space $(\Omega,\mcc,\nu)$ and a $\mcc$-measurable map $\mu^*:\mcc\to \mathcal M(X)$ such that:
	\begin{itemize}
	\item For all $w\ne w'$, $\mu^*(w)$ and $\mu^*(w')$ are mutually singular and $T$-invariant.
	\item For all $Y\in\mcb, \mu(Y)=\int\left(\mu^*(w)(Y)\right)d\nu(w)$.
	\item For almost all $w\in \Omega, (X,\mcb, \mu^*(w), T)$ is ergodic.
	\end{itemize}
Moreover the measure $\nu$ is unique up to sets of measure zero.
\end{theorem}

 Thus every invariant measure on $M$ can be written as a direct integral of ergodic measures.  We interpret this as saying that the ergodic measures are the building blocks of the invariant measures and focus the classification problem on the ergodic measures, or equivalently, on the ergodic measure preserving transformations.

 \paragraph{Notation} We use the notation $EMPT(X)$ for the collection of ergodic measure-preserving transformations on $(X,\mcb,\mu)$. 
 \medskip
 
 The property of being ergodic is a conjugacy invariant so it invariant under the \hypertarget{ergisGd}{conjugacy action of $MPT(X)$.} Moreover Halmos showed that  the ergodic measure-preserving transformations are a dense $\mcg_\delta$ set in the weak topology on $MPT(X)$. (See  \cite{halmosbook}.)

 \subsubsection{Rolling Dice: Bernoulli Shifts}\label{bs} Two classes of measure-preserving systems provide very nice examples of both  structure and a classification.  One is the full shift on a finite alphabet with the Bernoulli measure.  The other is Haar measure on a compact group.  The latter is discussed in the next section.
  
	\begin{definition}\label{Bshift} Let $\Sigma$ be a finite alphabet and
  	$\nu$ be a measure on $(\Sigma, P(\Sigma))$.  Let 
	$\mu$ be the product measure $\nu^\poZ$.  Then $\mu$ is a shift-invariant measure on 
	the full shift $\Sigma^\poZ$.  The system $(\Sigma^\poZ, \mcb, \mu, sh)$ is called a 
	\emph{Bernoulli Shift}.
  	\end{definition}
 Informally: the coordinates of the Bernoulli shift are independent, identically distributed (iid). In particular the ``time-zero" behavior is independent of the past.

\paragraph{Structure theory for Bernoulli shifts} \hypertarget{firstbernoulli}{There is a successful structure theory for Bernoulli shifts.}   Ornstein showed (\cite{FDimpsBern}) showed that being \emph{finitely determined} is equivalent to being isomorphic to a Bernoulli shift and that being \emph{very weak Bernoulli} implies being finitely determined.   Ornstein and Weiss (\cite{FDimpsBern}) proved the converse: being finitely determined implies being very weak Bernoulli. Feldman (\cite{feldman}) used this characterization to show:
\begin{theorem}(Feldman)\label{BorelBorel}
The collection of measure-preserving transformations on the unit interval that are isomorphic to a Bernoulli shift on a finite alphabet is  a Borel set.
\end{theorem}

Being isomorphic to a Bernoulli shift is  an \hyperlink{aandcoa}{analytic} condition, as we sketch in example \ref{nevertoolate}.  Feldman showed that the collection of very weak Bernoulli transformations is \hyperlink{aandcoa}{co-analytic}. The \hyperlink{suslintheorem}{ theorem of Suslin} (Corollary \ref{afterall}) says that a collection which is both analytic and co-analytic is Borel.  Thus the collection of Bernoulli shifts is Borel.

 \paragraph{Classification theory for Bernoulli shifts}
Kolmogorov and Sinai, building on the work of Shannon defined a notion of \emph{entropy} for measure-preserving transformations. Viewed as a  function defined on the ergodic measure preserving transformations, \hyperlink{mpt}{EMPT},  \emph{entropy} is a lower-semicontinuous.  Kolmogorov  showed that isomorphic Bernoulli shifts have the same entropy.

  This sets the stage for the following result of Ornstein, which is exposited in \cite{ornstein1}, \cite{ornstein2} and \cite{ornstein3}:
 
 \begin{theorem}\label{OrnErg}(Ornstein)
 Two Bernoulli shifts are isomorphic if and only if they have the same entropy.
 \end{theorem}
 Thus the Bernoulli shifts are a paradigm of both a successful structure and a successful classification theory: 
	 \begin{enumerate}
 	\item (Structure Theory) Bernoulli shifts can be identified in a Borel manner (they are the finitely determined transformations).
	\item (Classification Theory) Entropy is a complete invariant for isomorphism.
 	\end{enumerate}
 
 \medskip
 
 Among the various differences between $\poZ$-actions and the actions of other groups is that the question of whether Bernoulli shifts based on $F_2$-actions forms a Borel set among the $F_2$-actions is an open problem.  (See Open Problem \ref{OP3})
 
 \bigskip
 
  A notion related to Bernoulli shifts is that of a Kolmogorov Automorphism, or for short, a 
 $\mck$-automorphism.  Informally: Bernoulli Shifts are characterized by the present being independent from one moment earlier.  $\mck$-automorphisms are characterized by the present being asymptotically independent of the past.
 	\begin{definition}\label{Kaut} A $\mck$-automorphism is a measure-preserving transformation $T$  of a standard probability space $(X,\mcb, \mu)$ such that there is a sub-$\sigma$-algebra $\mck\subseteq \mcb$ with:
		\begin{enumerate}
		\item $\mck\subseteq T\mck$,
		\item the union of the algebras $T^n\mck$, for $n\ge 0$ generates $\mcb$ as a $\sigma$-algebra,
		\item $\bigcap_{n=1}^\infty T^{-n}\mck=\{\emptyset, X\}$.
		\end{enumerate}
	\end{definition}
There is a perfect set of non-isomorphic $\mck$-automorphisms, as was shown by Ornstein and Shields in \cite{ornsteinshields}. 
In section \ref{e0toK} it is shown  that the very elegant classification theory for Bernoulli shifts does NOT extend to the 
$\mck$-automorphisms and little is known about classifying the $\mck$-automorphisms. For example it is not known if the measure-isomorphism relation restricted to the $\mck$-automorphisms is Borel. (See Open Problem \ref{OP7}.)

\subsubsection{Symbolic systems as models for measure-preserving transformations.}

Let $(X, \mcb, \mu)$ be a standard probability space, and $T:X\to X$ be an ergodic 
measure-preserving transformation.  Then a set $\mca\subseteq \mcb$ is a \emph{generating set} if $\mcb$ is the smallest $T$-invariant $\sigma$-algebra containing $\mca$ (where we identify sets if their symmetric difference is zero).

If $\mca$ is a countable or finite generating set then we can make the elements of $\mca$ disjoint and keep it a generating set. For this reason, the terminology \emph{generating partitions} is frequently used. 
If $\mca=\la A_n:n\in\nn\ra$ is a  partition of a set of measure one, then almost every $x$ in $X$ determines a $\poZ$-sequence of \emph{names}
 $\la A_{n_k}:k\in\nn\ra$ where $T^k(x)\in A_{n_k}$. The sequence is called the $\mca$-name of $x$. If $\mca$  generates, then there is a set $S$ of measure one such that for all $x\ne y\in S$ the $\mca$-name of $x$ is different from the $\mca$-name of $y$.
 
 Hence on $S$ there is a \oto map $\phi$ from $X$ to $\mca^\poZ$ given by sending $x$ to its $\mca$-name. Letting $\nu=\phi^*(\mu)$ (so copying the measure on $X$ over to a measure on $\mca^\poZ$) we get an isomorphic copy of $X$ as a shift space $(\mca^\poZ, \nu)$.  
 \smallskip
 
 Rokhlin (\cite{therokh}) proved that there is always a countable generating set and Krieger (\cite{krieg1} and \cite{krieg2}) proved that if the entropy of $T$ is finite, then there is a finite generating set.  
 
 Krieger (\cite{krieg3}) showed the stronger result that if the original transformation $T$ is ergodic and $K\subseteq \mca^\poZ$ is the support of 
 $\nu$ then  $(K, sh)$ is uniquely ergodic. We have outlined the proof of the following theorem
 
\begin{theorem}\label{symbmpts}
Let $(X, \mcb, \mu)$ be an ergodic measure-preserving system. Then there is a countable alphabet 
$\mca$ and a closed shift-invariant set $\bk\subseteq \mca^\poZ$ carrying a unique shift-invariant measure 
$\nu$ such that $(\bk, \mcc, \nu, sh)$ is isomorphic to $\xbmt$. If $\xbmt$ has finite entropy then the 
alphabet $\mca$ can be taken to be finite.
\end{theorem}

\subsubsection{Koopman Operators} Let $\xbmt$ be a measure-preserving system.  Define a map
\[U_T:L^2(X)\to L^2(X)\]
by setting $U_T([f])=[f\circ T]$. Then $U_T$ is a well-defined unitary operator on $L^2(X)$ and moreover if $S$ and $T$ are isomorphic measure-preserving transformations then $U_S$ and $U_T$ are unitarily equivalent.

Thus it is possible to assign to each measure-preserving transformation a unitary operator in a way that sends conjugate transformations to unitarily equivalent operators.  These operators are called the \emph{Koopman Operators}.  The equivalence classes of Koopman operators are thus a collection of invariants, but except for very special cases (such as in the next subsection) they are not complete invariants.

\subsubsection{Translations on Compact Groups}\hypertarget{HvN}{We describe the \emph{Halmos-von Neumann Theorem}.} Let $G$ be a compact group. Then $G$ carries a Haar probability measure that is invariant under left multiplication.  Suppose that $g\in G$ and Haar measure is ergodic for the map $T_g:G\to G$ given by $T_g(h)=gh$.   Then $G$ must be the closure of the powers of $g$: $G$ is a monothetic abelian group.\footnote{The material presented in this example is explained in \cite{descview}.}

An ergodic transformation $T:X\to X$ is said to have \emph{discrete spectrum} if $L^2(U)$ is spanned by the 
eigenfunctions of $U_T$. The eigenfunctions can be taken to have constant modulus one and are closed under 
multiplication. Since $T$ is measure-preserving the eigenvalues have modulus one and form a subgroup 
of the unit circle.  Hence by Pontryagin duality their dual is a compact monothetic abelian group $G$.  The 
identity map from the set of eigenvalues into the unit circle is a generator of this group.  If $g$ is the identity character and $T_g$ is translation by $g$, then there is a unitary operator conjugating $U_{T_g}$ with $U_T$ that preserves multiplication of bounded functions in $L^2(G)$. From this one can compute a measure isomorphism between $(G, T_g)$ and $(X, T)$.

Again by Pontryagin duality, if $G$ is a compact group, then $L^2(G)$ is spanned by the characters of $G$.  An immediate calculation shows that if $G$ is monothetic, with generator $g$ and $T$ is translation by $g$, then the characters of $G$ are eigenvalues--hence $T$ has discrete spectrum.

The Halmos-von Neumann theorem (\cite{HvN}) gives both a structure theorem and a classification theorem.

\begin{theorem}(Halmos-von Neumann)\label{hvntheory} Let $T$ be an ergodic transformation.  Then:
	\begin{enumerate}
	\item (Structure) $T$ has discrete spectrum if and only if it is isomorphic to an ergodic translation on a compact group with Haar measure.
	\item (Classification) If $S, T$ have discrete spectrum, then $S$ and $T$ are conjugate by a measure-preserving transformation if and only if the Koopman operators associated to $S$ and $T$ have the same eigenvalues.
	\end{enumerate}
\end{theorem}

For a proof of this theorem and more on this topic see \cite{Walters}.

 \subsection{Presentations} Often transformations with various relevant properties one wants to understand can be presented in different contexts. 
 The context in which they are presented can be important. For example, two isomorphic measure-preserving 
 transformations, presented with measures on different spaces  can have different simplexes  of invariant measures. Another example is Theorem \ref{symbmpts} that shows that every measure-preserving transformation can be presented as a symbolic shift.

Striking examples include Bernoulli Shifts: there are Bernoulli shifts that are presented as  real-analytic measure-preserving transformations of the torus (with the Poulsen simplex of invariant measures). Bernoulli shfts can also be presented as uniquely ergodic symbolic systems. The issue of presentation is quite different than the issue of isomorphism, although topological conjugacy is quite sensitive to presentation.
 
 Presentations \emph{can} affect intrinsic measure theoretic properties.  For example Kushnirenko showed that a diffeomorphism of a compact manifold preserving a smooth volume element has finite entropy.  One of the most prominent open problems in the area is the converse: can every finite entropy measure-preserving transformation be presented as a $C^\infty$,  transformation on a compact manifold that preserves a volume form? (See Open Problem \ref{OP5}.)

 \subsection{Smooth dynamics}\label{smodyn} In the context of this paper \emph{smooth dynamics} will mean studying diffeomorphisms of compact manifolds.  The equivalence relation is that of conjugacy by homeomorphisms.  We repeat the \hyperlink{topconintro}{definitions from the introduction.} 
 
 Explicitly, we let $M$ be a smooth manifold, usually taken to be compact. We will be given a class $\mcc$ of diffeomorphisms of $M$  described with a structure theory. Two elements $f, g\in \mcc$ will be taken to be \emph{equivalent} if there is a homeomorphism $h:M\to M$ such that 
 \begin{equation} f\circ h=h\circ g. \label{homeoconj}
 \end{equation}
 
Smale (\cite{Smale1967}) calls this equivalence relation the \emph{qualitative theory} of diffeomorphisms. He says that ``the same phenomena and problems of the qualitative theory of  ordinary differential equations [\emph{are}] present in the simplest form in the diffeomorphism problem."
He argues that we should study the equivalence relation given when the conjugacy map $h$ is a homeomorphism, rather than being a diffeomorphism, because ``differential conjugacy is too fine" to reveal the relevant properties such as structural stability.

\paragraph{Some examples of classes of diffeomorphisms} There are many well-studied classes of diffeomorphisms.  For the purposes of this paper we will assume that they are $C^\infty$ and usually assume that the manifolds $M$ are compact. We start with two of the well-known classes of diffeomorphisms.
\medskip

Suppose $f:M\to M$ is a diffeomorphism and that $\Lambda\subseteq M$ is invariant.
Then $\Lambda$ is \emph{hyperbolic} for $f$ if there are constants $C>0$ and $0<\lambda<1$ and a decomposition of the tangent bundle restricted to $\Lambda$ as:
\begin{align*}
TM_\Lambda\ &=E^s\oplus E^u\\
Tf(E^u)&=E^u\\
Tf(E^s)&=E^s
\end{align*} 
and
for $v\in E^s, n>0$
  \[\|Tf^n(x)v\|\le C\lambda^n\|v\|\] 
and for $v\in E^u, n\ge 0$
\[\|Tf^{-n}(x)v\|\le C\lambda^n\|v\|\]

\begin{definition}\label{anosovdiff} A diffeomorphism $f:M\to M$ is \emph{Anosov} if the whole manifold $M$ is hyperbolic for $f$.
\end{definition}

 \medskip
 
 Suppose that $x$ is a periodic point of order $m$ and  $x, f(x), \dots f^{m-1}(x)$ are hyperbolic points. Then the \emph{stable} and \emph{unstable} manifolds of $x$ are:
\begin{align*}
W^s&=\{y\in M|f^{mk}(y)\to x\mbox{ as k}\to \infty\}\\
W^u&=\{y\in M|f^{-mk}(y)\to x\mbox{ as k}\to \infty\}
\end{align*}
\noindent These are immersed Euclidean spaces of dimension corresponding to the number of eigenvalues of the differential $Df^m(x)$ that are greater or less than one in absolute value.
\begin{definition}\label{morsesmaledef} A diffeomorphism $f:M\to M$ is \emph{Morse-Smale} if and only if there are a finite collection of periodic orbits $P_1,\dots P_l$ such that:
 	 \begin{enumerate}
 	 	\item $P_i$ is a hyperbolic periodic point for $i=1, \dots l$,
  		\item $\bigcup_{i=1}^l W^s(P_i)=M$,
  		\item $\bigcup_{i=1}^l W^u(P_i)=M$,
  		\item $W^u(P_i)$ and $W^s(P_j)$ are transversal for all $i,j$.
  	\end{enumerate}
 \end{definition}
 
 \medskip
 
\begin{enumerate}
	\item (Structure) Both the Anosov and the Morse-Smale diffeomorphisms are examples of \emph{structurally stable} diffeomorphisms.  The structurally stable diffeomorphisms are those $f$ for which there is a $C^1$-neighborhood $U$ such that every $g\in U$ is topologically conjugate to $f$. Thus 
	\begin{enumerate}
	\item the Anosov diffeomorphisms form an open set,
	\item the Morse-Smale diffeomorphisms form an open set
	\item the structurally stable diffeomorphisms form an open set. 
	\end{enumerate}
 It follows that there are only countably many classes of Anosov and Morse-Smale diffeomorphisms. 
 
 It would be a nice picture if the collection of structurally stable diffeomorphisms formed an open and 
 \emph{dense} subset of the diffeomorphisms, but this was refuted by Newhouse \cite{newhouse}, \cite{newhouse2}.
 
	\item (Classification) Each of the three classes form an open subset of the space of diffeomorphisms. Each of the open sets decompose into a union of open subsets corresponding to each of the equivalence classes. Thus for each class, there is a single countable list of representatives and radii that capture all of the members of the class. 
The representatives can be taken to be rational (in the appropriate sense) and the radii can be take to be of the form $1/n$.  Hence the countable lists can be viewed as members of a Polish space.  

The upshot is that for each of the three classes, one can assign complete numerical invariants in a Borel (even continuous) way. The equivalence relation is Borel reducible to the equality relation on a Polish space.
\end{enumerate}
	
\noindent	In dimension 2 more intuitive invariants can be assigned in a computable way.  For Anosov diffeomorphisms on the 2-torus these are hyperbolic elements of SL$_2(\mathbb Z)$ and for Morse-Smale diffeomorphisms, they are graphs with colored edges (\cite{peixoto}, \cite{shark}). Similar results hold for Morse-Smale diffeomorphisms in dimension 3 (\cite{bonatti}).

 An alert reader will notice a problem with both the structure and classification results just stated. The 
 program is to classify equivalence classes of diffeomorphisms, and the definition of Anosov, Morse-Smale and structurally stable diffeomorphisms are not invariant under topological conjugacy.
 Ideally, given a diffeomorphism $f$ that is topologically conjugate to an element of one of 
 these classes, there would be a way of constructing a $g$ equivalent to $f$ that actually belongs to 
 the class.  This seems to be an open problem. (See Open Problem \ref{OP6}.) 

\subsection{More detailed structure theory}
 So far the term ``structure theory" has been used solely for the purpose of determining membership in a class.  However, in several cases a structure theory gives more information--for example about the factor structure of a member of a class.  The Furstenberg structure theorem for topologically distal transformations and the Furstenberg-Zimmer theory for measure-preserving transformations are examples of this.  These are discussed in section \ref{intodistal}.
 \subsection{Summary} 
The previous sections   are intended to give a context for \emph{structure theory} and \emph{classification theory}. These are somewhat vague terms, but the intention is the following:
	\begin{itemize}
	\item Structure theory takes place in the context of an ambient Polish space and gives a method for determining membership in a class $\mcc$. It can also give information about the factors of a transformation or explicit information about the way the transformation moves elements in the space it acts on. 
	\item Classification theory considers a class $\mcc$ and an equivalence relation $E$ on that class. It is concerned with determining when two elements of $\mcc$ are $E$-equivalent.
	
	\end{itemize}

\noindent A structure theorem is concerned with a single transformation, and classification theory is concerned with pairs of transformations.

\section{Descriptive Complexity} 
For a structure theory or classification theory to be useful it is important that it be computable in some sense.  As the example given in section \ref{whatisdyn} that \hyperlink{notclassy}{used the Axiom of Choice} illustrated, classifications or structure theorems that are not in some sense computable are likely to be meaningless.  

The structure theoretical levels of complexity are easier to describe than the \emph{complexity benchmarks} for classification, as it is a series of yes-no questions. Recall from the discussion in the introduction the \hyperlink{bigthree}{three basic questions} about how effective a structure theory can be are:
	\begin{enumerate}
	\item Is the structure theory {computable with realistic resources assumptions?}
	\item Can it be carried out with inherently finite information?
	\item Can it be carried out with inherently countable information?
	\end{enumerate}

We will ignore the first question and say little about the second.  The {main} point of the paper is that there are natural examples 
 of dynamical systems fail even the third question.  
The third question can be rephrased as asking whether the structure theory shows that the class 
$\mcc$ is Borel.\footnote{Of course, once a structure theory has been shown to be Borel, the question arise to determine what level in the Borel hierarchy it lies in.  There is such a literature, but it is not addressed in this survey.}

 \begin{remark*}
 In this survey we only consider equivalence relations that are \hyperlink{aandcoa}{\emph{analytic}}.  
 These include all orbit equivalence relations of Polish groups. (See example \ref{nevertoolate}.)
 Some natural equivalence relations are not analytic (such the equivalence relation on distal flows of having the same distal height).  The theory is not nearly as well developed for these so they are not covered in this survey.
 
 \end{remark*}
\subsection{What is a reduction?}\label{reductions}
There are many jokes about mathematicians where the punchline is that the way they solved problem $A$ was to  ``reduce it to a problem $B$ they already knew how to solve."\footnote{A google search in 2018 gave over 115 million hits to the search ``reduce it to a problem you can solve joke."}

The idea behind the joke is that to solve a problem $A$ you (somewhat constructively) reduce it to a problem $B$ you already know how to solve. One can state this mathematically in the following way.  

\subsubsection{Reductions of sets}

Let $X$ and $Y$ be Polish spaces and $A\subseteq X$, $B\subseteq Y$.  Then $A$ is
\begin{center}{\underline{  \ \ \ \ \ \ \ \ \ \ \ \ \                  }- reducible}
\end{center}
\smallskip 

\noindent to $B$
provided that there is a function $f:X\to Y$ such that for all $x\in X$
	\begin{center}
	$x\in A$ if and only if $f(x)\in B$.
	\end{center}

\bfni{Fill in the blank} As the \hyperlink{notclassy}{example given earlier} showed, without some restrictions on the function $f$ this notion is not meaningful. If $B$ and $Y\setminus B$ are non-empty then every 
set $A\subset X$ is reducible to $B$: fix $b_0\in B$ and $b_1\in Y\setminus B$ and define $f(x)=b_0$ if $x\in A$ and $f(x)=b_1$ if $x\in X\setminus A$. Then $f$ reduces $A$ to $B$.
Thus one must fill in the blank with some condition on the concreteness of $f$.

In the complexity theory world of computer science, the blank is often filled in by some statement ``linearly",  ``polynomially", ``exponentially" or ``recursively." In the context of this paper we focus on Borel functions $f$.

Starting in the 1950's, in the subject then known as Recursion Theory (\cite{soare}) and now called Computability Theory, the joke was reversed: 
	\begin{quotation}
\noindent	To show a problem $B$ is not solvable, you find an unsolvable problem $A$ and reduce it to $B$ by a \underline{  \ \ \ \ \ \ \ \ \ \ \ \ \                  } function.
	\end{quotation}
If $B$ were solvable, and $A$ is reduced to $B$, $A$ would also be solvable, yielding a contradiction.
To make the definition complete, one has to fill in the blank, usually from one of the three types of functions: feasibly computable, recursively computable or Borel.

Indeed this has a long history: the unsolvability of the word problem for finitely presented groups was shown by reducing an incomputable set to the word problem using a computable function.
\medskip

The focus of this survey is  the notion of \emph{Borel Reducibility},\footnote{This idea originated in the late 1980s in the work of Friedman and Stanley (1989) related to model theory and independently Harrington, Kechris, and Louveau who used in it the context of analysis.}
 with some very short excursions into {continuous} reducibility. To repeat the official definition, we let $X$ and $Y$ be Polish spaces.  A set
	\begin{center} $A\subseteq X$ is \emph{Borel reducible} to  $B\subseteq Y$
	\end{center}
 if and only iff 
there is a Borel measurable function $f:X\to Y$ such that for all $x\in X$
	\begin{center}
	$x\in A$ if and only if $f(x)\in B$.
	\end{center}
Thus the inverse image of $B$ under $f$ is $A$. 
Because the definition of \emph{Borel measurable} implies that the inverse image of a Borel set is Borel, it follows that 
	\begin{center}
	if $A$ is not Borel, then $B$ is not Borel.
	\end{center}
	
\bfni{Notation:} We write $A\preceq^1_\mcb B$ if $A$ is Borel reducible to $B$.
In this notation, the previous remarks can be written for the record as: 
\begin{remark}\label{toshownotborel} Let $X, Y$ be Polish spaces. 
Let $A\subseteq X$ and $B\subseteq Y$, with $B$ Borel.  If $A\preceq_{\mathcal B}^1 B$ then $A$ is Borel.
\end{remark}

\subsubsection{Reductions of equivalence relations} For the classification problem a two dimensional version is  more useful. Let $X$ and $Y$ be Polish spaces and $E\subseteq X\times X$, $F\subseteq Y\times Y$ be equivalence relations.  Then
	\begin{center} $E$ is Borel reducible to $F$
	\end{center}
	 if and only if there is a (unary) {Borel} function $f:X\to Y$ such that for all $x_1, x_2\in X$
	 \begin{center}
	 $x_1Ex_2$ if and only if $f(x_1)Ff(x_2)$.
	 \end{center}
We use the notation $E\preceq^2_\mcb F$.
\medskip

It is important to recognize how the two definitions differ.  In the second, two dimensional definition, $f$ has domain $X$, not $X\times X$. So, in essence: 
	\begin{quotation}
	\noindent $f$ assigns an $F$-equivalence class to each $x$ in $X$ in a manner that two 
	elements of $X$ get assigned the same $F$-equivalence class if and only if they are in the 
	same $E$-equivalence class.
	\end{quotation}
\noindent Viewing the classes of $F$ as invariants, this is  paraphrasing of the statement that $f$ computes the invariants assigned to members of $X$.

Thus for $\preceq^2_\mcb$, $f$ is a function that assigns values to pairs $(x_1, x_2)$ by considering each $x_i$ separately. However viewed as a function from the Polish space $W=(X\times X)$ to the Polish space $Z=(Y\times Y)$ it is also a one-dimensional {reduction}.
\begin{remark}
Suppose that $E\subseteq X\times X$ and $F\subseteq Y\times Y$, then $E\preceq^2_\mcb F$ implies $E\preceq^1_\mcb F$.
\end{remark}

\subsubsection{Continuous reductions} An important special case of Borel reducibility is when the function $f$ is continuous.   This is stronger than being Borel reducible. Essentially every general statement below applies to continuous reductions as well as Borel reductions. If $A$ is continuously reducible to $B$ we write $A\preceq^1_C B$ and if $E$ is {continuously} reducible to $F$ we write $E\preceq^2_C F$.  If continuity is clear from context we omit the subscript $C$ and if the number of variables is clear from context we will omit the superscript in both cases.

\subsubsection{The pre-ordering}  For both Borel and continuous reductions and both the one and two dimensional cases, the relation of $\preceq_\mcb$ and $\preceq_C$ is transitive. To see this we note that if $f:X\to Y$ is a reduction (in the appropriate sense) and $g:Y\to Z$ is a reduction (in the same sense), then $g\circ f$ is a reduction in the same sense.

\hypertarget{bired}{As with all pre-orderings, we can have \emph{bi-reducible sets} and bi-reducible equivalence relations.} There are examples where $A\preceq_\mcb B$ and $B\preceq_\mcb A$ but there is no bijection $f:X\to Y$ such that $f$ reduces $A$ to $B$, and $f^{-1}$ reduces $B$ to $A$. (The analogue to the Cantor-Schroeder-Bernstein theorem fails.) We use the notation $A\sim_\mcb B$ (or $\sim_C$ or \dots) to mean the equivalence relation of bi-reducibility: $A\preceq_\mcb B$ and $B\preceq_\mcb A$.

Let $\sim$ be  the  relation $A\preceq_\mcb B$ and $B\preceq_\mcb A$.  Then $\preceq_B$ gives a partial ordering of $\sim$-classes of sets or equivalence relations depending on the dimension one is working in.

\begin{figure}[h]
\centering
\includegraphics[angle=90, origin=c, height=.9\textheight]{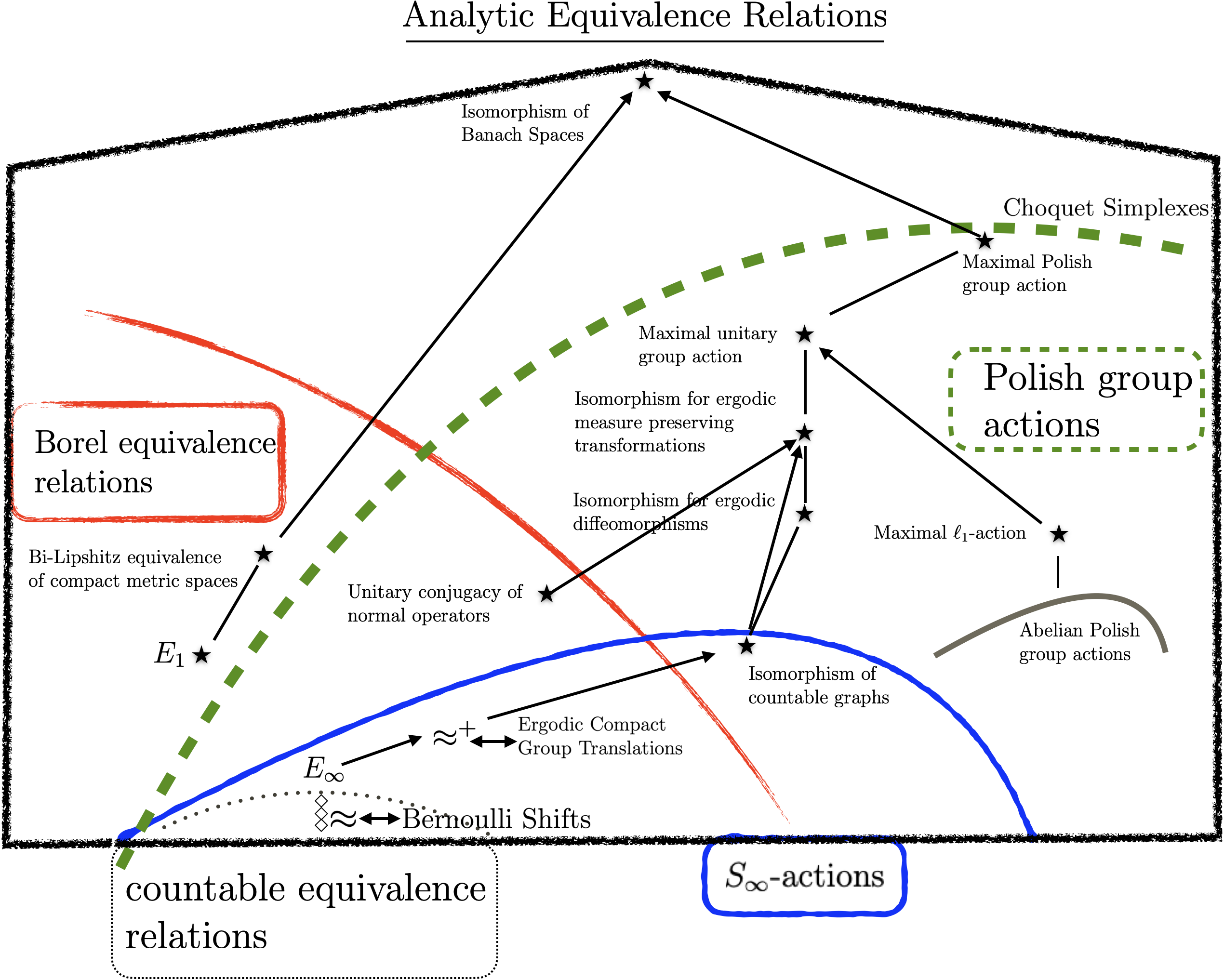}
\caption{measure-preserving systems }
\label{MPTs}
\end{figure}

\begin{figure}[h]
\centering
\includegraphics[angle=90, origin=c, height=.9\textheight]{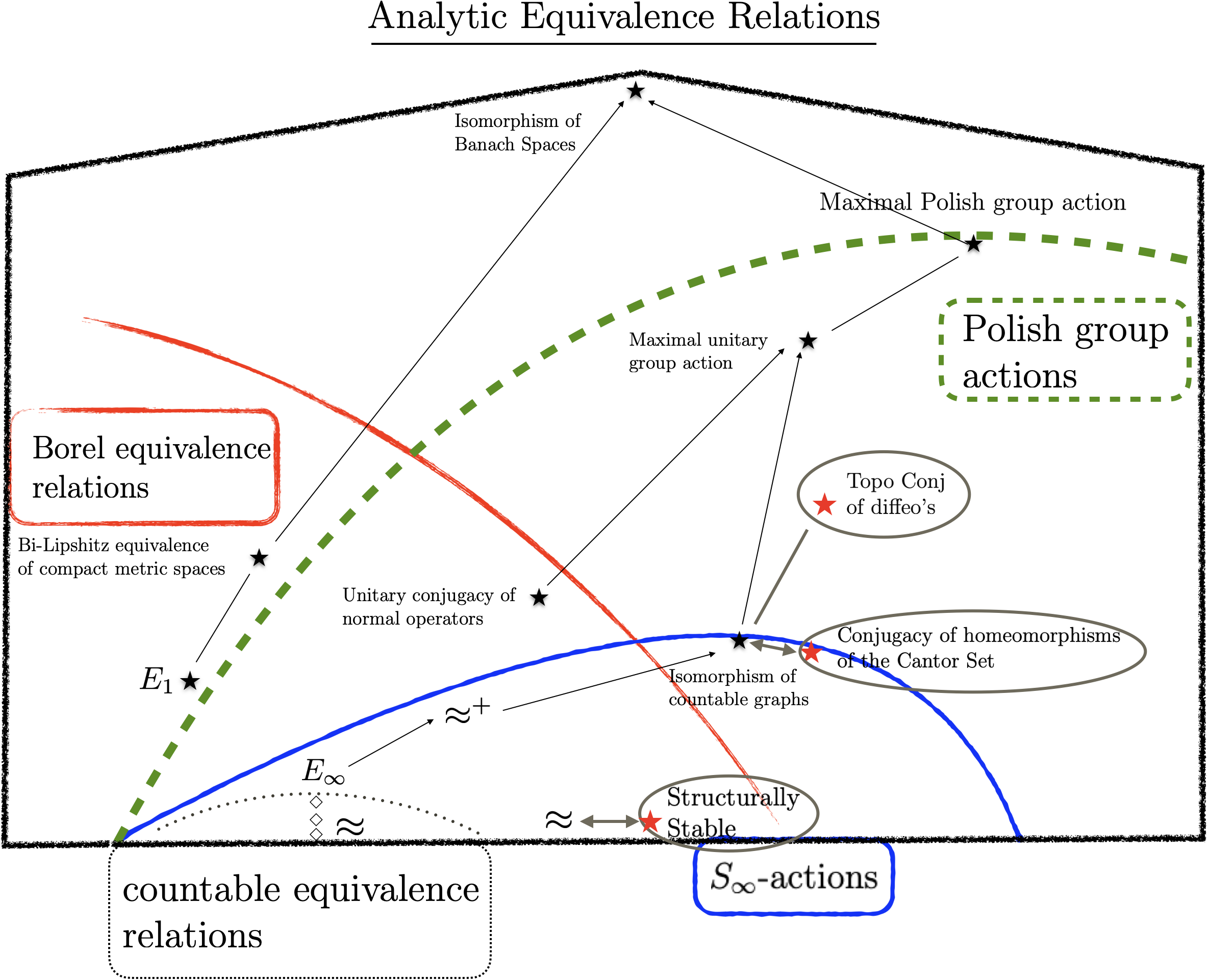}
\caption{Conjugacy by Homeomorphisms }
\label{topocon}
\end{figure}

\paragraph{Caveat!} There are several diagrams in this paper that give the relationship of reducibility between various classes of objects. The regions signify structures that are reducible to the type of object used to label the region.  So, for example, an object is in the region \emph{Polish group actions} if it can be reduced to a Polish group action. In particular this applies to figures \ref{yayayaya}, \ref{zoonodyna}, \ref{MPTs}, and \ref{topocon}.  In the diagrams, a line segment going upwards means \emph{Borel Reducible} and a pointed arrow means that the downwards reduction does not hold.  Not having an arrow means that bi-reducibility is open. An arrow at both ends indicates that bi-reducibility does hold.

\paragraph{Universality} Given a class $\mcc$ of either sets (in the one dimensional case) or equivalence relations (in the two dimensional case) an object $B\in \mcc$ is \emph{universal} for $\mcc$ if every $A\in \mcc$ is reducible to it. The interpretation of this is that $B$ is the most complex element of $\mcc$.   This is used for all of the notions: Borel reductions, continuous reductions and computable reductions. 

\paragraph{Terminology warning:} \emph{Maximal} is often used as a synonym of 
\emph{universal}. Less frequently \emph{complete} is also used. We will use \emph{maximal} 
when talking about the $\preceq^2_\mcb$ relation on equivalence relations, the 
two dimensional case.  We will use  \emph{complete} for $\preceq^1_\mcb$ restricted to sets, the one 
dimensional case.

Because of the transitivity of the notions of reducibility, if $B, C\in \mcc$  and $B$ is universal for $\mcc$ and $B\preceq_\mcb C$, then $C$ is universal for $\mcc$.  If follows that all universal sets are bi-reducible to each other.

\paragraph{The upshot} We can view the notions $\preceq_\mcb$ and $\preceq_C$ as measures of complexity.  If $A\preceq_\mcb B$, then $B$ is \emph{at least as complex} as $A$. In a given class $\mcc$ there may be {$\preceq$-maximal} elements.  They represent the equivalence class of the most complex elements of $\mcc$.

\subsection{Reducing ill-founded trees}\label{reducingillfounded}
There is a basic method of establishing a set is \emph{not} Borel. Because it arises in the proofs of several theorems we include it here rather than in the appendix \ref{DSTsum}. It gives a canonical example of a complete analytic set that is combinatorially convenient to work with.

\paragraph{The space of Trees} For applications, a common complete analytic, (and so non-Borel) set is 
the collection of ill-founded trees.
Let $\nn^{<\nn}$ be the collection of finite sequences of natural numbers, viewed as  being of the 
form $\sigma=\la n_0, n_1, \dots n_{l-1}\ra$.  The number $l$ is the length of $\sigma$.   The Polish space 
of trees is defined as follows:

	\begin{enumerate}
	\item A \emph{tree} is a set $\mct\subseteq \nn^{<\nn}$ that is closed under initial segments.  So 
	if $\sigma=\la n_0, n_1, \dots n_{l-1}\ra\in \mct$ and $k<l$ then 
	$\tau=\la n_0, n_1, \dots n_{k-1}\ra\in \mct$.\footnote{It is important to note that the trees are 
	\emph{not} required to be finitely branching. Every tree has the empty sequence in it.} 
	\item The  space $\trees$ is the collection of trees $\mct\subseteq \nn^{<\nn}$.  
	\item A basic open set in the topology on the space of trees is given by the finite sequence 
	$\sigma$ and $\la\sigma\ra=\{\mct:\sigma\in \mct\}$.
	\end{enumerate}
Recasting the definition, putting the discrete topology on $\{0, 1\}$, and the product topology on 
$\{0,1\}^{\nn^{<\nn}}$ gives a Polish space $X$ homeomorphic to the Cantor set. Each tree $\mct$ 
can be identified with its characteristic function
 $\Chi_\mct:\nn^{<\nn}\to \{0, 1\}$, and hence an 
element of $X$. It is not difficult to verify that the induced topology on 
$\{\Chi_\mct:\mct\in \trees\}\subseteq X$ gives a space homeomorphic to $\trees$

The first condition in the definition of a tree gives a countable sequence of requirements corresponding to open sets, showing that $\{\Chi_\mct:\mct\in \trees\}$ is a closed subset of $X$ and hence a Polish space.

	\begin{definition}
	Let $\mct\in \trees$ and $f:\nn\to \nn$.  Then $f$ is an \emph{infinite branch} through 
	$\mct$ if for all $l\in \nn$ the sequence $\sigma=\la f(0), f(1), \dots f(l-1)\ra\in \mct$. 
	\end{definition}
The K\"onig infinity lemma implies that every infinite, finitely branching tree has an infinite branch.  However in this context the trees in $\trees$ are allowed to be infinitely branching.

	\begin{remark*}
	In alternate language, not used in descriptive set theory, the space of trees is the space of rooted, labeled acyclic graphs (allowing countably infinite valence), and the collection of ill-founded trees are those graphs with non-trivial ends.
	\end{remark*}
	
Put the discrete topology on $\nn$ and consider $\nn^\nn$ with the product topology. This yields a Polish space called the \emph{Baire space}.  Then $B=\{(\mct, f): f\in \nn^\nn, \mct\in \trees\mbox{ and $f$ is a branch through }\mct\}$ is a closed subset of $\trees\times \nn^\nn$. Because 
\[\{\mbox{ill-founded trees}\}=\{\mct:\mbox{ for some }f\in \nn^\nn (\mct, f)\in B\}\]
it follows that the collection of ill-founded trees is an  \hyperlink{aandcoa}{analytic subset} of \trees.
\medskip

\noindent A tree is \emph{well-founded} if it is not ill-founded. Hence the well-founded trees are the complement of the ill-founded trees and thus \emph{co-analytic}.

\medskip

Fundamental results about the collection of ill-founded trees include:
	\begin{theorem}\label{treescomplete} 
	Let \trees\ be the space of trees and $I\!FT\subseteq \trees$ be the collection of ill-founded trees. Then:
		\begin{enumerate}
		\item If $X$ is an arbitrary Polish space and $A\subseteq X$ is analytic then $A$ is Borel reducible to $I\!FT$.
		\item $I\!FT$ is an analytic set that is not Borel.
		\end{enumerate}
	\end{theorem}
\begin{corollary}
Let $W\!F\subseteq \trees$ be the collection of well-founded trees.  Then 
\begin{enumerate}
\item If $X$ is an arbitrary Polish space and $C\subseteq X$ is co-analytic, then $C$ is Borel reducible to $W\!F$.
\item $W\!F$ is a co-analytic set that is not Borel.
\end{enumerate}
\end{corollary}	
	
In the language of section \ref{reductions}, the collection of ill-founded trees, as a subset of \trees, is a complete analytic set.  It also follows from Theorem \ref{treescomplete} that there \emph{are} non-Borel analytic sets.

\paragraph{Basic Technique for showing a set is analytic or co-analytic but not Borel.} 

\hypertarget{notBorel}{The basic technique for showing that a set $\mcc$ is not Borel} is reducing a known non-Borel set  to $\mcc$. To use this effectively one needs a starting point. The sets $I\!FT$ or $W\!F$ 
play this role for analytic and co-analytic sets.
\smallskip

Let $X$ be a  Polish space  and 
 $\mcc\subseteq X$ be analytic.  Then to show that $\mcc$  is complete analytic, and so not Borel, it suffices to build a reduction
\[f:\trees\to X\]
such that 
\[\mct \mbox{ is ill-founded if and only if }f(\mct)\in \mcc.\]
By Theorem \ref{treescomplete}, it follows that $\mcc$ is not Borel. 

Similarly the basic technique for showing that a co-analytic set $\mcc$ is complete co-analytic is to build a reduction
\[f:\trees\to X\]
such that 
\[\mct \mbox{ is well-founded if and only if }f(\mct)\in \mcc.\]

\paragraph{The most basic Benchmarks} The first collection of benchmarks is determining whether a classification is computable, Borel or analytic-and-not-Borel. These are the basic levels of descriptive complexity. If an example is Borel, one can also ask what level it lives in in the \hyperlink{BH}{Borel hierarchy.}

 Paraphrasing this informally, the question is whether every classification can be done using inherently finite techniques, inherently countable techniques, or whether it is impossible to do the classification using inherently countable techniques. If a classification is complete analytic or co-analytic then there is no possible classification that is feasible with countable resources.


\section{Complexity in Structure Theory}\label{CIST}

For our purposes structure theories for classes $\mcc\subset X$ are concerned with two questions:
	\begin{itemize}
	\item Can you determine whether an arbitrary $x\in X$ belongs to $\mcc$?
	\item How complicated is a given method for building elements of $\mcc$? 
	\end{itemize}

\paragraph{The complexity of the class}

The basic question is what level of descriptive complexity the class lies in relative to the basic benchmarks, given above.

\paragraph{An early example: weak mixing vs strong mixing} The first example of the use of 
descriptive set theory to solve a natural structure  problem in dynamical systems is due to Halmos and Rokhlin.   A prominent 
open problem  in the 1940's was whether the property of a 
transformation being \emph{weakly mixing} implied the property of being \emph{strongly 
mixing}. 

In 1944, in a paper called \emph{In general a measure-preserving transformation is mixing} 
(\cite{halmos1944}), Halmos proved that the collection of weakly mixing transformations is a dense 
$\mathcal G_\delta$ set in the weak topology.  In a paper published in 1948 with a title translated as 
\emph{A general measure-preserving transformation is not mixing} (\cite{rokhlin1948}), Rokhlin showed 
that the strong mixing transformations is meager (first category). The combination of the two papers show 
that there are weakly mixing transformations that are not strongly mixing. 

A proof that explicitly uses complexity shows that the collection of strongly mixing \hyperlink{BH}{transformations is not 
$\bpi^0_2$}, and hence cannot be the same class as the weakly mixing transformations. 
(In fact the class is $\bpi^0_3$.) In particular there is a weakly mixing transformation that is not strongly mixing as  it is complete $\bsigma^0_2$.

An explicit example of a weakly mixing transformation that is not strongly mixing was given using Gaussian systems was provided by Maruyama in 1949 (Theorem 11 of \cite{maruyama}). A very natural   example of such a transformation is due to Chacon (\cite{chacon}).

\subsection{Examples at three levels of complexity}

We now give examples of each level beyond recursively computable.

\paragraph{Non-recursive classes}
	\begin{example*}(Braverman, Yampolsky) \cite{BY} There is a computable complex number $s$ such that 
	membership in the  Julia set associated with 
	\[f(z)=z^2+s\]
	is not computable.
	\end{example*}

\subsection{Borel Classes}
\paragraph{Borel Classes of measure-preserving transformations}
As mentioned above, the weakly mixing and strongly mixing classes of measure-preserving transformations of the unit interval are Borel.  The multiple mixing counterparts are as well. 

In the topological category, the collections of uniformly continuous maps from the open unit interval to the open unit interval are easily seen to be $\bpi^0_2$ in the space of continuous transformations. 

Indeed most standard classes of transformations are Borel. As explained in the \hyperlink{firstbernoulli}{discussion above} the collection of Bernoulli shifts is Borel. However, proving this  uses \hyperlink{suslintheorem}{Suslin's Theorem} which involves analytic sets. 
\medskip

Other well-known  classes of measure-preserving transformations are Borel and the proof uses a significant amount of 
work. For example the class of rank-one-rigid transformations. (\cite{king})

\paragraph{Borel Classes of homeomorphisms of compact spaces}
 A potential problem with this collection is that there are a variety of natural spaces in which to set the various questions.  For example every compact $\poZ$-flow can be viewed as a compact subset of the $\mathbb H^\poZ$ or as a compact subset of $\mathbb H^2$ (where $\mathbb H$ is the Hilbert cube) or as an appropriate subset of $\poN^\poN$. All of the known natural codings of the collection of homeomorphisms of compact metric spaces agree on what collections of homeomorphisms are Borel. (This is discussed in \cite{BF1}.)

Examples include:
\begin{itemize}
	\item  $\{(X,T):T$ is topologically transitive and $X$ is Polish $\}$, 
	\item $\{(X, T):T$ is minimal and $X$ is Polish$\}$. 
\end{itemize}
Proofs of these results (and later calculations of an upper bound on the complexity of topologically distal transformations) were novel in the early 1990's and appear in \cite{BF1} (Lemma 3.3).  The second item is due to Kechris.

\paragraph{Borel Classes of diffeomeomorphisms} We have already seen that the collection of structurally stable diffeomorphisms is Borel--in fact it is open as are the collections of Anosov and Morse-Smale diffeomorphisms.
\medskip

\noindent To give a flavor of tools used for more complicated classes,  we use the following results in \cite{FoGo}:
\smallskip

Denote by $\text{Fix}(f)$ the set of fixed points of the map $f:M\to M$, i.e. $$\text{Fix}(f)=\{x\in M\ |\ f(x)=x\},$$ and by $\text{Per}_n(f)$ the set of periodic points of period $n$, i.e. 
\[\text{Per}_n(f)=\{x\in M\ |\ f^n(x)=x \mbox{ and for all } m<n, f^m(x)\ne x\}.\]

\begin{prop*}\label{p.pp} Let $M$ be a compact manifold.
For any $r=0, 1, 2, \ldots, \infty$ the following statements hold:
	\begin{enumerate}
	\item The set $\mathcal{U}_0=\{f\in \text{Diff}^{\,\, r}(M):\text{Fix}\,(f)=\emptyset\}$ is open;

	\item For each $m\in \mathbb{N}$ the set $\mathcal{U}_m=\{f\in \text{Diff}^{\,\, r}(M):|\text{Fix}\,(f)|=m\}$ is Borel; 

	\item For each $m, n\in \mathbb{N}$ the set $\mathcal{U}_m^n=\{f\in \text{Diff}^{\,\, r}(M):|\text{Per}_n(f)|=m\}$ is Borel;

    \item  For any $A:\nn\to\nn$  the set 
    \[\mathcal{U}_A=\{f\in \text{Diff}^{\,\, r}(M):|\text{Per}_n(f)|=A(n), \ n=0, 1, 2, \ldots \}\]
     is Borel;

    \item The set 
    \[\{(A, f): |\text{Per}_n(f)|=A(n), \ n=0, 1, 2, \ldots,\ \text{where}\ A:\nn\to \nn  \}\] in the space $\mathbb{N}^{\mathbb{N}}\times \text{Diff}^{\,\, r}(M)$ is Borel;

\end{enumerate}
    
\end{prop*}
We give the proof that appeared in  \cite{FoGo} for items 1-3. The first item is immediate: since $M$ is compact there is a minimum distance between $x$ and $f(x)$.  The set of diffeomorphisms whose minimum distance is bigger than $1/n$ is open, hence the set of diffeomorphisms with no fixed points is open. 

To see the second item: To see that $\mathcal{U}_1$ is Borel, let us choose a countable base for the topology on $M$, $\{V_k\}_{k\in \mathbb{N}}$, such that

\[\text{diam}\, V_1\ge \text{diam}\, V_2\ge \text{diam}\, V_3\ge \ldots\] 
{and} $\text{diam}\, V_k\to 0\ \ \text{as}\ k\to \infty.$
The set
\[
\mathcal{V}_k=\{f\in \text{\it Diff}^{\,\, r}(M)\ |\ f\ \ \text{has no fixed points outside of}\  V_k\}
\]
is open, and hence the set
$
\mathcal{U}_1=\left(\bigcap_{K\ge 1} \bigcup_{k=K}^\infty \mathcal{V}_k\right)\backslash \mathcal{U}_0
$
is Borel. 

To see that $\mathcal{U}_2$ is Borel, notice that each of the sets
\[\mathcal{V}_{k_1, k_2}=\left\{f\in \text{\it Diff}^{\,\, r}(M)\ |\ f\ \ \text{has no fixed points outside of}\  V_{k_1}\cup V_{k_2}\right\}.
\]
is open, and hence
\[
\mathcal{U}_2=\left(\bigcap_{K\ge 1} \bigcup_{k_1, k_2\ge K} \mathcal{V}_{k_1, k_2}\right)\backslash \left(\mathcal{U}_0 \cup \mathcal{U}_1\right)
\]
is Borel. It is clear how to continue for $n$ fixed points by induction.

\begin{example*}{\bf Artin-Mazur Diffeomorphisms}
It follows easily from the proposition that  the function giving the exponential rate of growth of the number of periodic points
$G:\text{Diff}^{\,\, r}(M) \to \mathbb{R}\cup\{\pm\infty\}$ defined by
    \[
    G(f)=\limsup_{n\to \infty}\frac{1}{n}\log |\text{Per}_n(f)|
    \]
 is Borel. In particular, the set of  \emph{Artin-Mazur} diffeomorphisms (defined as 
 $\{f:G(f)<\infty\}$) is Borel.
\end{example*}
	
\subsection{Natural Classes that are not Borel} \label{intodistal}

\paragraph{In the topological setting}
	A classical collection, due to Hilbert (see \cite{zippin}) is defined as follows. 
	\begin{definition}\label{topodistal}
	Let $(X,d)$ be a metric space and $T:X\to X$ be a homeomorphism.  Then $(X, T)$ is (topologically) \emph{distal} if
	for all distinct $x, y\in X$, $\inf_{n\in \poZ}\{d(T^nx, T^ny)\}>0$.
	\end{definition}

Clearly isometries are examples of distal transformations.  Moreover the distal transformations are closed under taking skew products with compact metric groups. Using these two facts one can generate a wide class of distal transformations.

\begin{fact*} The collection of distal transformations is a co-analytic set.  This is exposited in example \ref{distaliscoanalytic}.
\end{fact*}

\begin{definition}
Let $(X, d, T)$ and $(Y,d',S)$ be a compact separable metric spaces and $\pi:X\to Y$ be a factor map.  Then $X$ is an \emph{isometric} extension of $Y$ if there is a real valued function $\rho(x_1,x_2)$ defined for all pairs $(x_1,x_2)$ in
the same fibre of $X$ (i.e. whenever $\pi(x_1)=\pi(x_2)$) and such that
\begin{enumerate}
\item {$\rho(x_1,x_2)$ is continuous as a function on the subset of
$X\times X$ defined by the condition $\pi(x_1)=\pi(x_2)$.}

\item {for each $y\in Y$, $\rho(x_1,x_2)$ defines a metric on the fibre
$\pi^{-1}(y)$ under which this fibre is isometric to a fixed homogeneous metric
space.}

\item {$\rho(T x_1,T x_2)=\rho(x_1,x_2)$ for all $x_1, x_2$ in the same fibre of $X$ over $Y$.}
\end{enumerate}
\end{definition}

The next definition captures the process of iterating the isometric extensions transfinitely:
\begin{definition}
$(X,T)$ is a quasi-isometric extension of $(Y,T)$ if $(Y,T)$ is a factor of
$(X,T)$, and there is an ordinal $\eta$ and a factor $(X_{\xi},T)$ of $(X,T)$
for each $\xi\le \eta$ such that:
\begin{enumerate} 

\item {$(X_0,T)=(Y,T)$; $(X_{\eta},T)=(X,T)$}

\item {if $\xi<\xi'$ then $(X_{\xi},T)$ is a factor of $(X_{\xi'},T)$}

\item {$(X_{\xi+1},T)$ is an isometric extension of $(X_{\xi},T)$}

\item {if $\xi$ is a limit ordinal then $(X_{\xi},T)$ is the inverse limit of
$\langle (X_{\xi'},T):\xi'<\xi\rangle$.}
\end{enumerate}
\end{definition}

\noindent The Furstenberg Structure Theorem (\cite{topodistal}) for topologically distal flows says:

\begin{theorem}(Furstenberg)
Let $(X, d, T)$ be a minimal, topologically distal flow.  Then  $(X, d, T)$ is a quasi-isometric extension of the trivial flow. 
\end{theorem}

It is natural to ask what ordinals arise as $\eta$ in the Furstenberg Structure Theorem. Given a distal flow $(X,d,T)$ 
define the \emph{distal height} of the flow as the least ordinal $\eta$ such that $(X, d, T)$ is topologically equivalent to 
a quasi-isometric extension of  the trivial flow represented as a tower of extensions of height $\eta$.

It was shown in \cite{BF1} that the collection of distal heights of minimal distal flows on compact metric spaces is exactly the collection of countable ordinals. 
The following theorem (\cite{BF1}) summarizes these results:

	\begin{theorem}(Beleznay-Foreman) Let $\mathcal M$ be the Polish space of  minimal flows on 
	compact metric spaces. Then:
		\begin{enumerate}
		\item The collection of distal flows is a complete co-analytic set, in particular it is non-Borel.
		\item The natural function that associates to each distal  $(X, d, T)$ its distal height is a 
		$\bpi^1_1$-norm.
		\end{enumerate}
	\end{theorem}
	The first item is proved by showing that there is a Borel reduction from the \hyperlink{wo}{set of codes for well-orderings, $WO$} to the set of  
distal flows as a subset of the space of minimal $\poZ$-flows.

If $(X_0, d_0, T_0), (X_1, d_1, T_2)$ are minimal distal flows on compact spaces, set\\
$(X_0, d_0, T_0)\le_D (X_1, d_1, T_2)$ if the distal height of $(X_0, d_0, T_0)$ is less than or equal to the 
distal height of $(X_1, d_1, T_2)$. Then the theorem shows that  \hyperlink{pi11norm}{$\le_D$ is a true 
$\bpi^1_1$-norm} on the collection of distal flows. 
It follows that for all countable ordinals $\eta$ there is a distal flow of distal height $\eta$ and that  the
 collection of minimal distal homeomorphisms of distal height less than $\eta$ is a Borel set. 
 
\paragraph{In the setting of measure-preserving transformations} Separate work of Furstenberg and Zimmer (\cite{FUbook}, \cite{Zim}) produced analogues of the topological Furstenberg Structure Theorem that hold for measure-preserving transformations. We describe it here without giving details and with very slight inaccuracies.

The category we are working in consists of measure-preserving transformations of a standard measure space, which we can take to be the unit interval with Lebesgue measure. By \hyperlink{ergisGd}{Halmos' theorem, the ergodic diffeomorphisms} form a dense 
$\mcg_\delta$ set. Hence the weak topology makes the ergodic transformations a Polish space. We take our ambient Polish space $EMPT$ to be the set of ergodic measure-preserving transformations of the unit interval.

Let $\ycns$ be an ergodic measure-preserving transformation, and $G$ a compact group.  Let $H$ be a compact subgroup of $G$.  Then the  quotient space $G/H$ carries a left-invariant Haar measure. A compact co-cycle extension 
\[S_\phi:Y\times G/H\to Y\times G/H\] of $Y$ is determined by a measurable function $\phi:Y\to G/H$ and given by the formula:
\[S_\phi(y,[g]_H)=(S(y),[\phi(y)g]_H).\]
Then $S_\phi$ preserves the product of $\nu$ and Haar measure on $G/H$.

	\begin{definition}(See e.g. \cite{FUbook}) Let $\xbmt$ and $\ycns$ be ergodic measure-preserving systems. Then
		\begin{enumerate} 
		\item $X$ is a 	
			\emph{compact} extension of $Y$ if $X$ is measure equivalent to a compact co-cycle 
			extension of $Y$.
		\item Let $\pi:X\to Y$ be a measure-preserving extension. Let $X\times_Y X$  be the measure space which has 
		domain $\{(x_1, x_2):\pi(x_1)=\pi(x_2)\}$  and the canonical measure projecting to $\nu$ that makes 
		the two copies of $X$ relatively independent over $(Y, \mcc, \nu)$.   Let $T$ be  the 
			measure-preserving system taking $(x_0, x_1)$ to  $(T(x_0), T(x_1))$.  Then $Y$ canonically 
			a factor of $X\times_Y X$ and 
			 $X$ is a \emph{weakly mixing} extension of $Y$ if  $X\times_Y X$ is ergodic.
		\end{enumerate}
	\end{definition}
The theorem proved independently by Furstenberg and Zimmer is:
	\begin{theorem}
	Let $\xbmt$ be an ergodic measure-preserving system.  Then there is a countable ordinal $\eta$ and a system of measure-preserving transformations $\la (X_\alpha, \mcb_\alpha, \mu_\alpha, T_\alpha):\alpha\le \eta+1\ra$ such that
		\begin{enumerate}
		\item $(X_0, \mcb_0, \mu_0, T_0)$ is the trivial flow.
		\item For each $\alpha<\eta$, $X_{\alpha+1}$ is a compact extension of $X_\alpha$.
		\item If $\alpha $ is a limit ordinal then $X_\alpha$ is the inverse limit of $\la X_\beta:\beta<\alpha\ra$
		\item $X_{\eta+1}$ is either:
			\begin{itemize}
			\item[-] a trivial extension of $X_\eta$ (so $X_{\eta+1}\cong X_\eta$), or
			\item[-] a weakly-mixing extension of $X_\eta$. 
			\end{itemize}		
		\end{enumerate}
	\end{theorem}
\noindent Note that in each case, for $\alpha=1$, the transformation $X_\alpha$  has to be a translation on a compact quotient $G/H$.  In this case $G$ must be abelian so $G/H$ is also an abelian group.  Hence monothetic and an example of the form due to  \hyperlink{HvN}{Halmos and von Neumann.}

In analogy to the situation for topologically distal transformations, we give the following definition: 

	\begin{definition}\label{measdist}
	An ergodic measure-preserving system $\xbmt$ is \emph{measure distal} if it can be written in the form of the 
	previous theorem with $X_{\eta+1}=X_\eta$.
	\end{definition}

\noindent As in the topological case we define the \emph{distal height} of a measure distal transformation $\xbmt$ to be the least 
$\eta$ such that $\xbmt$ can be written satisfying 1-4 of the theorem.

\medskip
The results are now analogous to the topological case (see \cite{FB2}): 

\begin{theorem}(Beleznay-Foreman) Let $EMPT([0,1])$ be the Polish space of  ergodic measure-preserving transformations of the unit interval. Then:
		\begin{enumerate}
		\item The collection of measure distal flows is a complete co-analytic set, in particular it is non-Borel.
		\item The natural function that associates to each measure distal  $\xbmt$ its distal \hyperlink{pi11norm}{height is a $\bpi^1_1$-norm.}
		\end{enumerate}
	\end{theorem}
It follows from the second item  that if $\eta$ is a fixed countable ordinal then the collection of measure distal homeomorphisms of distal height less than $\eta$ is a Borel set.

\paragraph{An open problem in both settings}
The following problems are analogous and involve \hyperlink{overspill}{\emph{overspill}}. (See Open Problem \ref{OP4.5})
\begin{itemize}
\item Let $X$ be the collection of minimal homeomorphisms on compact metric spaces and $\mathcal D$ be the collection of topologically distal transformations.  Is there a Borel set $B\subseteq X\times X$ such that $B\cap (\mathcal D\times \mathcal D)$ is the relation of topological conjugacy?
\item Let $\mathcal{ MD}$ be the collection of measure distal transformations.  Is there a Borel set $B\subseteq EMPT\times EMPT$ such that $B\cap (\mathcal{ MD}\times\mathcal{MD})$ is the relation of measure conjugacy?
\end{itemize}
In unpublished work, Foreman showed that for the collection of  measure distal transformations with height \emph{bounded by a countable ordinal} 
$\alpha$, the answer to the second question is \emph{yes}.

\section{Complexity in Classification Theory}

\begin{figure}[h]
\centering
\includegraphics[height=.6\textheight]{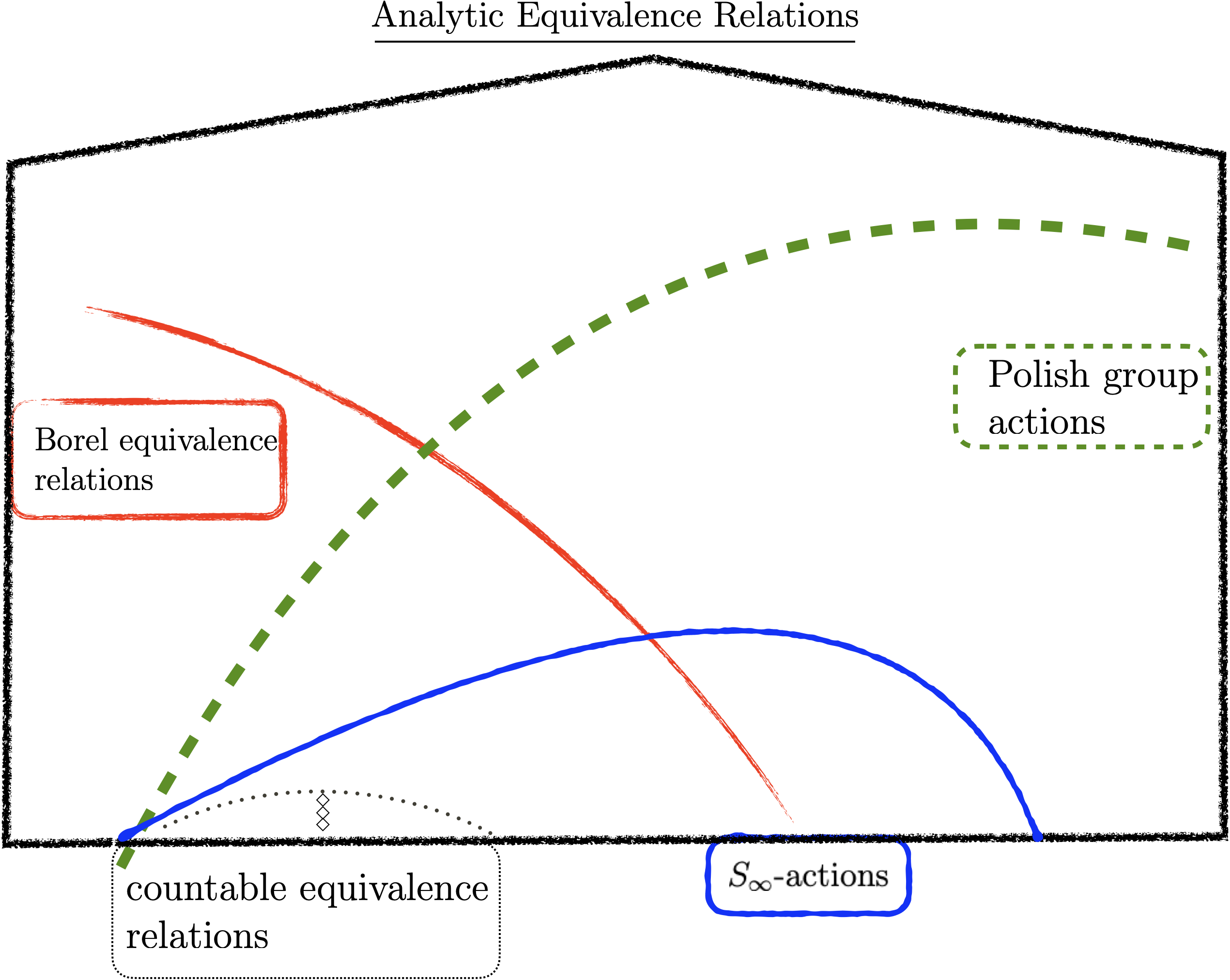}
\caption{Basic regions of complexity}
\label{yayayaya}
\end{figure}

In this section we list various benchmarks for complexity of classifications of equivalence relations $E$ that 
have roots outside of dynamical systems. The general regions used to classify equivalence relations  
are shown in figure \ref{yayayaya}. The most sensitive measures of complexity are the two dimensional 
questions about reductions \emph{between equivalence relations}.

\paragraph{A basic one-dimensional benchmark}
If $X$ is a Polish space  then $X\times X$ is also Polish. For any equivalence relation 
$E\subseteq X\times X$, one can ask the one-parameter question:  
	\begin{itemize}
	\item Is $E$ Borel as a subset of $X\times X$?
	\end{itemize}
If not, then determining whether a pair $(x_0, x_1)\in X\times X$ is in the relation $E$ cannot be accomplished with inherently countable resources.

\bigskip

Many important equivalence relations arise out of group actions. If $G$ acts on $X$  then the equivalence relation $x\sim y$ if and only if there is a $g\in G$ such that $gx=y$ is called the orbit equivalence relation. 

A \emph{Polish group action} is an action of a Polish group $G$ on a Polish space $X$ that is jointly continuous in each variable. Section \ref{S-inf} is concerned with actions of the group of permutations of the natural numbers, $S_\infty$. When we say \emph{$G$ acts on $X$} where $G$ and $X$ are both Polish, we will always be assuming that the action is a Polish group action.

\begin{remark*}
When we introduce a class of equivalence relations (such as the $S_\infty$-actions defined below), we really mean to take $\preceq^2_\mcb$ downwards closure.  So the class we label ``$S_\infty$-actions" means any equivalence relation reducible to an $S_\infty$-action.
\end{remark*}

\subsection{Equivalence relations with countable classes}\label{countable classes}
The equivalence relations that have finite or  countable classes are, perhaps, the simplest equivalence relations. 
There is an extensive literature on this, for example see \cite{kechcount}, and \cite{kechmiller}. Equivalence relations with countable classes are referred to as \emph{countable equivalence relations.}

The orbit equivalence relation of any countable group action has countable classes, so there is a plethora of examples.

Let $\la \phi_n:n\in\nn\ra$ be a collection of Borel bijections of a Polish space to itself.  Let $\Gamma$ be the group generated by $\la \phi_n:n\in\nn\ra$. Then each orbit of $\Gamma$ is countable and the orbit equivalence relation is Borel.  Conversely,  a theorem of Feldman and Moore (\cite{FeldmanMoore}) states:

\begin{theorem}(Feldman-Moore)
Suppose that $E\subseteq X\times X$ is a countable Borel equivalence relation on a Polish space $X$. Then there is a countable group $\Gamma$ and a Borel action of $\Gamma$ on $X$ such that $E$ is the orbit equivalence relation of the action.
\end{theorem}

\paragraph{A $\preceq^2_B$ maximal countable Borel equivalence relation} Let $G$ be a countable discrete group and 
$\phi:G\times X\to X$ be a Borel action on a Polish space $X$.  Then the orbit equivalence relation is a countable Borel equivalence relation on $X$.

Following \cite{kechcount} we let $s_G$ be the shift 
action of $G$ on the space $X^G$, consisting of functions $f:G\to X$. Let $E(G,X)$ be the orbit equivalence relation of 
this action. Dougherty, Jackson and Kechris (\cite{DJKhyperfinite}) showed that if $G$ is the free group on two 
generators, $F_2$, then $E(F_2,X)$ is a maximal Borel countable equivalence relation. The following is immediate from the Dougherty-Jackson-Kechris theorem.

\begin{theorem} \hypertarget{einfty}{Let $E_\infty$ be any countable equivalence relation \hyperlink{bired}{bi-reducible to $E(F_2, X)$} for any Polish space $X$.}  Then every countable Borel equivalence relation is reducible to $E_\infty$.
\end{theorem}

We show more about countable equivalence relations in section \ref{countable classes}.

 \subsection{The Glimm-Effros Dichotomy: $E_0$}\label{GE}
We now turn to equivalence relations with uncountable  classes. 
  Because complete numerical invariants are so useful, they are the baseline benchmark of complexity. If $Y$ is a Polish space we let $=_Y\ \subseteq Y\times Y$ be the  equivalence relation of equality on $Y$. 
   
  If $Y$ is a Polish space then there is a Borel injection $g:Y\to \mathbb R$.
(See \cite{naivedst}.). Let $f$ be a Borel reduction of any equivalence relation $E\subseteq X\times X$ to $(Y,=_Y)$. Composing $f$ with $g$ we see that $g\circ f$ is a   reduction of $(X, E)$ to $(\mathbb R, =_{\mathbb R})$. Moreover $g\circ f$ is Borel  if $f$ is. Thus we can assume that Borel reductions to any $=_Y$ can be changed to Borel reductions to equality on the real numbers.
  
 \begin{definition}\label{smoothdef} Let $E$ be an equivalence relation on a Polish space $X$. Then $E$ is \emph{smooth} if 
 $E\preceq_\mcb^2 \ =_{\mathbb R}$. 
\end{definition}  
 
In  1990, Harrington, Kechris and Louveau (\cite{HKL}), extending earlier work of Glimm \cite{glimm} and Effros \cite{effros}, showed a very general dichotomy theorem.

\begin{definition}\label{E0} \hypertarget{ezero}{Let $E_0$ be the equivalence relation on the Cantor set 
$\{0,1\}^\nn$ defined by setting $fE_0\ \!g$ if and only if there is an $N$ for all $m>N, f(m)=g(m)$.}
\end{definition}

\begin{proposition}
 The equivalence relation $=$ on $\{0,1\}^\nn$ is continuously reducible to $E_0$.
 \end{proposition}
 
 \pf  Fix a bijection $\la \cdot, \cdot\ra:\nn\times \nn\to \nn$.  Define $R:\{0,1\}^\nn\to \nn$ by setting\footnote{This simple, elegant proof was taken from notes from a lecture of Su Gao. See citation in the bibliography.}
  \[R(f)(\la m, n\ra)=f(n).\]
Since $R$ is well-defined, $f=g$ implies $R(f)E_0R(g)$.  On the other hand suppose for all $m>N, R(f)(m)=R(g)(m)$.  Let $n\in \nn$ and choose an $M$ such that $\la M, n\ra>N$.  Then 
	\begin{eqnarray*}
	f(n)&=&R(f)(\la M, n\ra)\\
	&=&R(g)(\la M, n\ra)\\
	&=&g(n).
	\end{eqnarray*}
 \qed

The Harrington-Kechris-Louveau dichotomy is the following theorem.
\begin{theorem} \label{HKL}(\cite{HKL}) Let $X$ be a Polish space and $E$ a Borel equivalence relation on $X$. Then either:
	\begin{enumerate}
	\item $E$ is \emph{smooth}	or
	\item $E_0\preceq^2_\mcb E$.	
	\end{enumerate}
\end{theorem}
\begin{corollary}
$=\ \preceq^2_\mcb E_0$ but $E_0\not\preceq^2_\mcb\ =$. Thus $E_0$ is strictly more complicated than the identity equivalence relation.
\end{corollary}

This theorem suggests an orderly stair-step pattern to the benchmarks, but this is misleading, as we shall see.
\smallskip

\noindent We prove the easy direction of Theorem \ref{HKL}, because it is most useful in dynamical systems.\footnote{
 The proof of Theorem \ref{E0doom} in an earlier draft of this paper used zero one laws for sets with the Property of Baire. The author would like to thank the referee for suggesting the simpler proof given here.}

	\begin{enumerate}
	\item Put a measure on 
 $\{0,1\}$ by weighting each element 1/2.  The product measure $\mu$ on $\{0,1\}^{\nn}$ is often called the 
 \emph{coin-flipping measure}.  
 	\item Let $G=\Sigma_{n\in \nn}\poZ_2$ be the direct sum of infinitely many copies of 
 $\poZ_2$ (with finite support).   	
	 \item The coordinatewise action of $G$ on $(\{0, 1\}^{\nn}, \mcb, \mu)$ preserves 
	 $\mu$.  
 	\item A standard $0-1$ law for this action says that if $A\subseteq \{0, 1\}^{\nn}$ is invariant under the $G$-action then $A$ either has $\mu$-measure $0$ or $\mu$-measure $1$. 
	 \end{enumerate}

\begin{theorem}\label{E0doom}
Suppose that $E$ is an equivalence relation on an uncountable Polish space $X$ and $E_0\preceq^2_\mcb E$. Then $E$ is not smooth.
\end{theorem}

 \pf  Towards a contradiction, assume that $E$ is Borel reducible to the identity equivalence relation on $\mathbb R$ by a map $S:X\to \mathbb R$. Let $R:\{0, 1\}^\nn\to X$ be a reduction of $E_0$ to $E$.   By composing we see that $S\circ R$ is a reduction of $E_0$ to the identity relation on $[0,1]$. It suffices to show that there can be no such reduction from $E_0$ to the identity relation $=$ on $[0,1]$.
 
 Suppose that $S$ is such a reduction. Because $G$ preserves $E_0$, for each rational number $r\in \mathbb R$,  $L_r=\{f:S(f)<r\}$ is $G$ invariant, and so has measure zero or one.  Let 
 \[x=\sup\{r:L_r\mbox{ has measure one}\}.\]
 Then $I=\{f:S(f)=x\}$ has measure one.  Hence there are $f, g\in I$ such that $f\not\!\! E_0 g$, a contradiction. \qed

\subsection{$=^+$ and the Friedman-Stanley jump operator} \label{fsbs}
\hypertarget{sinf}{Let $S_\infty$ be the group of permutations of the natural numbers.}\footnote{Some authors reserve  the notation $S_\infty$ for the group of permutations of $\nn$ that move only finitely many points and use $S^\infty$ for the whole group of permutations. The convention in this paper is different and  is becoming universal.}  Then $S_\infty$ is a $\mcg_\delta$ subset of $\nn^\nn$ and so carries a Polish topology.  This topology is compatible with the group structure and so we view $S_\infty$ as a Polish group.
\medskip

The following operation on analytic equivalence relations was defined by Friedman and Stanley 
(\cite{Fried}).
\smallskip

Let $E\subseteq X\times X$ be an analytic equivalence relation on a Polish space. Define $E^+$ to be the equivalence on $X^\nn$ given by setting $\la [x_n]_E:n\in \nn\ra$ $E^+$-equivalent to $\la [y_n]_E:n\in\nn\ra$ if and only if
	\begin{quotation} \noindent for all $n$ there is an $m$ $x_nEy_m$ and for all $m$ 	
	there is an $n$ $y_mEx_n$.
	\end{quotation}

\medskip
Because the collection of Borel sets (and respectively analytic sets) are closed under countable intersections and unions, the form of the definition 
makes it clear that if $E$ is Borel, then $E^+$ is Borel (and similarly for $E$ being analytic).  If $X$ and $Y$ are perfect Polish spaces then $=_X^+$ and $=_Y^+$ are Borel bi-reducible. Moreover starting with the identity relation $=$ 
iterating the ``+" operation any countable ordinal  number $\alpha$ of times leads to a sequence of equivalence 
$\la E_\beta:\beta<\alpha\ra$ such that $\beta<\beta'$ implies $E_\beta\preceq_\mcb E_{\beta'}$ but $E_{\beta'}\not{\!\! \preceq_\mcb} E_\beta$ (\cite{Fried}). Thus this is a strict hierarchy of height $\omega_1$.

If $E$ is Borel then $\{\vec{x}\in X^\nn: [x_n]_E\ne [x_m]_E$ for all $n\ne m\}$ is a Borel set.  Restricting the ``$+$-operation" to this set at each stage also gives a proper hierarchy in $\preceq^2_\mcb$.

The following result is proved with a similar approach to remark \ref{E0doom}. See \cite{GaoLecture} for a simple elegant proof of the first reduction.
\begin{theorem}
The following relations hold: 
\[E_0\preceq^2_\mcb\ =^+\mbox{ and \ } {=^+ \not\preceq^2_\mcb}\ E_0.\]
  In particular, 
$=^+\not{\!\! \preceq^2_\mcb} =$.
\end{theorem}
\subsection{$S_\infty$-actions}\label{S-inf}

We now discuss $S_\infty$-actions, where $S_\infty$ is the Polish group of all permutations of $\nn$ that was defined  in the previous section.
The $S_\infty$-actions are a fundamental benchmark because their actions are ubiquitous:  every non-Archimedean Polish group is isomorphic to a closed subgroup of $S_\infty$. (See \cite{BK}.)
\medskip

Fix an infinite  perfect Polish space $X$ and consider $P\subseteq X^\nn$ defined as the collection of infinite sequences $\la x_n:n\in\nn\ra\in X^\nn$ such that
\[\mbox{for all $n\ne m, x_n\ne x_m$}.\] 
Then $P$ is a perfect closed subset of $X^\nn$, $S_\infty$ acts on $P$ coordinatewise and for 
$\vec{x}, \vec{y}\in P$ we have 
\[\vec{x}=^+\vec{y} \mbox{ if  and only if  there is a }\phi\in S_\infty, \phi\vec{x}=\vec{y}.\]
Thus the equivalence relation on $P$ induced by $=^+$ is the orbit equivalence relation of an $S_\infty$-action.

\begin{remark}\label{borelSinfty}
The equivalence relation $=^+$ restricted to $P$ is a Borel equivalence relation. Thus there are faithful $S_\infty$-actions inducing Borel equivalence relations.
\end{remark}
Remark \ref{borelSinfty} is immediate: for $\vec{x}$ and $\vec{y}$ in $P$, $\vec{x}=^+\vec{y}$ if and only if for all $n$ there is an $m$ with $x_n=y_m$ and for all $m$ there is an $n$ with $y_m=x_n$.
\smallskip

\paragraph{Notation} The only example of the use of the ``$+$-operation" in this paper will be when $E$ is the identity equivalence relation, $=$, on the unit circle. 
We will restrict $=^+$ to one-to-one sequences. In an abuse of notation we will continue to refer to this relation as $=^+$. 

\medskip

\paragraph{Classification by countable structures} We now give an example to illustrate the importance of $S_\infty$-actions for creating invariants for equivalence relations.
\smallskip
 
Let $(X, E)$ be an analytic equivalence relation that has complete invariants consisting of countable 
groups. This means that a countable group is associated to each element of $X$ by a Borel function  in  
such a way that $xEy$ just in case the group associated with $x$ is isomorphic to the group associated 
with $y$. 
In the language we are using here, the equivalence relation $E$ is being reduced to the equivalence 
relation on pairs of countable groups given by  isomorphism. The next example makes this precise. 
\begin{example*} Countable groups are determined by the characteristic function of their multiplication. 
Explicitly, for $G$ an infinite countable group we can assume that its domain is $\nn$. Define 
\[\Chi_G(l,m,n)= \begin{cases}
1 &\mbox{ if }l*_Gm=n\\
0 & \mbox{otherwise}
\end{cases}\]
Then each $\Chi_G\in \{0, 1\}^{\nn^3}$ and $\{\Chi_G: G$ is a countable group$\}$ is a $\mcg_\delta$ subset of $\{0, 1\}^{\nn^3}$, and hence form a Polish space. Call it $CG$. 
\smallskip

Let $S_\infty$ act on $CG$ by setting $(\phi\cdot\Chi_G)(l, m, n)=\Chi_G(\phi(l), \phi(m), \phi(n))$.   If $G, H$ are 
isomorphic, let $\phi:G\to H$ be the isomorphism.  Then $\phi$ takes the multiplication table of $G$ to the multiplication table of $H$. Hence  
$\phi\cdot\Chi_G=\Chi_H$. 
Similarly if $\phi\cdot\Chi_G=\Chi_H$, then $\phi:G\to H$ is an isomorphism.
\end{example*}

This example is clearly not specific to groups--it can be adapted to any countable structures.  In particular it applies to countable graphs viewed as elements of $\{0, 1\}^{\nn^2}$.  Let $E_{groups}$ be the equivalence relation of isomorphism of countable groups and $E_{graphs}$ be the equivalence relation of isomorphism of countable graphs. 

\begin{theorem}\label{thetwogrs}
Let $S_\infty$ act on an arbitrary Polish space $X$ and $F$ be the resulting orbit equivalence relation. Then:
	\begin{enumerate}
	\item (Mekler, \cite{Mekler}) $F\preceq^2_\mcb E_{groups}$
	\item (Friedman, Stanley, see \cite{Gaobook}) $F\preceq_\mcb^2 E_{graphs}$
	\end{enumerate}
Thus both $E_{groups}$ and $E_{graphs}$ are $\preceq_\mcb^2$-maximal among $S_\infty$-actions.
\end{theorem}

These examples illustrate why equivalence relations that are reducible to $S_\infty$-actions are referred to as those that can be \emph{classified by countable structures}, the name for them given in \cite{Hj}. We note that there is varying terminology for equivalence relations that are  $\preceq^2_\mcb$ bi-reducible to maximal $S_\infty$-actions. They are called \emph{maximal} or \emph{universal} $S_\infty$-actions and also \emph{Borel Complete} equivalence relations.

To illustrate the large collection  of equivalence relations that arise from $S_\infty$-actions and are 
$\preceq^2_\mcb$ maximal among $S_\infty$-actions we cite results of  Camerlo and Gao (\cite{CG}).

Recall that \emph{almost finite} commutative $C^*$-algebras (commutative``AF" $C^*$-algebras) are those that are inverse limits of finite dimensional commutative $C^*$ algebras.  They have a classification in terms of \emph{Bratelli diagrams}. (See \cite{CG}.) 

\begin{theorem}(Camerlo-Gao)\label{camgao}
The isomorphism relation among commutative AF $C^*$-algebras is  $\preceq^2_\mcb$-maximal among $S_\infty$-actions.
\end{theorem}

\paragraph{A maximal $S_\infty$-action  is NOT Borel}
One property that makes equivalence relations such as $E_{graphs}$ and $E_{groups}$ very useful is that maximal $S_\infty$ relations are not Borel.  We record this here:

\begin{fact}\label{SinftynotBorel} $S_\infty$-actions that are {$\preceq^2_\mcb$}-maximal have complete analytic orbit equivalence 	relations. Hence they are not Borel.  
	\end{fact}
In particular, \emph{isomorphism of countable groups}, \emph{isomorphism of countable graphs} and \emph{isomorphism of AF $C^*$-algebras} are not Borel equivalence relations.

For some applications in dynamical systems the particular case of ``isomorphism of countable graphs" is a particularly convenient example of an $\preceq^2_\mcb$-maximal $S_\infty$ orbit equivalence relation.

\paragraph{A caution} We do note that not all non-trivial  ``classifications by algebraic structures" give equivalence relations that are $\preceq^2_\mcb$-maximal $S_\infty$ actions. An example, due to Friedman and Stanley (\cite{Fried}) is the equivalence relation of \emph{isomorphism of countable abelian torsion groups.}  This is a complete analytic equivalence relation that determined by an $S_\infty$-action, but is not $\preceq^2_\mcb$-maximal among such actions.

\paragraph{Benchmarks to this point}
Here is how these relations are related:

\pagebreak
\begin{figure}[!h]
	\centering
	\includegraphics[height=.50\textheight]{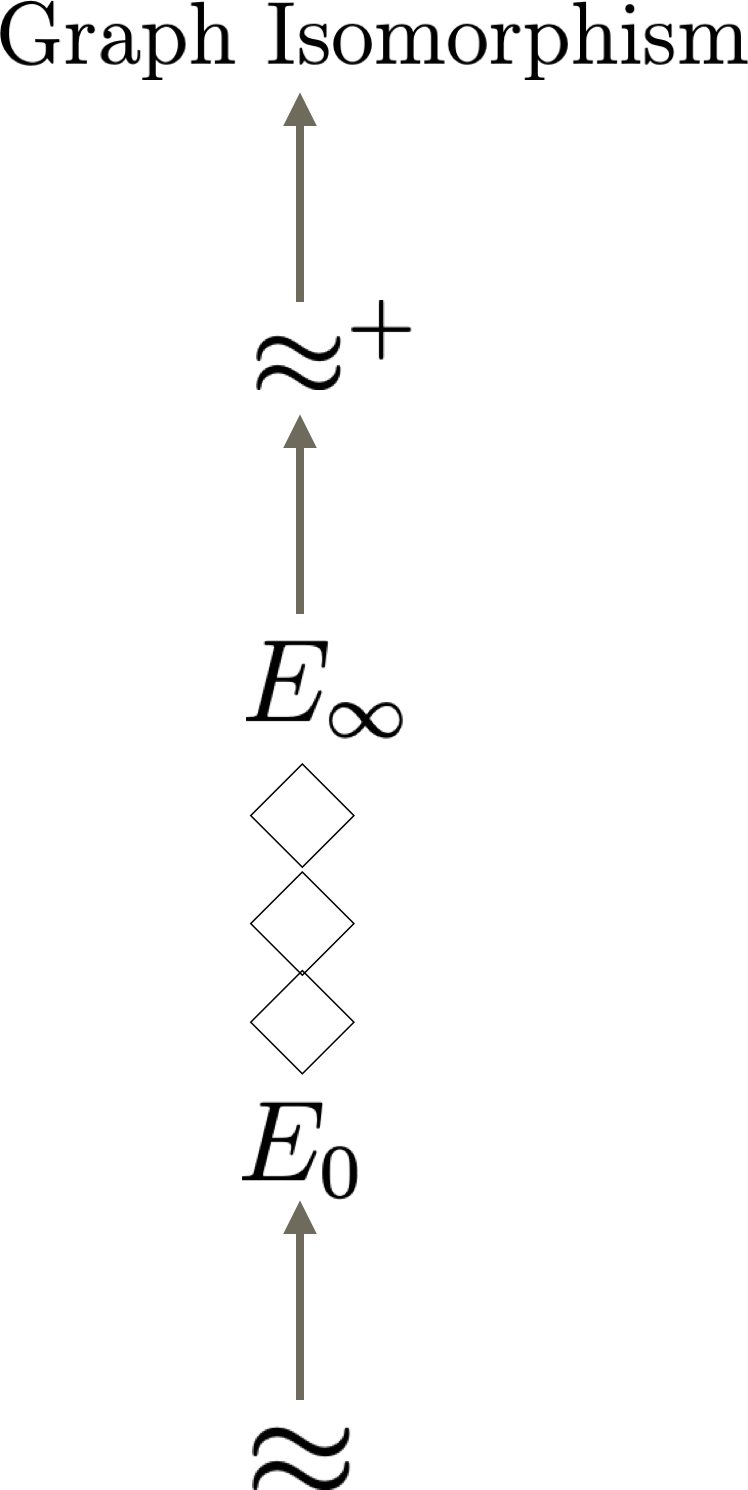}
	\caption{Benchmarks }
	\label{tosinfinty}
	\end{figure}

\subsection{Turbulence}\label{turbulence}
To establish that an equivalence relation is reducible to an $S_\infty$-action, one exhibits a reduction.  Hjorth 
(\cite{Hj}) devised a tool for showing that an equivalence relation is \emph{not} reducible to an $S_\infty$-action.\footnote{The tool was discovered before there were widespread interactions between descriptive set theory and dynamical systems, so the term \emph{turbulance} does not signify any direct connection with chaotic behavior in the sense of dynamics.}

\begin{definition} Let $G$ be a Polish Group acting continuously on a
Polish Space $X$. Then the action is \emph{turbulent} iff:
\begin{enumerate}
\item Every orbit is dense.
\item Every orbit is meager.
\item For all $x,y\in X, V\subset X, U\subset G$ open, with $x\in V, 1\in
U,$ there exists $y_0\in [y]_G$ and $\langle
g_i\rangle_{i\in\mathbb N}\subset U,
\langle x_i\rangle_{i\in\mathbb N}\subset V$ with
\begin{enumerate}
\item $x_0=x$,
\item  $x_{i+1}=g_ix_i$
\item there is a subsequence $i_n$ such that $x_{i_n}$ converges to $y_0$.
\end{enumerate}
\end{enumerate}
\end{definition}

Hjorth proved the following result: 
\medskip

\begin{theorem}\label{turbthm} Suppose that a  continuous action of $G$ on $X$ is
turbulent. Then the orbit equivalence relation cannot be reduced to an $S_\infty$-action.
\end{theorem}
Further results of Hjorth show that \emph{turbulence} almost exactly captures the property of not being reducible to an $S_\infty$-action. (See \cite{Hj}.)

\begin{remark*} The property of being turbulent is generic in the sense that if $Y\subset X$ is a 
$\mathcal G_\delta$ comeager set that is invariant under the $G$-action then the orbit equivalence relation induced by $G$ on $Y$ is also turbulent.
\end{remark*}




\subsection{Polish Group actions}\label{pgas}

The next major benchmark is to ask whether a given equivalence relation is reducible to the equivalence relation given by a Polish Group action.  Because $S_\infty$ is a Polish group this is a more general collection of equivalence relations.

Because many natural equivalence relations are given by Polish group actions this is an enormous collection.  Here are some facts that help navigate the subject. An excellent book on the subject is due to Becker and Kechris (\cite{BK}), although there have been many developments since the book was published. 
	\begin{theorem}\label{1234go} Fix a Polish group $G$.
	\begin{enumerate}
	\item There is a $\preceq^2_\mcb$ maximal equivalence relation given by Polish $G$-actions.
	\item There is a $\preceq^2_\mcb$ maximal equivalence relation given among all  Polish group actions.
	\item If $H\subseteq G$ is a closed subgroup, then the maximal Polish $H$-action is continuously reducible to the maximal Polish $G$-action.
	\end{enumerate}
	\end{theorem}
\noindent Items 1 and 2 are due to Becker and Kechris. Item 3 is essentially due to Mackey in the Borel case and proved by Hjorth in the continuous case.

Uspenski\u{\i} proved that the group of homeomorphisms of the Hilbert cube is a universal Polish group: every Polish group is isomorphic to a closed subgroup of this  group. Thus item 2 follows by applying items 1 and 3 to the group Uspenski\u{\i} used \cite{Gaobook}, \cite{US}).

\section{Standard Mathematical Objects in each region}
\label{benchmarks}

The classes with classifications in the lower regions of the diagram are extremely plentiful so we don't expand on systems with numerical invariants or on equivalence relations \hyperlink{ezero}{with $E_0$} embedded in them, though we give examples of each for dynamical systems. For $S_\infty$ we have already given an example of a Borel action (Remark \ref{borelSinfty}) and $\preceq^2_\mcb$-maximal $S_\infty$-actions (Theorems \ref{thetwogrs} and \ref{camgao}).

\subsection{Countable equivalence relations}\label{ctblequiv}
Recall that an equivalence relation is said to be \emph{countable} if it has countable classes (see Section \ref{countable classes}.)

In section \ref{ctblequiv} we saw that \hyperlink{einfty}{Dougherty, Jackson and Kechris} (\cite{DJKhyperfinite}) showed that there is a $\preceq^2_\mcb$-maximal countable equivalence relation they call $E_\infty$. So understanding countable equivalence relations is, in part understanding their relationship to $E_\infty$.
Rather than try to cover a large literature, we focus on work of Thomas (\cite{thomasclsasprob},\cite{thomasp-local}) that is  related to dynamical systems and gives  natural example of a 
$\preceq_B^2$-maximal equivalence relation with countable classes.

Note that the following sets are naturally viewed as Polish spaces endowed with countable equivalence relations:
	\begin{itemize}
	\item the space of finitely generated groups, with the equivalence relation of isomorphism denoted $\cong_{fin}$.
	\item the space of torsion-free abelian groups of rank $n$,
	with the equivalence relation of isomorphism denoted $\cong_n$
	\item the space of torsion-free, p-local  abelian groups of rank $n$, with the equivalence relation of isomorphism 
	related $\cong_n^p$ (for $p$ a prime number).
	\end{itemize} 

Using Zimmer's notion of \emph{super-rigidity}, and the Kurosh-Mal'cev $p$-adic localization technique, Thomas (\cite{thomasclsasprob},\cite{thomasp-local}) proved the following results except 1 which is due to Baer. The fifth item is joint work of Thomas and Velickovic \cite{TVshow}:\footnote{Hjorth showed that for $n>1$, $\cong_n\not\preceq_\mcb^2\cong_1$ using elementary techniques.  (See \cite{hjcongn}.)}

\begin{figure}[!h]
	\centering
	\includegraphics[height=.50\textheight]{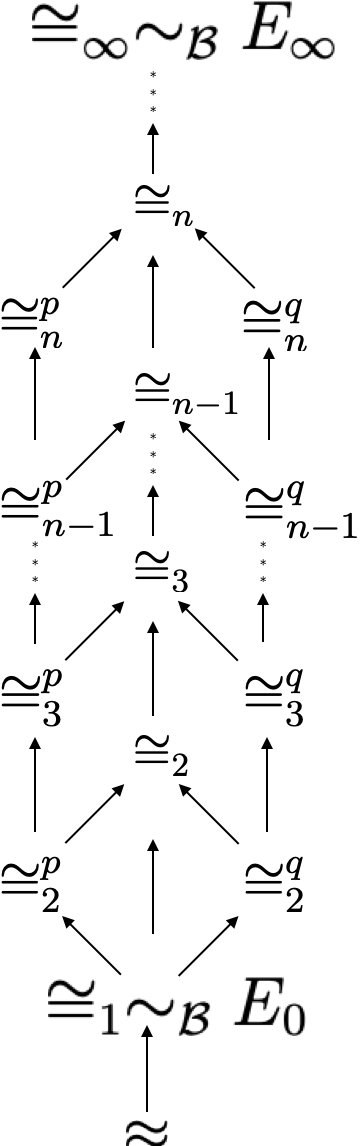}
	\caption{The Thomas diagram}
	\label{thomasdiagram}
	\end{figure}

	\begin{enumerate}
	\item $\cong_1\sim_\mcb E_0$,
	\item for $n\ge 1, \cong_n$ is strictly  $\preceq^2_\mcb$ below $\cong_{n+1}$.
	 In symbols:\\
	 $\cong_n\preceq^2_\mcb\ \cong_{n+1}$,
	\item for distinct primes $p, q$ and $n\ge 2$, the relations $\cong_n^p, \cong_n^q$ are $\preceq^2_\mcb$-incomparable,
	\item for each prime $p$ and each $n\ge 1$, $\cong^p_n\ \preceq^2_\mcb \ \cong_{n+1}$ in the strict sense,
	\item If $\cong_\infty$ is the relation of isomorphism for finitely generated groups, then 
	$E_\infty\sim_\mcb (\cong_\infty)$.
	\end{enumerate}
These results can be best summarized in  figure \ref{thomasdiagram}.

\subsection{General Polish Group actions} We first give an example of a Polish group action that is not reducible to an $S_\infty$-action but is still Borel.

\begin{example*}
Let $\mathcal U$ be the collection of unitary operators on a separable Hilbert Space $H$ and  
$\mathcal N$ be the collection of normal operators on $H$.  Let $E$ be the equivalence relation on 
$\mathcal U$ given by conjugacy of $\mathcal N$. The Spectral Theorem shows that equivalence relation of 
unitary conjugacy of normal operators on a Hilbert space can be reduced to mutual absolute continuity of 
measures on the unit circle.  The latter can be checked to be Borel using the 
\hyperlink{suslintheorem}{Suslin  Theorem}.\footnote{This was extended by Ding and Gao (\cite{DG}) to the 
conjugacy action of all bounded linear operators on the unitary operators.}

However it was proved by Kechris and Sofronides   that this equivalence relation is not reducible to an $S_\infty$-action. They use the turbulence of the $\ell^2$-action on $\mathbb R^\nn$. (See \cite{KeSo}.)
\end{example*}

A related example is the \emph{maximal unitary group action}.  This was shown by  Pestov and Gao (\cite{GP}) to lie   $\preceq^2_\mcb$ above  the maximal $\ell_1$-action, which lies above all Abelian Polish group actions. 
In particular it is not reducible to any $\ell_1$-action. Kanovei (\cite{kanovei}) showed that the maximal unitary group action is not reducible to any $\ell_1$-action.

The collection of Koopman operators that arise from measure-preserving transformations on a measure space $(X, \mcb, \mu)$ is a closed subgroup of the group of unitary transformations on $L^2(X)$. Thus it follows from the Mackey-Hjorth theory (Theorem \ref{1234go}) that the isomorphism action of the measure-preserving transformations on the ergodic transformations is reducible to the maximal unitary group action.

\paragraph{Choquet Simplexes} An example of a natural equivalence relation of maximal 
$\preceq^2_\mcb$-complexity among Polish group actions is due to Sabok.  

Recall that a Chouqet simplex is a compact, separable, affine space $X$ such
 that any point $x\in X$
 is the integral of a unique measure concentrated at the extreme points of $X$. There is  a Choquet 
 simplex, the Poulsen simplex, which is universal in the sense that every Choquet simplex can be 
 embedded as a face. This puts a Polish topology on the collection of Choquet simplexes. 
 Let $(\mcc, \tau)$ be this Polish space. 
 
Set two Choquet simplexes $S_1, S_2$ \emph{equivalent} if there is an affine homeomorphism taking 
 $S_1$ to $S_2$. Let $AFF($Poulsen$)$ be the group of affine homeomorphisms of the Poulsen simplex. Then  $AFF($Poulsen$)$ is a Polish group and its natural action on  $\mcc$ maps faces to faces, so Aff(Poulsen) acts on $\mcc$ by  a Polish action. 
 
Lindenstrauss, Olsen and Sternfeld (\cite{loy}) showed that any affine homeomorphism between two proper subfaces of the Poulsen extends to an affine homeomorphism of the whole simplex.  It follows that the action of $AFF($Poulsen$)$ on $\mcc$ coincides with the relation of being affinely homeomorphic. Thus we can view the equivalence relation of being affinely homeomorphic as the orbit equivalence 
relation of a Polish group action.
 
M. Sabok (\cite{Sabok}) showed that affine homeomorphism of Choquet simplices is 
$\preceq^2_\mcb$-maximal among  equivalence relations induced by Polish group actions:

\begin{theorem}(Sabok) Every orbit equivalence relation of a Polish group action is $\preceq^2_\mcb$-reducible to 
  the equivalence relation of affine homeomorphism of Choquet simplices.
\end{theorem}

\subsection{Borel equivalence relations that don't arise from a Polish Group Actions}
Kechris and Louveau defined a Borel equivalence relation $E_1$ on $\{0,1\}^{\nn\times\nn}$ 
 that is not reducible to any Polish group action and is moreover $\preceq^2_\mcb$-minimal with these properties.
 
 View $\{0, 1\}^{\nn\times \nn}$ as the collection of  infinite sequences $\la x_n:n\in\nn\ra$ where $x_n\in \{0,1\}^\nn$.  Set $\la x_n\ra_nE_1\la y_m\ra_m$ if and only if there is an $N$ for all $n>N, x_n=y_n$. 
 Note that this is an analogue of $E_0$ where elements of $\{0, 1\}$ have been replaced by elements of $\{0,1\}^\nn$ and the equivalence relation is eventual equality.

\begin{theorem}(Kechris-Louveau, \cite{kechlouv})
$E_1$ is Borel and not reducible to a Polish group action. 
\end{theorem}

We now give another example arising in topology. Let $\mathcal X$ be the collection of compact separable metric spaces with $\nn$ as a dense subset. Then $\mathcal X$ carries a natural Polish topology. 

We define the notion of \emph{bi-Lipschitz equivalence}.
For  compact separable metric $\mck_1, \mck_2$ set $\mck_1 R\mck_2$ if and only if there is a bijection 
$f:\mck_1\to \mck_2$ such that both $f, f^{-1}$ are Lipschitz.  The following theorem appears in \cite{Rosendal}.

\begin{theorem}(Rosendal)
The equivalence relation $R\subseteq \mcx\times \mcx$ is Borel but \emph{not} reducible to a Polish group action.
\end{theorem}

The Rosendal theorem is proved by reducing $E_1$ to  bi-Lipschitz equivalence and applying Theorem \cite{kechlouv}. It is conjectured that $E_1$ is $\preceq^1_\mcb$-below every Borel equivalence relation that is not reducible to a Polish group action. This is Open Problem \ref{OP14}.

\subsection{The maximal analytic equivalence relation} While it  was a folklore result that  there is a maximal equivalence relation among all analytic equivalence relation, the construction was quite abstract.   In \cite{FLR} Ferenczi, Louveau and Rosendal gave a natural example. It is an equivalence relation on the Polish space of  separable Banach spaces, which can be viewed as the collection of closed subspaces of 
$\mathbb C([0,1])$ endowed with the Effros Borel structure. 

\begin{theorem}(Ferenczi, Louveau, Rosendal) The relation of isomorphism between separable Banach spaces is maximal among all analytic equivalence relations. 
\end{theorem}
Since there \emph{are} equivalence relations that are not reducible to Polish group actions, and the Ferenczi, Louveau, Rosendal example is maximal, it follows that the Ferenczi, Louveau, Rosendal example does not arise as a Polish group action.

We can now put these examples into the basic benchmark regions:

\newpage

\begin{figure}[h]
\centering
\includegraphics[height=.6\textheight]{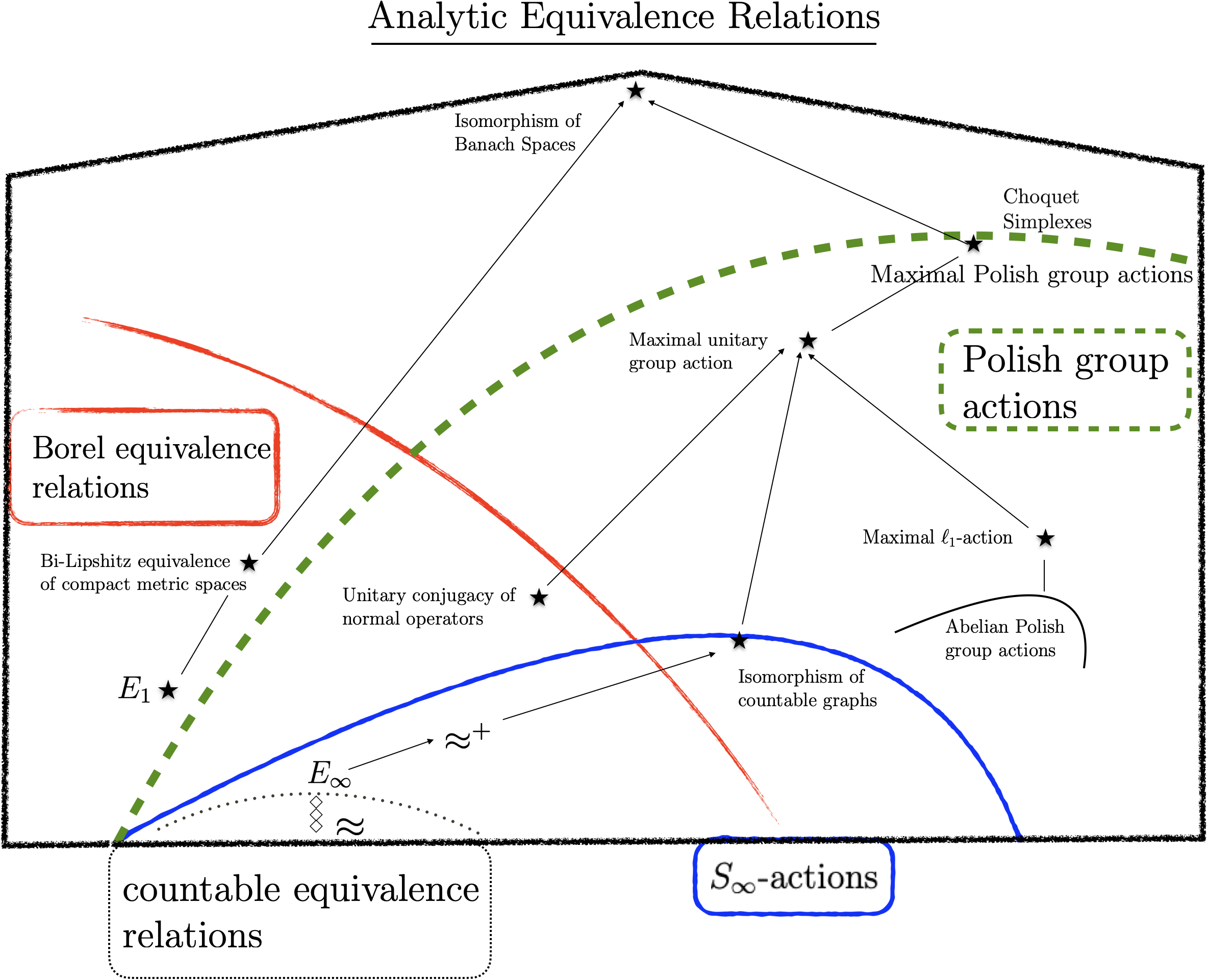}
\caption{Classes placed in regions of complexity }
\label{zoonodyna}
\end{figure}

\section{Placing dynamical systems in each section}
\subsection{Measure Isomorphism}\label{measureiso}
The relationship between measure isomorphism and descriptive complexity is better understood than the relationship between topological conjugacy and descriptive set theory.

Work of Mitin (\cite{Mitin}) and Ryzhikov (\cite{Ryzhikov}) studies measure-preserving transformations and show that they encode the theory of Peano Arithmetic.  As a result many questions about them are recursively undecidable.  In this section we will focus on Borel reducibility.

\subsubsection{Systems with Numerical Invariants}
The clearest example here is the class of Bernoulli Shifts. As noted in Section \ref{bs}, the collection of measure-preserving transformations isomorphic to a Bernoulli shift has a structure theory and Ornstein's Theorem (Theorem \ref{OrnErg}) shows that entropy is a complete numerical invariant.

\subsubsection{Reducing $E_0$ to dynamical systems} \label{e0toK} A remarkable theorem of Feldman from a paper published in 1974 (\cite{feldman}), shows the following result:

\begin{theorem}
The equivalence relation $E_0$ is continuously reducible to isomorphism of $\mck$-automorphisms.
\end{theorem}
Feldman's paper predated the language and definitions given in this paper, however his proof is easily translated into the proof given in section \ref{GE} after theorem \ref{HKL}.

Surprisingly it is not known if the conjugacy relation for $\mck$-automorphisms is Borel. (See Open Problem \ref{OP7}.)  Nor is it known if it can be reduced to an $S_\infty$-action.

\subsubsection{A class bi-reducible to $=^+$.}
In theorem \ref{hvntheory}, Halmos and von Neumann showed that two discrete spectrum measure-preserving transformations are isomorphic if and only if their associated Koopman operators have the same eigenvalues. In tracing the proof (see \cite{descview}) one sees that this is a reduction of isomorphism 
to $=^+$ on the complex unit circle. 

The reduction in the other direction also uses the Halmos-von Neuman theory.  First observe that there is a dense $\mcg_\delta$ subset $I$ of  the complex unit circle $S^1$ such that any $\xi_1, \xi_2, \dots \xi_n$ which belong to $I$ are algebraically independent. 

To each $\vec{\xi}=\la \xi_n:n\in\nn\ra$ with $\xi_n\in I$,  let $G_{\vec{\xi}}$ be the multiplicative subgroup of $S^1$ generated by $\vec{\xi}$ and consider its Pontryagin dual $\widehat{G_{\vec{\xi}}}$ consisting of all characters of the group $G_{\vec{\xi}}$.  Let $\Chi$ be the identity map from $G_{\vec{\xi}}$ to the unit circle.  Then $\widehat{G_{\vec{\xi}}}$ is a monothetic abelian group generated by $\Chi$. Multiplication by $\Chi$ in $\widehat{G_{\vec{\xi}}}$ preserves Haar measure and is an ergodic measure-preserving transformation $T_{\vec{\xi}}$.  The Koopman operator of $T_{\vec{\xi}}$ has $G_{\vec{\xi}}$ as its group of eigenvalues.

Summarizing: the map that takes $\vec{\xi}$ to $(\widehat{G_{\vec{\xi}}}, \mcb, \mbox{Haar}, T_{\vec{\xi}})$ is a reduction of the equivalence relation $=^+_I$ to isomorphism of discrete spectrum measure-preserving transformation. The upshot is the following which appeared in \cite{descview}:

\begin{theorem}(Foreman-Louveau)
Isomorphism of discrete spectrum measure-preserving transformations is Borel bi-reducible to $=^+$.
\end{theorem}

\subsubsection{Isomorphism is not reducible to any $S_\infty$-action.}
Hjorth \cite{Hjorth} showed that there is a turbulent equivalence relation that is $\preceq^2_\mcb$ reducible to isomorphism for ergodic height 2 measure distal transformations (see section \ref{turbulence} for a discussion of turbulence). It follows that:

\begin{theorem}(Hjorth)
Isomorphism for ergodic measure-preserving transformations is not reducible to an $S_\infty$-action.
\end{theorem}

Foreman and Weiss (\cite{ananti}) improved Hjorth's result by showing that the relation of isomorphism of ergodic measure transformations is \emph{itself} a turbulent equivalence relation.

\begin{theorem}(Foreman-Weiss) Let $MPT([0,1])$ be the group of measure-preserving actions on $[0,1]$ with the usual Polish topology and $EMPT([0,1])$ be the collection of ergodic transformations.  Then the conjugacy action of $MPT([0,1])$ on $MPT([0,1])$ is turbulent. 
\end{theorem}
From the remark following Theorem \ref{turbthm},  that generic classes of measure-preserving transformations, such as the ergodic or weakly mixing transformations, are also turbulent. Hence isomorphism is not reducible to an $S_\infty$-action. 

\begin{corollary}\label{nosinf}
Let $\mcc$ be a dense $\mcg_\delta$ collection of measure-preserving transformations in $[0,1]$. Then the isomorphism relation on $\mcc$ is not reducible to an $S_\infty$-action.
\end{corollary}

\subsubsection{Isomorphism is not Borel}
In \cite{Hjorth} Hjorth showed that the isomorphism relation between measure-preserving transformations is not Borel. However his proof used non-ergodic transformations in an essential way.

Foreman, Rudolph and Weiss showed the following results (\cite{FRW}):

\begin{theorem}(Foreman, Rudolph, Weiss)\label{FRWtrans}
Let $G$ be the group of measure-preserving transformations of $[0,1]$ with the standard Polish topology and $EMPT$ be the ergodic transformations in $G$.  The conjugacy action on $G$ is complete analytic. In particular it is not Borel.
\end{theorem}

Using the method of suspensions on can lift these results to $\mathbb R$-actions (\cite{FRW}).
\begin{corollary}
The isomorphism relation between ergodic measure-preserving $\poR$-flows is complete analytic, and hence not Borel.
\end{corollary}

In unpublished work, Foreman was able to establish the following:
\begin{theorem}\label{ctbletomeas}
Isomorphism of countable graphs is reducible to  isomorphism for measure-preserving transformations.
\end{theorem}

Combining these results with Corollary \ref{nosinf}, we get:
\begin{corollary}
The equivalence relation of isomorphism of measure-preserving transformations is strictly above $S_\infty$-actions in $\preceq^2_\mcb$.
\end{corollary}
Theorem \ref{ctbletomeas} also provides a second proof of the Kechris-Sofronidis result that the maximal unitary action is strictly above the maximal $S_\infty$-action.

One can ask for more information in placing the relation of measure isomorphism on $EMPT$ in the 
$\preceq^2_\mcb$ relation.  It was proved by Kechris and Tucker-Drob in \cite{kechdrob} that the relation of \emph{unitary equivalence of normal operators} is reducible to measure isomorphism on $EMPT$.  Since the unitary equivalence is Borel and isomorphism of ergodic measures is not, unitary equivalence is strictly $\preceq^2_\mcb$ below isomorphism of ergodic measures.

\begin{theorem}(Kechris, Tucker-Drob)
The relation of \emph{unitary equivalence of normal operators} is $\preceq^2_\mcb$ below \emph{isomorphism of ergodic measure-preserving transformations}. By Theorem \ref{FRWtrans} the relation is strict.
\end{theorem}

Shortly after this paper was written B. Weiss was proved the following result, as yet unpublished:

\begin{theorem*}(Weiss) Let $G$ be a countable group. Let $X$ be the space of homomorphisms of $G$ to 
$MPT([0,1])$ such that the action is free and ergodic. If $G$ is of the form $H\times \poZ$, then there is a Borel reduction from the isomorphism relation on $X$ to the isomorphism relation on $\poZ$-actions.
\end{theorem*}

\begin{corollary*}
If $G=\poZ^d$ for $d\ge 1$ then there is a Borel reduction from free measure preserving $G$-actions to measure preserving $\poZ$-actions.
\end{corollary*}

\subsubsection{Smooth measure-preserving transformations}
Much of the motivation for studying measure isomorphism is to classify the quantitative behavior of solutions to ordinary differential equations. The reductions built in Theorems \ref{FRWtrans} and \ref{ctbletomeas} take their ranges in the collection of ergodic measure-preserving transformations with non-trivial odometer factors.  There are no known measure-preserving transformations on compact manifolds that have odometer factors. This is a special case of the classical \emph{Realization Problem} (Open Problem \ref{OP5}). A statement of the problem appears in 1975 in a paper of Ornstein (\cite{ornstein4}) and the problem is attributed to Katok by Hasselblatt in \cite{hassel}.

The natural question to ask is whether the measure isomorphism relation between volume preserving diffeomorphisms on a compact manifold is Borel.  Does the additional structure provided by being $C^\infty$ give enough information about measure-preserving transformations to determine isomorphism? 

\noindent The answer is no:

\begin{theorem}(Foreman, Weiss)\label{noforsmoothmpts} The isomorphism relation between Lebesgue measure-preserving $C^\infty$ transformations on the 2-torus $\bt^2$ is complete analytic.  In particular it is not Borel.
\end{theorem}
\begin{remark*}The techniques of Theorem \ref{ctbletomeas} adapt exactly to prove the analogous result for smooth measure-preserving transformations. However, because the isomorphism relation on smooth measure-preserving transformations is not a group action, the technique of \emph{turbulence} is not directly applicable. Thus it is possible that the measure isomorphism relation on $C^\infty$ transformations on the 2-torus is $\preceq^2_\mcb$ bi-reducible with the maximal $S_\infty$-action. This is the content of Open Problem \ref{OP8}.
\end{remark*}

The proof of Theorem \ref{noforsmoothmpts} involves three steps:
\begin{enumerate}
	\item Identify a class of symbolic shifts (the \emph{circular systems}) that are realizable using the Anosov-Katok \emph{ABC} method. 
	\item Show that the class of measure-preserving transformations containing an odometer has the same factor and isomorphism structure as the class of circular systems. They are ``isomorphic" categories by an isomorphism $\Phi$.
	\item Build a reduction from the space of ill-founded trees to the odometer systems, use the isomorphism $\Phi$ to transfer the range to the circular systems and then the ABC method to realize the results as diffeomorphisms. 
\end{enumerate}
The proofs appear in the three papers: \cite{circsys}, \cite{globstruct}, \cite{noclass}.

We note that these results give no upper bound on the complexity of the isomorphism relation for ergodic measure-preserving transformations, other than the obvious one--that of the maximal Polish group action. Hence we ask the question:
	\begin{quotation}
	\noindent Is the maximal Polish group action reducible to isomorphism of ergodic measure-preserving transformation?
	\end{quotation}
\noindent This is Open Problem \ref{OP9}.

\subsubsection{Kakutani Equivalence}\label{kaku}
An important relation in ergodic theory is Kakutani equivalence (see \cite{Thouvenot} for a nice exposition). 

Let $T$ be a measure-preserving transformation on a standard measure space $(X, \mcb, \mu)$.  
Let $A\in \mcb$ be a set of positive measure and define $\mu_A:P(A)\cap \mcb\to [0,1]$ by setting 
\[\mu_A(B)=\frac{\mu(B)}{\mu(A)}.\] 

By Poincare recurrence, for almost all $x\in A$, there is an $n>0$ such that $T^{n}(x)\in A$. Let $n(x)$ be the least such $x$ and $T_A(x)=T^{n(x)}(x)$. Then $T_A:A\to A$ and is a measure-preserving transformation. If $T$ is ergodic then $T_A$ is ergodic.

\begin{definition}\label{Kakeqdef}
Two ergodic transformations $\xbmt$ and $\ycns$ are said to be Kakutani equivalent if there are positive measure sets $A\in \mcb$ and $B\in \mcc$ such that $T_A$ is isomorphic to $S_B$.
\end{definition}

Gerber and Kunde \cite{gerberkunde} proved the following theorem, following the general steps of the proof of Theorem \ref{noforsmoothmpts}, but with  significant improvements. 

\begin{theorem}\label{kakisbad}
The equivalence relation of Kakutani equivalence on ergodic measure-preserving diffeomorphisms of the 2-torus is complete analytic. Hence it is not a Borel equivalence relation.
\end{theorem}
Because Kakutani equivalence  for measure-preserving diffeomorphisms is $\preceq^2_\mcb$ reducible to Kakutani equivalence for general measure-preserving transformations, the following is immediate:\footnote{Chronologically, the corollary preceded the theorem.}
\begin{corollary}\label{which came first?}
The equivalence relation of Kakutani equivalence on the space of ergodic measure-preserving transformations on $[0,1]$ is complete analytic and hence not Borel.
\end{corollary}

\noindent Gerber and Kunde also proved an $\mathbb R$-flow version of Corollary \ref{which came first?}. 

\begin{definition}Let $F=\{f_t:t\in \mathbb R\}$  be a measure-preserving $\mathbb R$-flow on 
$(X,\mu)$. A \emph{reparametrization}  of $F$ is a jointly  measurable function $\tau:X\times \mathbb R\to \mathbb R$ such that for almost every $x$ the map $t\mapsto \tau(x,t)$ is a homeomorphism of $\mathbb R$ and the map $(x,t)\mapsto f_{\tau(x,t)}(x)$ is a measure-preserving flow on $X$. 

If $G=\{g_t:t\in \mathbb R\}$ is another measure-preserving flow on $X$, then $F$ and $G$ are \emph{isomorphic up to a homeomorphic time change} if there is a reparameterization $\tau$ of $F$ and a measure-preserving transformation $\phi:X\to X$ such that for all $t$ and almost every $x$:
\[g_t(\phi(x))=\phi(f_{\tau(x,t)}(x)).\]
\end{definition}

\noindent The Gerber-Kunde techniques extend to show:

\begin{corollary}
The equivalence relation of \emph{isomorphism up to a homeomorphic time change} for $C^\infty$ measure-preserving flows is complete analytic and hence not Borel.
\end{corollary}

The proofs of the Gerber-Kunde theorems do not give information about $\preceq^2_\mcb$ reducibility to any of the benchmark analytic equivalence relations. (This is Open Problem \ref{OP10}.)

\subsubsection{The summary diagram}
The following diagram summarizes this section.  The $\mck$-automorphisms are not on the diagram 
because, other than being above $E_0$, their position is not known.  For the same reason the Kakutani 
equivalence relation is not on the diagram.

\newgeometry{left=1cm,bottom=3cm}

\restoregeometry


\subsection{Topological Conjugacy}\label{homeoclass}
In this section we discuss the \emph{qualitative behavior} as proposed by Smale. We do this in both the zero-dimensional and the smooth contexts.  ``Equivalence" in this context  means \emph{conjugate by a homeomorphism}.

\paragraph{After this paper was finished} Joint work of the author and Gorodetski improved the results for smooth transformations on manifolds with the equivalence relation of \emph{topological conjugacy} to all dimensions at least one.  These will appear in a future paper.

\subsubsection{Smooth transformations with numerical invariants}
As discussed in section \ref{smodyn}, the  collection of smooth transformations that is structurally stable has complete numerical invariants.  These include the classes of Anosov  and Morse-Smale diffeomorphisms.  It is not clear how far, if at all, these classes can be extended and still have complete numerical invariants. 
\subsubsection{Diffeomorphisms of the torus above $E_0$.}
Foreman and Gorodetski showed that there are diffeomorphisms of any manifold of dimension at least 2 that are $\preceq^2_\mcb$ above $E_0$. It follows that  topological conjugacy for general diffeomorphisms is not numerically classifiable.

\begin{theorem}(Foreman, Gorodetski) \label{fogo1} 
Let $M$ be a $C^k$ manifold of dimension at least 2 and $1\le k\le \infty$. Then there is a continuous reduction $R$ from $E_0$ to the collection of $C^k$ diffeomorphisms with the relation of topological conjugacy.  
\end{theorem}
\begin{remark*}
In fact, the construction creates a reduction that takes values on the boundary of the Morse-Smale diffeomorphisms of the torus.
\end{remark*}
\noindent We give here a nearly honest proof, which is relatively easy to illustrate with a picture.

\begin{figure}[h]
\centering
\includegraphics[height=.4\textheight]{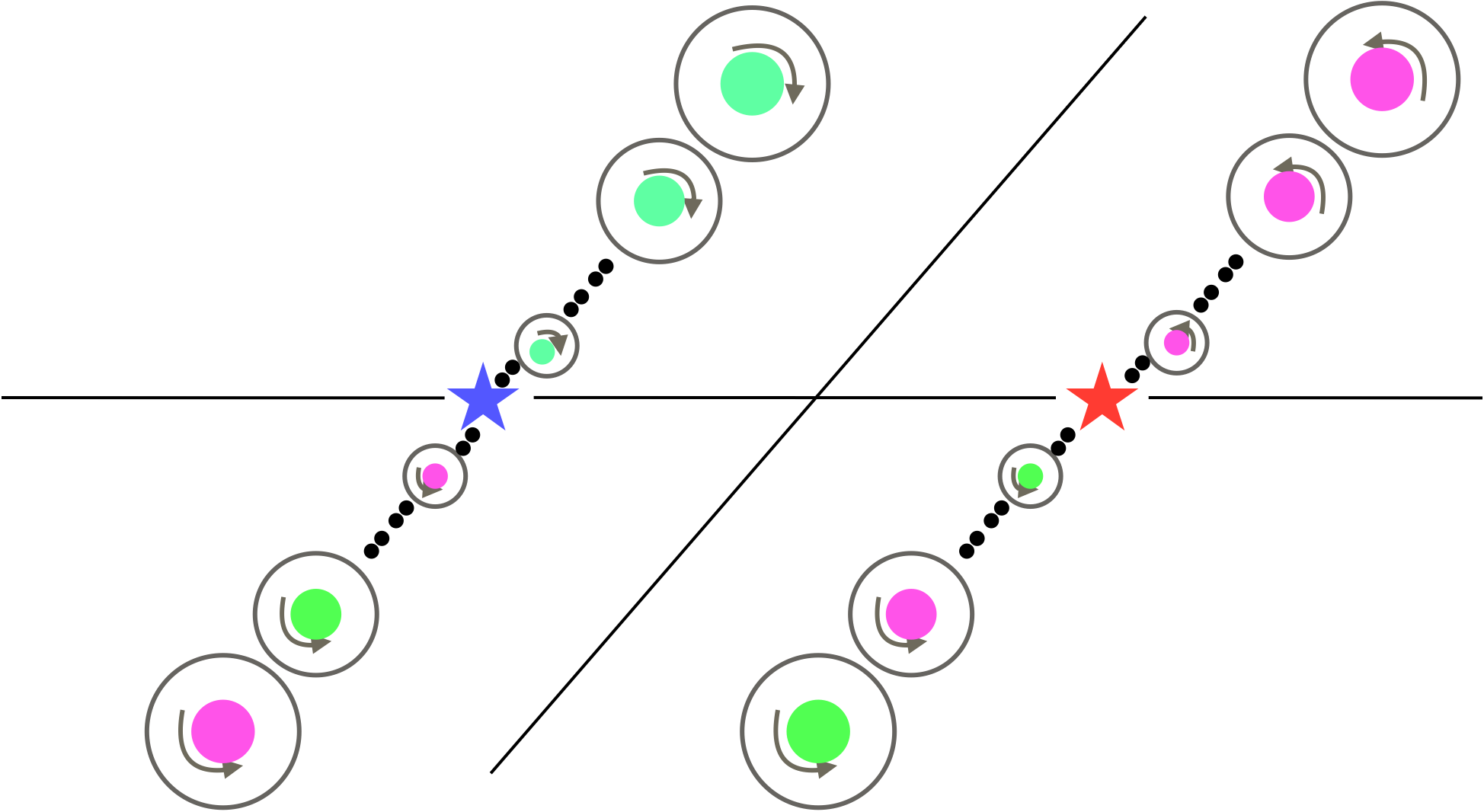}
\caption{How to embed $E_0$ }
\label{E0pic}
\end{figure}

Let $M$ be a manifold of dimension at least 2. We take it to have  dimension 2 for convenience (see figure \ref{E0pic}).
We work inside a neighborhood diffeomorphic to $\mathbb R^2$.   Choose two points (the left star and the right star) and four sequences of disks.  The disks' diameters go to zero at the same relatively fast rate and two of them converge to the left point and two of them converge to the right point.   There are two lines, that separate the region into four parts. Each part contains exactly one of the sequences of disk.  One of the lines passes through the left and right points, giving a notion of upper and lower sequences.

	\begin{enumerate}
	\item The support of the diffeomorphism is the closure of the interiors of the disks. 
	\item The diffeomorphism rotates each interior circle around the center of each disk and 
	each disk contains a subdisk that the diffeomorphism rotates by a constant amount. (The 
	subdisks are indicated by the shading.)	
	\item Each of the sequence of disks above the horizontal line is identified with a fixed sequence of  
	rotation numbers tending to zero 
	and these are different for the two sequences.
	\item There is a fixed sequence of different rotation numbers with $|r_n|<1/2$  and 
	$\la r_n:n\in\nn\ra$  tending to zero, with 
	infinitely many
	 positive and infinitely many negative. 
	 \item The lower subdisks have rotation numbers $\la \pm r_n:n\in\nn\ra$. 
	 \item The lower subdisks then code an infinite sequence of 0's and 1's depending on whether they 
	 agree on the sign of the rotation in the $n^{th}$ subdisk.
	 \end{enumerate}

The reduction $R$ maps $\{0, 1\}^\nn$ to diffeomorphisms built with the two points and the sequences of disks.  Given $f:\nn\to \{0,1\}$,  $R(f)$ is the diffeomorphism that 
\begin{itemize}
	\item rotates the $n$-th subdisc of lower right sequence in the positive direction and lower left sequence in the negative direction if $f(n)=1$
	\item rotates the $n^{th}$ subdisk of the lower left sequence in the positive direction and the $n^{th}$ 
	subdisk of the lower right sequence in the negative direction if $f(n)=0$.
\end{itemize}

By construction, on the range of $R$, each conjugating homeomorphism of $\mathbb R^2$ must take the sequences of upper disks to themselves, and only swap finitely many of the lower disks. Thus $R$ is a reduction of $E_0$ to topological conjugacy.

This example is, in many ways, unsatisfying.  It is the identity on large portions of the space. This raises myriad questions. A simple one to frame is the following (See Open Problem \ref{OP11}.)

\begin{quotation}
\noindent Does $E_0$ reduce to the collection of topologically minimal diffeomorphisms of the 2-torus with the relation of  topological conjugacy?
\end{quotation}

\subsubsection{A dynamical class  that is $\preceq^2_\mcb$ maximal for countable equivalence relations}\label{dynocount}
A well-studied class of dynamical systems is the collection of \emph{finite subshifts}.  These are closed, shift-invariant subsets of $\Sigma^\poZ$ for a finite alphabet $\Sigma$.

Let $\mcc$ be the class of finite subshifts. If $\Sigma$ is the alphabet then a \emph{finite code} or \emph{block code}  is a 
function 
$c:\Sigma^l\to \Sigma$ for some $l\in \nn$. Then $c$ determines a mapping 
$c^*:\Sigma^\poZ\to \Sigma^\poZ$, where $c^*(f)$ 
is achieved by  ``starting at zero and sliding $c$ along $f$ in both directions."  At each location $f$ restricted to 
the $l$-sized window is input into $c$ which outputs a value in $\Sigma$.  

If $\bk$ and $\bk'$ are compact subshifts and conjugate by a homeomorphism $h$, then that 
homeomorphism is given by a pair of  finite codes, one for $h$ and one for $h^{-1}$ . (See \cite{dougmarc}.)
Since there are only countably many finite codes, it follows that the equivalence relation has only 
countable classes.

Clemens (see \cite{clemens}) showed that general problem of homeomorphism is 
$\preceq^2_\mcb$-maximal among equivalence relations with countable 
classes.

\begin{theorem}(Clemens)
The relation of topological conjugacy between subshifts of $\Sigma^\poZ$ (with $\Sigma$ finite) is 
$\preceq^2_\mcb$ maximal among analytic equivalence relations with countable classes.
\end{theorem}

\subsubsection{An example of a maximal $S_\infty$-action from dynamical systems.}\label{camerlogao}
Since every perfect Polish space contains a homeomorphic copy of the Cantor set, it is natural to ask what the complexity is of the collection of homeomorphisms of the Cantor set with the relation $E$ of topological conjugacy. This was answered by Camerlo and Gao (\cite{CG}). 

\begin{theorem}(Camerlo-Gao)
Topological conjugacy of homeomorphisms of the Cantor set is bi-reducible with a maximal $S_\infty$-action.
\end{theorem}

As in the case of the  smooth examples that exhibit high complexity (Theorems \ref{fogo1}, \ref{fiveandabove}) the transformations in the range of the reduction are far from being minimal.  Thus we have the corresponding question: is topological conjugacy of minimal homeomorphisms of the Cantor set bi-reducible  with the maximal $S_\infty$-action?  What about topologically transitive? (See Open Problem \ref{OP12.5}.)

\subsubsection{Topological conjugacy of diffeomorphisms is not Borel}
It might be  expected that for objects as concrete as a $C^k$ diffeomorphism of a compact manifold, topological conjugacy would be Borel.  However, in dimensions 5 and above it is known not to be the case.

\begin{theorem}\label{fiveandabove} (Foreman, Gorodetski) Let $M$ be a $C^k$ smooth manifold of dimension at least five and $1\le k\le \infty$. Then the relation of \emph{Graph Isomorphism} on countable graphs is reducible to topological conjugacy of $C^k$ diffeomorphisms of $M$.
\end{theorem}
Because \emph{Graph Isomorphism} is a $\preceq^2_\mcb$-maximal $S_\infty$-action, it is complete analytic. Hence the following corollary is immediate.
	\begin{corollary}
	Let $M$ be a smooth manifold of dimension at least five. Then the equivalence relation of topological conjugacy on the space of diffeomorphisms of $M$ is complete analytic. In particular it is not Borel.
	\end{corollary}
\noindent 	
The reduction used in the  proof of Theorem \ref{fiveandabove} can be explained with pictures 
(although a rigorous proof requires some work).

Fix a graph $G$ with vertices $\nn$. Consider the solid unit ball in $\mathbb R^3$ and choose a set $I=\la p_n:n\in\nn\ra$ of countably many algebraically independent points tending to the south pole. Draw line segments between each pair of points and between each point and the south pole. The result is a complete graph on countable many points: the independent points together with the south pole. 

\begin{figure}[h]
\centering
\includegraphics[height=.2\textheight]{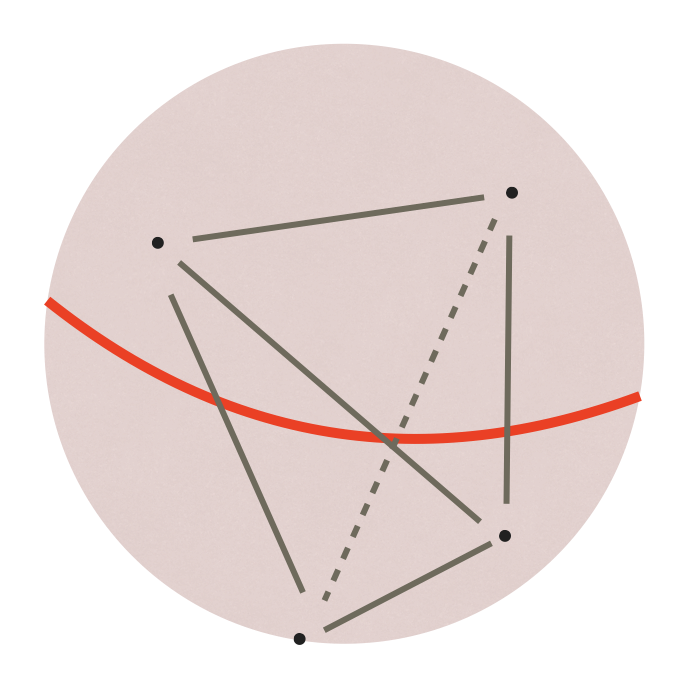}
\caption{The unit ball with edges through it connecting independent points.}
\label{uball}
\end{figure}

Let $G=(\nn, E)$ be an arbitrary countable graph.  Start the reduction by identifying $\nn$ with the countably many points, but not the south pole. (Figure \ref{uball}.)

The coding device is to put disjoint five dimensional  ``cigars" around each of the lines between the points in $I$.  Given a graph $G$ with vertices $\nn$, the support of the diffeomorphism $f_G$ that will be the image of $G$ under that reduction, will be the closure of the cigars. The diffeomorphism $f_G$  either flows towards the center of the cigar if the points connected by the cigar are connected in $G$ or flows away from the center if they are not connected.

\begin{figure}[h]
\centering
\includegraphics[height=.2\textheight]{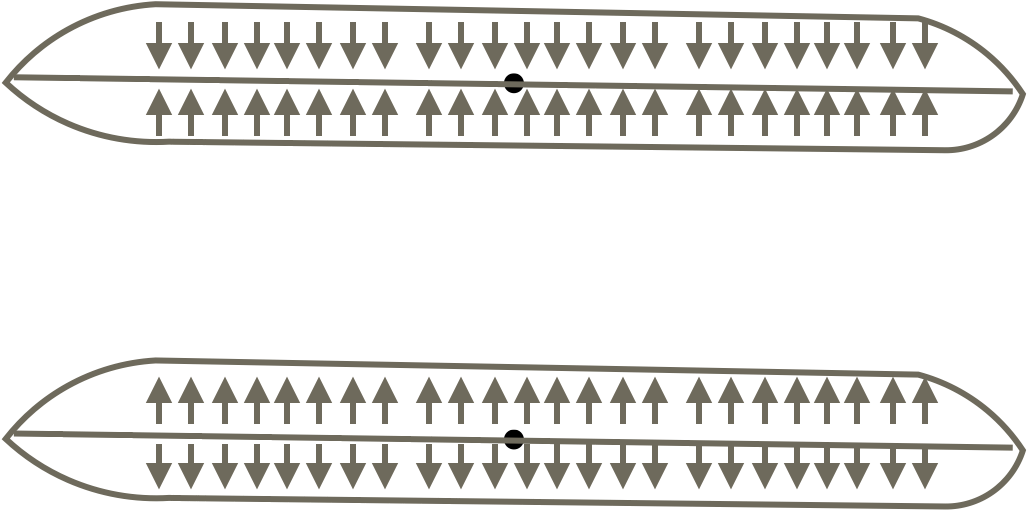}
\caption{Flowing towards and away from the centerline.}
\label{cigars}
\end{figure}

To show this is a reduction one must show that if we are given two graphs $G=(\nn, E)$ and $H=(\nn, F)$, then the reduction produces diffeomorphisms $f_G$ and $f_H$ such that 
	\begin{itemize}
	\item $G$ is isomorphic to $H$
	 \item[] if and only if
	\item $f_G$ is topologically equivalent to $f_H$.
	\end{itemize}
If $f_G$ is topologically equivalent to $f_H$ then the witnessing homeomorphism fixes the south pole and  induces a permutation of the  points in $I$.  This, in turn, gives a permutation of the natural numbers that is an isomorphism between the graphs $G$ and $H$.

\begin{figure}[h]
\centering
\includegraphics[height=.4\textheight]{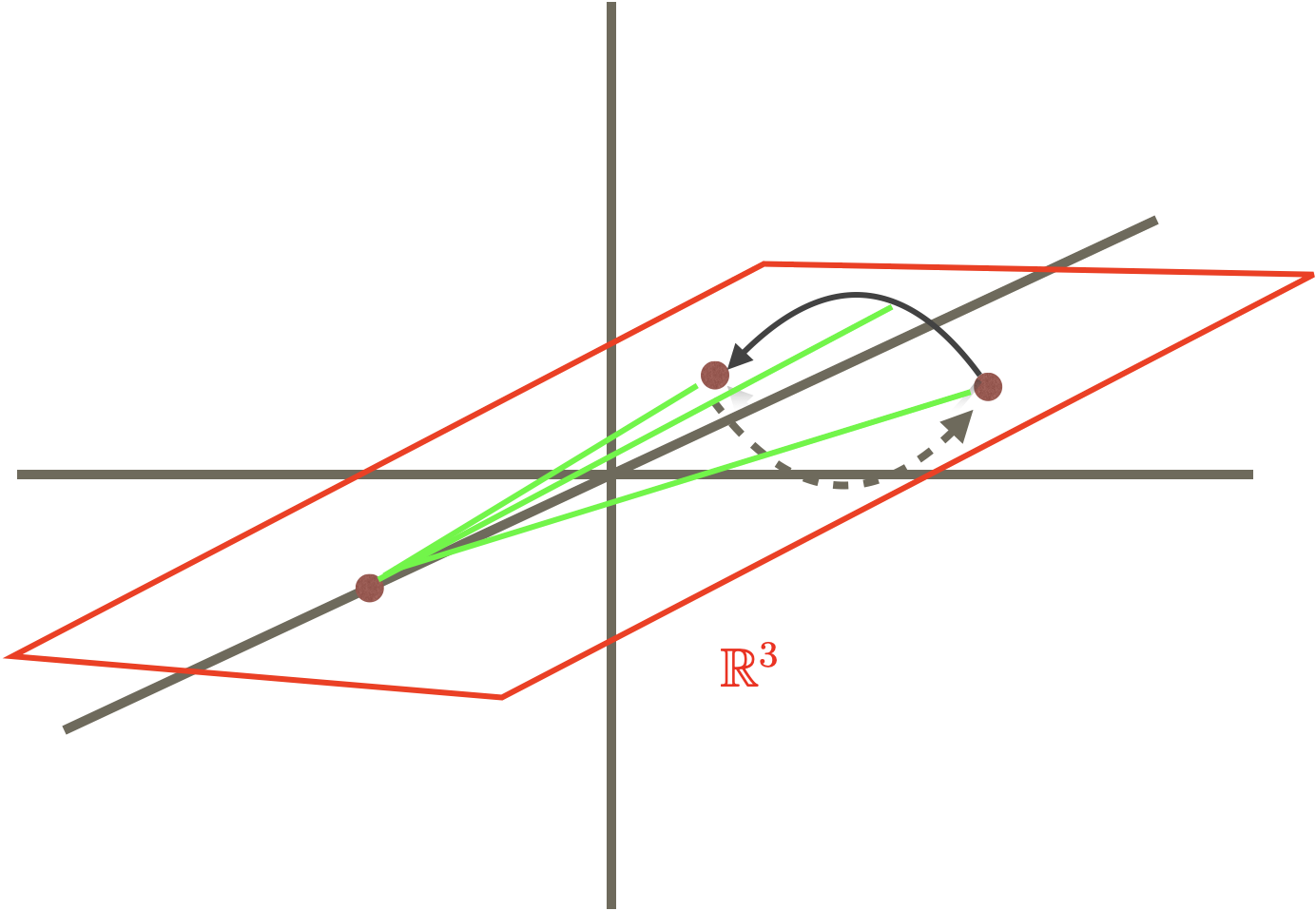}
\caption{Using the fourth dimension to swap two elements of the unit ball in $\mathbb R^3$.}
\label{swap}
\end{figure}

The difficult direction is the opposite. Suppose that we are given a bijection $\phi:\nn\to \nn$ that is an isomorphism between $G$ and $H$.  We need to find a homeomorphism that  mimics $\phi$'s behavior on the countable set $I$ and fixes the south pole of the ball in $\mathbb R^3$.  For a single transposition, this is easy in four dimensions. (See figure \ref{swap}.)  However even two transpositions can interfere with each other.  Having five dimensions gives  sufficient room to prevent this.

\subsubsection{The summary diagram}
The following diagram summarizes the situation for topological conjugacy of diffeomorphisms, as of July 1, 2022. The section above describes what was known on February 1, 2022 and hence does not exactly correspond with the diagram.

\pagebreak
\newgeometry{left=1cm,bottom=3cm}
\pagenumbering{gobble}

\restoregeometry


\appendix
\section{Descriptive Set Theory facts}\label{DSTsum}

We give a very quick review of descriptive set theoretic facts used in this article.  For a  survey that includes complete 
proofs of the facts stated here, see the article \emph{Naive Descriptive Set Theory} (\cite{naivedst}).  That article does 
not assume any familiarity with logic.

We assume the reader is familiar with the ordinals (see Halmos' \emph{Naive Set Theory}) which are canonical 
representatives for any well-ordering.  Ordinality is concerned with order, as opposed to cardinality which is concerned 
with size as determined by bijections. The start with the natural numbers  $0, 1, 2, \dots$ which form the ordinal 
$\omega$, and continue by putting a point on top to get $\omega+1$:

\begin{figure}[h]
\centering
\includegraphics[height=.1\textheight]{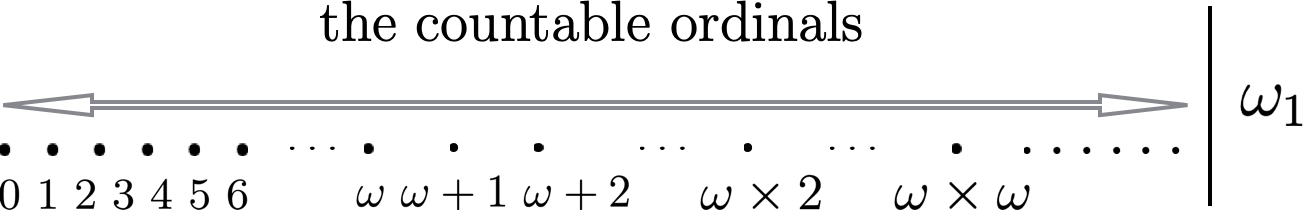}
\caption{A few ordinals}
\label{ordinals}
\end{figure}

After this a second  copy of $\omega$ to get $\omega+\omega=\omega \times 2$. Every ordinal $\alpha$ has an 
immediate successor, $\alpha+1$, and every collection of ordinals has a supremum. Every ordinal is either a 
successor ordinal or a limit ordinal.  Informally ordinals are defined to be \emph{the set of smaller ordinals}. So $0$ is 
the empty set, $1=\{0\}$ and so forth. (See \cite{halmosnaive}, \cite{Jech} or \cite{Levy} for explanation.) We will use 
this convention here: $\alpha=\{\beta:\beta\mbox{ is an ordinal and }\beta<\alpha\}$. 
\medskip

In ZFC, the ordinals give canonical representatives of the isomorphism types of every well-ordering. The Axiom of 
Choice is equivalent to the statement that every set can be well-ordered. It follows that    there are ordinals of every 
cardinality. The smallest uncountable ordinal is called 
$\omega_1$.  Because the collection of real numbers is uncountable, it has cardinality at least the cardinality of 
$\omega_1$ and Cantor's \emph{Continuum Hypothesis} says that the two cardinalities are the same.

\medskip

\hypertarget{pspaces}{Descriptive set theory is largely concerned with \emph{Polish Spaces} $X$--these are 
topological spaces whose topology is compatible with a separable complete metric.}  A Polish space is \emph{perfect} 
if it contains no isolated points. A measure on a Polish space is \emph{standard} if open sets are measurable, it is non-
atomic, complete and separable.

Very very roughly, the main topic of descriptive set theory is understanding what can be done without the use of the 
uncountable Axiom of Choice. For example, the Polish spaces contain the hierarchy of \emph{Borel Sets}, the smallest 
$\sigma$-algebra of subsets $X$ containing the open sets. This algebra is built by induction on the ordinals in a 
hierarchy of length $\omega_1$. Heuristically it contains all sets that can be built with inherently countable information. 

For classifying dynamical systems, Borel sets are not sufficient.  One needs \emph{analytic} and \emph{co-analytic} 
sets. These are extensions of the Borel sets that are nonetheless universally measurable.

{We now give an inductive construction of the Borel} sets that begins with the open sets. To keep track of the level of 
complexity 
of the set, we introduce the following  \hypertarget{BH}{``logical" notation:}

\begin{definition}\label{borel hi}
Let $(X, \tau)$ be a Polish topological space. The levels of the Borel hierarchy are as follows:
\begin{enumerate}
	\item $\bsigma_{1}^{0}$ sets are the open subsets of $X$. 
	\item $\bpi_{1}^{0}$ sets are the closed sets.
	\item Suppose the hierarchy has been defined up to (but not including) the ordinal level $\alpha <\omega_1$. 
Then $A\in$ $\bsigma_{\alpha}^{0}$ if
	and only if	there is a sequence ${\la}B_{i}\ |\ i\in\omega{\ra}$ of subsets of $X$ with $B_{i}\in\bpi^0_{\beta_i}$ 
such that for each $i\in\omega$, $\beta_i<\alpha$
	and
	\[A=\bigcup_{i\in\omega} B_i.\]
	\item $B\in\bpi_{\alpha}^{0}$ if and only if there is a  subset  of $X$, $A\in$ $\bsigma_{\alpha}^{0}$ such that 
$B=X\setminus A$.
	\item We write $\bdelta_{\alpha}^{0}=\bsigma_{\alpha}^{0}\cap\bpi_{\alpha}^{0}$.
\end{enumerate}
\end{definition}
The second and fourth items imply that the collection defined this way is closed under complements and the third item says, in particular, that $\bsigma^0_\alpha$ sets are constructed by taking arbitrary  countable unions of sets of lower complexity. 
	Thus $\bsigma^0_\alpha$ sets are closed under countable unions and $\bpi^0_\alpha$ sets are closed under countable intersections.
But none of the $\bsigma^0_\alpha$'s are closed under countable intersections, nor are the $\bpi^0_\alpha$'s closed 
under unions. Thus these definitions organize the Borel sets in a hierarchy of length $\omega_1$. Because 
$\omega_1$ has uncountable cofinality, 
$\bigcup_{\alpha<\omega_1}\bsigma^0_{\alpha}=	\bigcup_{\alpha<\omega_1}\bpi^0_{\alpha}$ forms a $\sigma$-algebra.

 Here are some facts:
	\begin{fact*} Let $X$ be an uncountable Polish space. Then for all countable $\alpha$:
	\begin{itemize}
	\item  $\bsigma^0_{\alpha+1}\ne\bpi^0_{\alpha+1}$.
	\item Both $\bsigma^0_\alpha$ and $\bpi^0_\alpha$ are proper subsets of $\bdelta^0_{\alpha+1}$.
	\item In turn, each $\bdelta^0_\alpha$ is a proper subset of each $\bsigma^0_{\alpha+1}$ and each 
	$\bpi^0_{\alpha+1}$.
	\end{itemize}
	\end{fact*}

\begin{figure}[!h]
	\centering
	\includegraphics[height=.17\textheight]{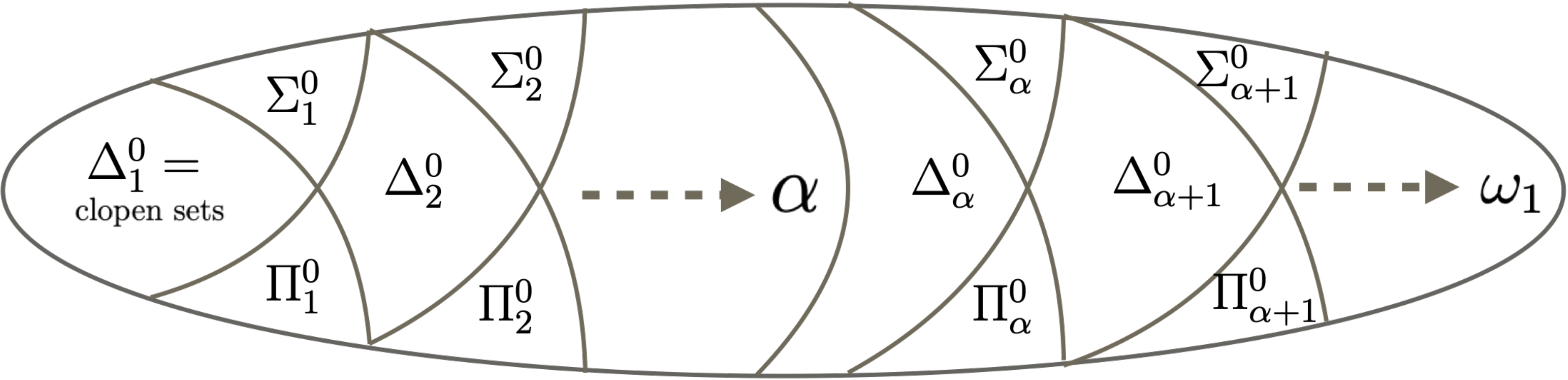}
	\caption{The construction of the Borel Sets.}
	\label{whale}
	\end{figure}

These facts are summarized in  figure \ref{whale}.

\paragraph{Analytic and co-Analytic sets}\hypertarget{aandcoa}{We now turn to the next type higher.  We define two collections of universally measurable sets using projections.}

Let $X$ and $Y$ be Polish spaces and $B\subseteq X\times Y$ be a Borel set.  The set 
	\begin{equation}
	A=\{x:\mbox{(for some }y) (x,y)\in B\}\label{analdef}
	\end{equation}
is called an \emph{analytic} set. It's complement
	\begin{equation}
	C=\{x:\mbox{(for all }y)(x, y)\notin B\}\label{coanaldef}
	\end{equation}
is called a \emph{co-analytic} set.
Note that we could get an equivalent definition of co-analytic by taking a Borel set $B^*$ and setting 
	\begin{equation}
	C=\{x:\mbox{(for all }y)(x, y)\in B^*\}.\label{coanalyaltdef}
	\end{equation}
To be sure the definition is clear--the analytic sets are all sets of the form of $A$ in equation \ref{analdef} with $B$ a Borel set.  The co-analytic sets are all sets of the form in equation \ref{coanaldef}, or equivalently all sets of the form in equation \ref{coanalyaltdef} with $B$ a Borel set. 
For purposes of showing that a set is analytic or  co-analytic it is useful to use logical notation: ``for some y" can be be written $\exists y$ and ``for all y" can be written $\forall y$. Informally analytic sets are ``$(\exists y)$Borel" and co-analytic sets are ``$(\forall y)$Borel." Because products of Polish spaces are Polish, a set that can be written ``$(\exists y)(\exists z)$Borel" is analytic and ``$(\forall y)(\forall z)$Borel" is co-analytic.

\paragraph{Notation} \hypertarget{bfpisigma}{The ``$\bsigma^1_1$ sets" are the analytic sets and the ``$\bpi^1_1$-sets" are the co-analytic sets.}

\begin{fact*} Every perfect Polish space contains a non-Borel analytic set. Moreover, 
the analytic sets are closed under countable intersections and unions. Hence the co-analytic sets are also closed under unions and intersections.

\end{fact*}
\noindent In many ways the analytic and co-analytic sets function similarly to the open and closed sets.

\medskip

\hypertarget{suslinthm}{Here is a remarkable theorem about the relationship between analytic and co-analytic sets due to Lusin with a very important corollary due to Suslin.(See \cite{lusin}, or section 3.5 of \cite{naivedst}. Reference  \cite{monks} gives a popular history account.)}
	\begin{theorem}(Lusin's Separation Theorem) 	Suppose that  $X$ is a Polish space and $A, B$ are disjoint analytic subsets of $X$.
	Then there is a Borel set $C$ with $A\subset C$ and $C\cap B=\emptyset$. 
	(So $C$ \emph{separates} $A$ from $B$.)
	\end{theorem}
\hypertarget{suslintheorem}{The next corollary is immediate.}
	\begin{corollary}\label{afterall} (Suslin's Theorem)
	{Let $A\subseteq X$.  If both $A$ and $X\setminus A$ are analytic, then $A$ is Borel.}
	\end{corollary}

\begin{example}\label{nevertoolate} For Borel group actions, the orbit equivalence relation is always analytic. If $G$ is acting on $X$ then 
\[B=\{(x_0, x_1, g):gx_0=x_1\}\]
is a Borel set.  The pair $(x_0, x_1)$ lies in the same $G$-orbit if and only if $(x_0, x_1)\in A$ where 
\[A=\{(x_0, x_1):(\mbox{for some }g)(x_0, x_1, g)\in B\}.\]
Thus being in the same orbit is analytic.
\medskip

\noindent Similarly if  $\mcc$ and $\mathcal D$ are Polish spaces of measure-preserving transformations and $E\subseteq \mathcal D$ is a Borel set then 
\[A=\{T\in \mcc:(\mbox{for some measure isomorphism $\phi$ and some }S\in  E)\ \phi\circ T=S\circ\phi\ (a.e.)\}\]
is analytic. 

From this one can see that the collection of measure-preserving transformations of $[0,1]$ that are isomorphic to a Bernoulli shift is an analytic set (details in \cite{descview}). 
These results are weaker than Feldman's paper where he shows that the collection of measure preserving transformations isomorphic to a Bernoulli shift is Borel (\cite{feldman}).  This is discussed in section  \ref{bs}.\end{example}

\begin{example}\label{distaliscoanalytic} Let $X$ be a compact metric space. Then the collection of distal homeomorphisms of $X$ is a {co-analytic set}
\end{example}
\pf Because \[\{(T, x, y):x\ne y \mbox{ and (for some $q$)\! (for all }n), d(T^n(x),T^n(y))>1/q\}.\]
is Borel, 
\[\{T: (\mbox{for all }x\ne y)(\mbox{for some $q$ for all }n, d(T^n(x),T^n(y))>1/q\}.\]
is in the form of equation \ref{coanalyaltdef}. It follows that the collection of distal homeomorphisms of $X$ is co-analytic.\qed

\paragraph{Reductions and hierarchies}  

	\begin{fact*} 
	Let $X, Y$ be  Polish spaces and $A\subseteq X$, $B\subseteq Y$. Let $f:X\to Y$ reduce $A$ to $B$.
		\begin{itemize}
		\item If $f$ is continuous and $B\in \bsigma^0_\alpha$ then $A\in \bsigma^0_\alpha$. 
		\item If $f$ is continuous and $B\in \bpi^0_\alpha$ then $A\in \bpi^0_\alpha$.
		\item If $f$ is Borel and $B$ is analytic, then $A$ is analytic.
		\item If $f$ is Borel and $B$ is co-analytic, then $A$ is co-analytic.
		\end{itemize}
	\end{fact*}
\paragraph{Orderings and pre-orderings}
Let $(I, \le_I)$ be a linear ordering of a set of cardinality less than or equal to the cardinality of the real numbers and $X$ a perfect Polish space.  Let $C\subseteq X$ and $\pi:C\to I$ be a surjective function.  Then $I$ is coded by the set of pairs 
$\le_\pi=\{(x,y):\pi(x)\le_I \pi(y)\}$ in the sense that the isomorphism type of $I$  can be recovered from 
$\le_\pi$. 
 Note that 
$\le_\pi\subseteq X\times X$, and hence can be studied as a subset of a Polish space, with no direct mention of $I$ itself.   

Of particular interest is the special case where $I$ is an ordinal $\alpha$ and $\le_I$ is the well-ordering of the ordinals less than $\alpha$. In this case $\pi$ itself can be recovered by $\le_\pi$.  When $I$ is a well-ordering, 
the relation $\le_\pi$ is called a \emph{pre-well-ordering} because it is a pre-ordering   and taking the quotient by the equivalence relation 
\[x\sim y\mbox{ if and only if } (x\le_\pi y \mbox{ and } y\le_\pi x)\]
 results in a well-ordering of the $\sim$ equivalence classes that is isomorphic to $\alpha$.  When $I$ is a well-ordering the  function $\pi$  is called a \emph{norm}.

\begin{example*}
\hypertarget{wo}{Let $WO\subseteq \{0, 1\}^{\nn\times \nn}$ be the set of characteristic functions of well-orderings of $\nn$. } Then $WO$ is a 
complete $\bpi^1_1$-set.
\end{example*}
\noindent The set $WO$ is frequently called the collection of \emph{codes} for well-orderings.

For notational convenience when $\pi:C\to \alpha$ is a norm, it can be extended to all of $X$ by adding a point, denoted $\infty$ to the range of $\pi$  and extending $\pi$
 by setting $\pi(x)=\infty$ for every element of $X\setminus C$. \smallskip

\hypertarget{pi11norm}{The following definition is difficult to digest, but the technicalities turn out to be essential.}

\begin{definition} \label{pi11normsdef} Let $X$ be a Polish space, $\alpha$ an ordinal, $A$  a $\bpi^1_1$-set and $\pi:X\to \alpha\cup \{\infty\}$ be such that $x\in A$ iff $\pi(x)<\infty$. Then $\pi$ is a $\bpi^1_1$-norm iff 
the relations $\le^*_\pi$ and $<^*_\pi$ are both $\bpi^1_1$ where:
\begin{enumerate}
\item $x\le^*_\pi y$ iff $x\in A$ and $(\pi(x)\le \pi(y))$,
\item $x<^*_\pi y$ iff $x\in A$ and $(\pi(x)<\pi(y))$.
\end{enumerate}
\end{definition}
\noindent A motivated reader might ask \emph{why $\bpi^1_1$}? The answer is that an analogous $\pi$ on a $\bsigma^1_1$ set can be shown not to exist.  
\smallskip

The next proposition explains, in part, the importance of $\bpi^1_1$-norms.

\begin{prop}
Suppose that $C\subseteq X$ is a $\bpi^1_1$-set, with $X$ Polish and  $\pi:C\to \omega_1$ is a 
$\bpi^1_1$-norm. Then $B\subseteq C$ is Borel implies that there is a countable ordinal $\alpha$ such that $\pi[B]\subseteq \alpha$.
\end{prop}
\hypertarget{overspill}{A corollary of this proposition is that if $C$ is a $\bpi^1_1$ subset of a Polish space $X$, 
$\pi:X\to \omega_1\cup \{\infty\}$ is a $\bpi^1_1$-norm on $C$, and $A\subseteq X$ is
a Borel set with $\pi$ mapping a subset of $A$ to a cofinal subset of $\omega_1$, then $A\not\subseteq C$.  This phenomenon is called \emph{overspill}.}

\paragraph{The reduction property of co-analytic sets}
The reduction property for $\bpi^1_1$-sets is related to $\bpi^1_1$-norms. It can be used to establish Lusin's theorem.
	\begin{theorem}(The reduction property for $\bpi^1_1$-sets.) Let $X$ be a perfect Polish space. 
	Suppose that $A, B\subseteq X$ are $\bpi^1_1$-sets.  Then there are $A_0\subseteq A, B_0\subseteq B$ with:
		\begin{itemize}
		\item $A_0, B_0\in \bpi^1_1$
		\item $A_0\cap B_0=\emptyset$
		\item $A_0\cup B_0=A\cup B$. 
		\end{itemize}
	\end{theorem}

\section{Open Problems}
The author has made no attempt to make a complete list of open problems--there are dozens. Nor have they attempted correct attributions.  These are not even the most important problems, although some, such as Open Problem \ref{OP5}, are very classical and well studied.  However they ask questions raised by this survey that may not appear in other places. 
\begin{open}\label{OP1}
 Are there complete numerical invariants for orientation preserving diffeomorphisms of the circle with irrational rotation number up to conjugation by orientation preserving  diffeomorphisms? 
 \end{open}
 
 \bfni{Comment:} After the first draft of this paper was written, Kunde solved this negatively  by showing that the equivalence relation $E_0$ is reducible to diffeomorphisms with Liouvillean rotation number under the equivalence relation of conjugation by diffeomorphisms. 
 
 \begin{open}\label{OP2}
 Is there an algorithm for determining whether two subshifts of finite type are conjugate by homeomorphisms?
 \end{open}
 
 \begin{open}\label{OP3} Is the analogue of Theorem \ref{BorelBorel} true for groups besides $\poZ$.  For example: 
 let $X$ be the space of Lebesgue measure preserving actions of $F_2$, the free group  on two generators, on $[0, 1]$. Is the collection of actions measure isomorphic to Bernoulli $F_2$-actions a Borel set? 
  \end{open}
 
 \begin{open}\label{OP7}
Is measure isomorphism restricted to $\mck$-automorphisms a Borel equivalence relation? Can it be reduced to an $S_\infty$-action? 
\smallskip

\noindent What about measure isomorphism for weakly mixing transformations?  Or other classes of transformations?
\end{open}

\begin{open}\label{snuck in1} 
Consider various classes of ergodic probability measure preserving transformations such as weakly mixing, mixing, 0-entropy, and $\mck$-automorphisms.  Is there a\ \ $\preceq^2_\mcb$ relationship between any of these classes, or between one of these classes and the class of arbitrary ergodic probability measure preserving transformations.
\end{open}

\begin{open}\label{snuck in2} Is isomorphism for ergodic measure preserving transformations on $\sigma$-finite measure spaces Borel?
\smallskip

\noindent What is the $\preceq^2_\mcb$ relationship between isomorphisms of ergodic probability measure preserving transformations and isomorphism of ergodic $\sigma$-finite measure preserving transformations?
\end{open}

\begin{open}\label{OP4.5} Do the structure theorems for distality give  Borel criteria for isomorphism? More precisely:
\begin{enumerate}
\item[A.] Let $X$ be the collection of minimal homeomorphisms on compact metric spaces and $\mathcal D$ be the collection of topologically distal transformations.  Is there a Borel set $B\subseteq X\times X$ such that $B\cap (\mathcal D\times \mathcal D)$ is the relation of topological conjugacy?
\item[B.] Let $\mathcal{ MD}$ be the collection of ergodic measure distal transformations.  Is there a Borel set $B\subseteq EMPT\times EMPT$ such that $B\cap (\mathcal{ MD}\times\mathcal{MD})$ is the relation of measure conjugacy?
\end{enumerate}

\end{open}
 
  \begin{open}\label{OP5}
 Suppose that $T:X\to X$ is an ergodic finite entropy transformation. Is there a compact manifold $M$ with a smooth volume element $\nu$ and a measure-preserving diffeomorphism $S$ such that $\xbmt$ is measure isomorphic to  $(M, \mcc, \nu, S)$?
 \end{open}
 
 \begin{open}\label{OP6} Let $\mathcal C$ be the collection of diffeomorphisms of a compact manifold $M$ that are topologically conjugate to a structurally stable diffeomorphism.  Let $E$ be the relation of topological conjugacy restricted to $\mathcal C$.  Is $E$ Borel reducible to $=$?  In other words does $E$ have Borel computable, complete numerical invariants? 
 \smallskip
 
\noindent This can also be asked for any standard class of structurally stable transformations, such as the Morse-Smale transformations. Is the conjugacy relation Borel?
\smallskip
 
\noindent Similarly in classes where there are existing invariants (such as the \emph{schemes} that are used in dimension 3 in \cite{bonatti}) one can ask whether those invariants are themselves Borel. 
  \end{open}

\begin{open}
Is the relation of \emph{topological conjugacy} Borel when restricted to Axiom A diffeomorphisms of a compact surface?  Can $E_0$ be embedded into topological conjugacy for Axiom A transformations?
\end{open}

\begin{open}\label{OP8} Let $M$ be a smooth manifold with a smooth volume element $\nu$.  Let ${\mathcal{SE}}$ be the collection of smooth, ergodic, $\nu$-measure-preserving transformations of $M$. Is the measure isomorphism relation on ${\mathcal{SE}}$ bi-reducible with the maximal $S_\infty$-action?
\end{open}

\begin{open}\label{OP9}(Sabok's Conjecture) Is the measure isomorphism relation on ergodic measure-preserving transformations of the unit interval $\preceq^2_\mcb$ bi-reducible with the maximal Polish group action?
\end{open}

\begin{open}\label{OP10} Where does the Kakutani equivalence relation on ergodic measure-preserving transformations sit among the analytic equivalence relations in $\preceq^2_\mcb$? Is it reducible to an $S_\infty$-action?  It is known to be a complete analytic set, hence not Borel (see Theorem \ref{kakisbad}). 
\smallskip

\noindent In particular what is the $\preceq^2_\mcb$ relationship between the equivalence relation of Kakutani equivalence and isomorphism  between ergodic probability measure preserving transformations? Between Kakutani equivalence and  isomorphism of ergodic measure preserving transformations on $\sigma$-finite measure spaces?
\end{open}

\begin{open}\label{OP11} Does $E_0$ reduce to the collection of topologically minimal diffeomorphisms of the 2-torus with the relation of  topological conjugacy? What about topologically transitive diffeomorphisms?
\end{open}

\begin{open}\label{OP12.5} 
Is topological conjugacy of minimal homeomorphisms of the Cantor set bi-reducible  with the maximal $S_\infty$-action?  What about topologically transitive homeomorphisms? 
\end{open}

\begin{open}\label{OP14}
Suppose that $F$ is a Borel equivalence relation on a Polish space that is not Borel reducible to the orbit equivalence relation of a  Polish group action. Is 
$E_1\preceq^2_\mcb F$?
\end{open}

\section{Acknowledgements}
The author received help from many people in writing and editing this paper.  He particularly wants to acknowledge Filippo Calderoni, Marlies Gerber,  Anton Gorodetski,   Philipp Kunde, Andrew Marks, Ronnie Pavlov. Alexander Kechris made invaluable corrections on an early text and correspondence with Su Gao was essential in the completion of the text.   As always, Benjamin Weiss provided important suggestions and comments and provided very helpful references.
\medskip 

\paragraph{Expository Books and articles on analytic equivalence relations}  
\begin{enumerate} 
\item The book by Gao \cite{Gaobook}, 
\item The books by Kechris \cite{kechris}, Kechris-Miller  \cite{kechmiller} and Becker-Kechris \cite{BK},
\item The paper by Friedman and Stanley \cite{Fried},
\item Dave Marker's website \cite{marker},
\end{enumerate}
\noindent The author would like to acknowledge partial support from NSF grant DMS-2100367 

\bibliography{cldscit.bib}

\providecommand{\bysame}{\leavevmode\hbox to3em{\hrulefill}\thinspace}
\providecommand{\MR}{\relax\ifhmode\unskip\space\fi MR }
\providecommand{\MRhref}[2]{%
  \href{http://www.ams.org/mathscinet-getitem?mr=#1}{#2}
}
\providecommand{\href}[2]{#2}
\begin{thebibliography}{10}

\bibitem{Aaronson}
Jon Aaronson, \emph{An introduction to infinite ergodic theory}, Mathematical
  Surveys and Monographs, vol.~50, American Mathematical Society, 1997.

\bibitem{Barreira}
Luis Barreira and Claudia Valls, \emph{Dynamical systems, an introduction},
  Universitext, Springer, London, 2013.

\bibitem{BK}
Howard Becker and Alexander~S. Kechris, \emph{The descriptive set theory of
  polish group actions}, London Mathematical Society Lecture Note Series,
  Cambridge University Press, 1996.

\bibitem{BF1}
Ferenc Beleznay and Matthew Foreman, \emph{The collection of distal flows is
  not {B}orel}, American Journal of Mathematics \textbf{117} (1995), no.~1,
  203--239.

\bibitem{FB2}
Ferenc Beleznay and Matthew Foreman, \emph{The complexity of the collection of
  measure-distal transformations}, Ergodic Theory Dynam. Systems \textbf{16}
  (1996), no.~5, 929--962.
%

\bibitem{bonatti}
C.~Bonatti, V.~Grines, and O.~Pochinka, \emph{Topological classification of
  {M}orse-{S}male diffeomorphisms on 3-manifolds}, Duke Math. J. \textbf{168}
  (2019), no.~13, 2507--2558.

\bibitem{brayamp}
Mark Braverman and Michael Yampolsky, \emph{Computability of {J}ulia sets},
  Algorithms and Computation in Mathematics, vol.~23, Springer-Verlag, Berlin,
  2009.

\bibitem{BY}
\bysame, \emph{Computability of {J}ulia sets}, Algorithms and Computation in
  Mathematics, vol.~23, Springer-Verlag, Berlin, 2009.

\bibitem{CG}
Riccardo Camerlo and Su~Gao, \emph{The completeness of the isomorphism relation
  for countable {B}oolean algebras}, Trans. Amer. Math. Soc. \textbf{353}
  (2001), no.~2, 491--518.

\bibitem{chacon}
R.~V. Chacon, \emph{Weakly mixing transformations which are not strongly
  mixing}, Proc. Amer. Math. Soc. \textbf{22} (1969), 559--562.

\bibitem{clemens}
John~D. Clemens, \emph{Isomorphism of subshifts is a universal countable
  {B}orel equivalence relation}, Israel J. Math. \textbf{170} (2009), 113--123.

\bibitem{DG}
Longyun Ding and Su~Gao, \emph{Is there a spectral theory for all bounded
  linear operators?}, Notices Amer. Math. Soc. \textbf{61} (2014), no.~7,
  730--735.

\bibitem{DJKhyperfinite}
R.~Dougherty, S.~Jackson, and A.~S. Kechris, \emph{The structure of hyperfinite
  {B}orel equivalence relations}, Trans. Amer. Math. Soc. \textbf{341} (1994),
  no.~1, 193--225.

\bibitem{effros}
Edward~G. Effros, \emph{Transformation groups and {$C^{\ast} $}-algebras}, Ann.
  of Math. (2) \textbf{81} (1965), 38--55.

\bibitem{feldman}
Jacob Feldman, \emph{Borel structures and invariants for measurable
  transformations}, Proc. Amer. Math. Soc. \textbf{46} (1974), 383--394.

\bibitem{FeldmanMoore}
Jacob Feldman and Calvin~C. Moore, \emph{Ergodic equivalence relations,
  cohomology, and von {N}eumann algebras. {I}}, Trans. Amer. Math. Soc.
  \textbf{234} (1977), no.~2, 289--324.

\bibitem{FLR}
Valentin Ferenczi, Alain Louveau, and Christian Rosendal, \emph{The complexity
  of classifying separable {B}anach spaces up to isomorphism}, J. Lond. Math.
  Soc. (2) \textbf{79} (2009), no.~2, 323--345.

\bibitem{descview}
Matthew Foreman, \emph{A descriptive view of ergodic theory}, Descriptive set
  theory and dynamical systems ({M}arseille-{L}uminy, 1996), London Math. Soc.
  Lecture Note Ser., vol. 277, Cambridge Univ. Press, Cambridge, 2000,
  pp.~87--171.

\bibitem{naivedst}
\bysame, \emph{Naive descriptive set theory, \\
  \url{arXiv:2110.08881}},  (2010), 1--53.

\bibitem{whatis}
\bysame, \emph{What is a {B}orel reduction?}, Notices Amer. Math. Soc.
  \textbf{65} (2018), no.~10, 1263--1268.

\bibitem{FoGo}
Matthew Foreman and Anton Gorodetski, \emph{Anti-classification results for
  smooth dynamical systems,   \url{arXiv:2206.09322}}, (2022), 1--56.
  
\bibitem{FRW}
Matthew Foreman, Daniel~J. Rudolph, and Benjamin Weiss, \emph{The conjugacy
  problem in ergodic theory}, Ann. of Math. (2) \textbf{173} (2011), no.~3,
  1529--1586.

\bibitem{ananti}
Matthew Foreman and Benjamin Weiss, \emph{An anti-classification theorem for
  ergodic measure preserving transformations}, J. Eur. Math. Soc. (JEMS)
  \textbf{6} (2004), no.~3, 277--292.

\bibitem{globstruct}
\bysame, \emph{From odometers to circular systems: a global structure theorem},
  J. Mod. Dyn. \textbf{15} (2019), 345--423.

\bibitem{circsys}
\bysame, \emph{A symbolic representation for {A}nosov-{K}atok systems}, J.
  Anal. Math. \textbf{137} (2019), no.~2, 603--661.

\bibitem{noclass}
\bysame, \emph{Measure preserving diffeomorphisms of the torus are
  unclassifiable}, J. Eur. Math. Soc. (JEMS) (To Appear), 1 -- 80.

\bibitem{Fried}
Harvey Friedman and Lee Stanley, \emph{A {B}orel reducibility theory for
  classes of countable structures}, J. Symbolic Logic \textbf{54} (1989),
  no.~3, 894--914.

\bibitem{topodistal}
H.~Furstenberg, \emph{The structure of distal flows}, Amer. J. Math.
  \textbf{85} (1963), 477--515.

\bibitem{FUbook}
H. ~Furstenberg, \emph{Recurrence in ergodic theory and combinatorial number
  theory}, Porter Lectures, Princeton University Press, 1981.

\bibitem{GaoLecture}
Su~Gao, \emph{{D}ynamical {S}ystems and {C}ountable {S}tructures, {L}ecture
  {N}otes:\\
  \url{https://www.birs.ca/workshops/2022/22w5134/files/Su\%20Gao/banff\%20lecture\%20gao\%20nopause.pdf}}.

\bibitem{Gaobook}
\bysame, \emph{Invariant descriptive set theory}, A Series of Monographs and
  Textbooks, vol. 293, Taylor and Francis, 2009.

\bibitem{GP}
Su~Gao and Vladimir Pestov, \emph{On a universality property of some abelian
  {P}olish groups}, Fund. Math. \textbf{179} (2003), no.~1, 1--15.

\bibitem{gerberkunde}
Marlies Gerber and Philipp Kunde, \emph{Anti-classification results for the
  {K}akutani equivalence relation, \url{{https://arxiv.org/abs/2109.06086v1}}},
   (2021), 1 -- 72.

\bibitem{glimm}
James Glimm, \emph{Locally compact transformation groups}, TAMS \textbf{101}
  (1961), 124--138.

\bibitem{monks}
Loren Graham and Jean-Michel Kantor, \emph{Naming infinity}, Belknap Press,
  2009.

\bibitem{halmos1944}
Paul~R. Halmos, \emph{In general a measure preserving transformation is
  mixing}, Ann. of Math. (2) \textbf{45} (1944), 786--792.

\bibitem{halmosbook}
\bysame, \emph{Lectures on ergodic theory}, Chelsea Publishing Company,
  New York, 1956.

\bibitem{halmosnaive}
\bysame, \emph{Naive set theory}, Undergraduate Texts in Mathematics,
  Springer-Verlag, 1974.
  
  \bibitem{HvN} Paul~R. Halmos and John von Neumann, \emph{Operator methods in 
  classical mechanics II}, Ann. of Math. (2) \textbf{43} (1942), 332-350.

\bibitem{HKL}
L.~A. Harrington, A.~S. Kechris, and A.~Louveau, \emph{A {G}limm-{E}ffros
  dichotomy for {B}orel equivalence relations}, J. Amer. Math. Soc. \textbf{3}
  (1990), no.~4, 903--928.

\bibitem{hassel}
Boris Hasselblatt, \emph{Problems in dynamical systems and related topics},
  Dynamics, ergodic theory, and geometry, Math. Sci. Res. Inst. Publ., vol.~54,
  Cambridge Univ. Press, Cambridge, 2007, pp.~273--324.

\bibitem{Hjorth}
G.~Hjorth, \emph{On invariants for measure preserving transformations}, Fund.
  Math. \textbf{169} (2001), no.~1, 51--84.

\bibitem{hjcongn}
\bysame, \emph{Around nonclassifiability for countable torsion free abelian
  groups}, Abelian Groups and Modules: International Conference in Dublin,
  August 10--14, 1998 (Paul~C. Eklof and R{\"u}diger G{\"o}bel, eds.),
  Birkh{\"a}user Basel, 1999, pp.~269--292.

\bibitem{Hj}
\bysame, \emph{Classification and orbit equivalence relations},
  Mathematical Surveys and Monographs, vol.~75, American Mathematical Society,
  2000.

\bibitem{kechcount}
S.~Jackson, A.~S. Kechris, and A.~Louveau, \emph{Countable {B}orel equivalence
  relations}, J. Math. Log. \textbf{2} (2002), no.~1, 1--80.



\bibitem{Jech}
Thomas~J. Jech, \emph{Set theory: The third millenium edition}, Springer, 2003.

\bibitem{kanovei}
Vladimir Kanovei, \emph{Borel equivalence relations: structure and
  classification}, University Lecture Series, vol.~14, American Mathematical
  Society, 2008.

\bibitem{HasKat}
Anatole Katok and Boris Hasselblatt, \emph{Introduction to the modern theory of
  dynamical systems}, Encyclopedia of Mathematics and its applications,
  vol.~54, Cambridge University Press, 1995.
  
  \bibitem{kechris}
Alexander Kechris, \emph{Classical descriptive set theory}, 1 ed., Graduate
  Texts in Mathematics, Springer, New York, NY, 1995.
  
\bibitem{kechrisglobal} 
\bysame, \emph{Global Aspects of Ergodic Group Actions}, Mathematical Surveys and Monographs, vol. 160; Am. Math. Soc., 2010.



\bibitem{kechlouv} ~A. ~S. ~Kechris and A. ~Louveau, \emph{
The Classification of Hypersmooth Borel Equivalence Relations}, :
J. ~American Mathematical Society, Vol. 10, No. 1 (Jan., 1997), pp. 215-242


\bibitem{kechmiller}
Alexander Kechris and Ben~D. Miller, \emph{Topics in orbit equivalence},
  Lecture Notes in Mathematics, Springer Berlin Heidelberg, 2004.

\bibitem{KeSo}
A.~S. Kechris and N.~E. Sofronidis, \emph{A strong generic ergodicity property
  of unitary and self-adjoint operators}, Ergodic Theory Dynam. Systems
  \textbf{21} (2001), no.~5, 1459--1479.

\bibitem{kechdrob}
Alexander~S. Kechris and Robin~D. Tucker-Drob, \emph{The complexity of
  classification problems in ergodic theory}, Appalachian set theory
  2006--2012, London Math. Soc. Lecture Note Ser., vol. 406, Cambridge Univ.
  Press, Cambridge, 2013, pp.~265--299.

\bibitem{Khintchine1949}
A.~I. Khinchin, \emph{Mathematical {F}oundations of {S}tatistical {M}echanics},
  Dover Publications, Inc., New York, N.Y., 1949, Translated by G. Gamow.

\bibitem{king}
Jonathan King, \emph{The commutant is the weak closure of the powers, for
  rank-{$1$} transformations}, Ergodic Theory Dynam. Systems \textbf{6} (1986),
  no.~3, 363--384.

\bibitem{krieg1}
Wolfgang Krieger, \emph{On entropy and generators of measure-preserving
  transformations}, Trans. Amer. Math. Soc. \textbf{149} (1970), 453--464.

\bibitem{krieg2}
\bysame, \emph{Erratum to: ``{O}n entropy and generators of measure-preserving
  transformations''}, Trans. Amer. Math. Soc. \textbf{168} (1972), 519.

\bibitem{krieg3}
\bysame, \emph{On unique ergodicity}, Proceedings of the {S}ixth {B}erkeley
  {S}ymposium on {M}athematical {S}tatistics and {P}robability ({U}niv.
  {C}alifornia, {B}erkeley, {C}alif., 1970/1971), {V}ol. {II}: {P}robability
  theory, 1972, pp.~327--346.

\bibitem{Levy}
Azriel Levy, \emph{Basic set theory}, Dover Publications, Inc., Mineola, NY,
  2002, Reprint of the 1979 original [Springer, Berlin; MR0533962 (80k:04001)].

\bibitem{dougmarc}
Doug Lind and Brian Marcus, \emph{An introduction to symbolic dynamics and
  coding}, Cambridge University Press, 1995.

\bibitem{loy}
Joram Lindenstrass, Gunnar Olsen, and Y.~Sternfeld, \emph{The {P}oulsen
  simplex}, Annales de l'institute Fourier \textbf{28} (1978), no.~1, 91--114.

\bibitem{lusin}
Nicolas Lusin, \emph{Sur les classes des constituantes des compl\'{e}mentaires
  analytiques}, Ann. Scuola Norm. Super. Pisa Cl. Sci. (2) \textbf{2} (1933),
  no.~3, 269--282.

\bibitem{marker}
David Marker, \emph{Descriptive set theory, \\
  \url{http://homepages.math.uic.edu/~marker/math512/dst.pdf}},  (2002), 1--56.

\bibitem{maruyama}
Gisir{o} Maruyama, \emph{The harmonic analysis of stationary stochastic
  processes}, Mem. Fac. Sci. Ky\={u}sy\={u} Univ. A \textbf{4} (1949), 45--106.

\bibitem{Mekler} 
A.~ H. Mekler, \emph{Stability of nilpotent groups of class 2 and prime exponent}, J. {S}ymbolic {L}ogic, {\bf 46} (1981), no. 4, 781--788.


\bibitem{Mirza}
Maryam Mirzakhani and Tony Feng, \emph{Introduction to ergodic theory, \\
  \url{https://www.mit.edu/~fengt/ergodic_theory.pdf}},  (2014), 1 -- 56.

\bibitem{Mitin}
A.~V. Mitin, \emph{Undecidability of the elementary theory of groups of
  measure-preserving transformations}, Mat. Zametki \textbf{63} (1998), no.~3,
  414--420.

\bibitem{moscho}
Yiannis~N. Moschovakis, \emph{Descriptive set theory: Second edition}, 2 ed.,
  Mathematical Surveys and Monographs, no. 155, American Mathematical Society,
  2009.

\bibitem{newhouse}
Sheldon~E. Newhouse, \emph{Nondensity of axiom {${\rm A}({\rm a})$} on
  {$S^{2}$}}, Global {A}nalysis ({P}roc. {S}ympos. {P}ure {M}ath., {V}ol.
  {XIV}, {B}erkeley, {C}alif., 1968), Amer. Math. Soc., Providence, R.I., 1970,
  pp.~191--202.

\bibitem{newhouse2}
\bysame, \emph{Diffeomorphisms with infinitely many sinks}, Topology
  \textbf{13} (1974), 9--18.

\bibitem{ornstein1}
Donald Ornstein, \emph{Bernoulli shifts with the same entropy are isomorphic},
  Advances in Math. \textbf{4} (1970), 337--352.

\bibitem{ornstein2}
\bysame, \emph{The isomorphism problem for measure-preserving transformations},
  Functional {A}nalysis ({P}roc. {S}ympos., {M}onterey, {C}alif., 1969),
  Academic Press, New York, 1970, pp.~71--74. \MR{0262463}

\bibitem{ornstein3}
\bysame,  \emph{Ergodic theory, randomness, and dynamical systems},
  Yale Mathematical Monographs, No. 5, Yale University Press, New Haven,
  Conn.-London, 1974, James K. Whittemore Lectures in Mathematics given at Yale
  University.

\bibitem{ornstein4}
\bysame, \emph{Some {Open} {Problems} in {Ergodic} {Theory}}, Publications
  math\'ematiques et informatique de Rennes (1975), no.~S4.

\bibitem{ornsteinshields}
Donald~S. Ornstein and Paul~C. Shields, \emph{An uncountable family of
  {$K$}-automorphisms}, Advances in Math. \textbf{10} (1973), 63--88.

\bibitem{FDimpsBern}
Donald~S. Ornstein and Benjamin Weiss, \emph{Finitely determined implies very
  weak {B}ernoulli}, Israel J. Math. \textbf{17} (1974), 94--104.

\bibitem{shark}
A.~A. Oshemkov and V.~V. Sharko, \emph{On the classification of {M}orse-{S}male
  flows on two-dimensional manifolds}, Mat. Sb. \textbf{189} (1998), no.~8,
  93--140.

\bibitem{peixoto}
M.~M. Peixoto, \emph{On the classification of flows on {$2$}-manifolds},
  Dynamical systems ({P}roc. {S}ympos., {U}niv. {B}ahia, {S}alvador, 1971),
  1973, pp.~389--419.

\bibitem{peterson}
Karl~E. Petersen, \emph{Ergodic theory}, Cambridge Studies in Advanced
  Mathematics,, Cambridge University Press, 1983.

\bibitem{rokhlin1948}
V.~Rohlin, \emph{A ``general" a measure-preserving transformation is not
  mixing}, Doklady Akad. Nauk SSSR \textbf{60} (1948), 349--351.

\bibitem{therokh}
\bysame, \emph{Lectures on the entropy theory of transformations with
  invariant measure}, Uspehi Mat. Nauk \textbf{22} (1967), no.~5 (137), 3--56.

\bibitem{Rosendal}
Christian Rosendal, \emph{Cofinal families of {B}orel equivalence relations and
  quasiorders}, J. Symbolic Logic \textbf{70} (2005), no.~4, 1325--1340.

\bibitem{Ryzhikov}
V.~V. Ryzhikov, \emph{Representation of transformations preserving the
  {L}ebesgue measure, in the form of a product of periodic transformations},
  Mat. Zametki \textbf{38} (1985), no.~6, 860--865, 957.

\bibitem{Sabok}
Marcin Sabok, \emph{Completeness of the isomorphism problem for separable {$\rm
  C^\ast$}-algebras}, Invent. Math. \textbf{204} (2016), no.~3, 833--868.

\bibitem{Smale1963}
S.~Smale, \emph{Dynamical systems and the topological conjugacy problem for
  diffeomorphisms}, Proc. {I}nternat. {C}ongr. {M}athematicians ({S}tockholm,
  1962), Inst. Mittag-Leffler, Djursholm, 1963, pp.~49--496.

\bibitem{Smale1967}
\bysame, \emph{Differentiable dynamical systems}, Bull. Amer. Math. Soc.
  \textbf{73} (1967), 747--817.

\bibitem{soare}
Robert Soare, \emph{Recursively enumerable sets and degrees.}, Perspectives in
  Mathematical Logic, Springer Berlin Heidelberg, 1987.

\bibitem{thomasp-local}
Simon Thomas, \emph{The classification problem for torsion-free abelian groups
  of finite rank}, J. Amer. Math. Soc. \textbf{16} (2003), no.~1, 233--258.

\bibitem{thomasclsasprob}
\bysame, \emph{The classification problem for {$S$}-local torsion-free abelian
  groups of finite rank}, Adv. Math. \textbf{226} (2011), no.~4, 3699--3723.

\bibitem{TVshow}
Simon Thomas and Boban Velickovic, \emph{On the complexity of the isomorphism
  relation for finitely generated groups}, J. Algebra \textbf{217} (1999),
  no.~1, 352--373.

\bibitem{Thouvenot}
Jean-Paul Thouvenot, \emph{Entropy, isomorphism and equivalence in ergodic
  theory}, Handbook of dynamical systems, {V}ol. 1{A}, North-Holland,
  Amsterdam, 2002, pp.~205--238.

\bibitem{US}
V.~V. Uspenski\u{\i}, \emph{A universal topological group with a countable
  basis}, Funktsional. Anal. i Prilozhen. \textbf{20} (1986), no.~2, 86--87.

\bibitem{Walters}
Peter Walters, \emph{An introduction to ergodic theory}, Graduate Texts in
  Mathematics, Springer New York, 1982.

\bibitem{Weihrauch}
Klaus Weihrauch, \emph{Computable analysis}, Texts in Theoretical Computer
  Science. An EATCS Series, Springer-Verlag, Berlin, 2000, An introduction.

\bibitem{wikipedia:HamSys}
Wikipedia, \emph{{H}amiltonian system, \\
  \url{https://en.wikipedia.org/wiki/Hamiltonian_system}}.

\bibitem{wikipedia:rotation}
\bysame, \emph{{R}otation number, \\
  \url{https://en.wikipedia.org/wiki/Rotation_number}}.

\bibitem{Zim}
Robert~J. Zimmer, \emph{Ergodic actions with generalized discrete spectrum},
  Illinois J. Math. \textbf{20} (1976), no.~4, 555--588.

\bibitem{zippin}
L.~Zippin, \emph{Transformation groups}, Lectures in Topology; the University
  of Michigan Conference of 1940, Ann Arbor, University of Michigan Press,
  2013.

\end{thebibliography}
\bibliographystyle{amsplain}

 \end{document}